\documentclass[10pt,legno]{article}
\usepackage{amsmath, amssymb}
\usepackage{latexsym}
\usepackage{amssymb,epsfig,color}
\usepackage{graphicx,epsfig}
\usepackage[toc,page]{appendix}
\begin{document}
\newcommand{\qed}{\hfill \ensuremath{\square}}
\newtheorem{thm}{Theorem}[section]
\newtheorem{cor}[thm]{Corollary}
\newtheorem{lem}[thm]{Lemma}
\newtheorem{prop}[thm]{Proposition}
\newtheorem{defn}[thm]{Definition}
\newcommand{\proof}{\vspace{1ex}\noindent{\em Proof}. \ }
\newtheorem{pro}[thm]{proof}
\newtheorem{ide}[thm]{Idee}
\newtheorem{rem}[thm]{Remark}
\newtheorem{ex}[thm]{Example}
\bibliographystyle{plain}
\numberwithin{equation}{section}
%-------------------------------------------------------------------------
\numberwithin{equation}{section}
\newcounter{saveeqn}
\newcommand{\subeqn}{\setcounter{saveeqn}{\value{equation}}%
 \stepcounter{saveeqn}\setcounter{equation}{0}% =
\renewcommand{\theequation}{\mbox{\arabic{section}.\arabic{saveeqn}\alph{=
equation}}}} %\alph, \roman

\newcommand{\reseteqn}{\setcounter{equation}{\value{saveeqn}}%
\renewcommand{\theequation}{\arabic{section}.\arabic{equation}}}

% MATH =-------------------------------------------------------------------
\def\nm{\noalign{\medskip}}
\def\u{{\mathbf u}}
\def\e{{\mathbf e}}
\def\T{{\mathbf T}}
\def\aa{{\mathbf a}}
\def\A{{\mathbf A}}
\def\B{{\mathbf B}}
\def\v{{\mathbf v}}
\def\g{{\mathbf g}}
\def\I{{\mathbf I}}
\def\f{{\mathbf f}}
\def\n{{\mathbf n}}
\def\T{{\mathbf T}}
\def\N{{\mathbf N}}
\def\M{{\mathbf M}}
\def\w{{\mathbf w}}
\def\KK{{\mathbf K}}
\def\P{{\mathbf P}}
\def\QQ{{\mathbf Q}}
\def\f{{\mathbf f}}
\def\x{{\mathbf x}}
\def\F{{\mathbf F}}
\def\GG{{\mathbf G}}
\def\H{{\mathbf H}}
\def\U{{\mathbf U}}
\def\y{{\mathbf y}}
\def\t{{\mathbf t}}
\def\T{{\mathbf T}}
\def\bvarphi{\boldsymbol{\varphi}}
\def\bpsi{\boldsymbol{\psi}}
\def\bphi{\boldsymbol{\phi}}
\def\ta{\boldsymbol{\tau}}
\def\G{{\mathbf \Gamma}}
\def\KK{{\mathbf K}}
\def\PP{{\mathbf P}}
\def\QQ{{\mathbf Q}}
\def\LL{{\mathbf L}}
\def\etaa{{\boldsymbol \eta}}
\newcommand{\Om}{\Omega}
\newcommand{\om}{\omega}
\newcommand{\Real}{\mathbb{R}}
\newcommand{\nuu}{\tilde{{\nu}}}
\newcommand{\xe}{\tilde{x}}
\newcommand{\ye}{\tilde{{y}}}
\newcommand{\bohm}{{\partial}{\ohm}}
\newcommand{\la}{\langle}
\newcommand{\ra}{\rangle}
\newcommand{\ms}{\mathcal{S}_\ohm}
\newcommand{\mk}{\mathcal{K}_\ohm}
\newcommand{\mks}{\mathcal{K}_\ohm ^{\ast}}
\newcommand{\grad}{\bigtriangledown}
\newcommand{\ds}{\displaystyle}
\newcommand{\pf}{\medskip \noindent {\sl Proof}. ~ }
\newcommand{\p}{\partial}
\renewcommand{\a}{\alpha}
\newcommand{\z}{\zeta}
\newcommand\q{\quad}
\newcommand{\pd}[2]{\frac {\p #1}{\p #2}}
\newcommand{\pdl}[2]{\frac {\p^2 #1}{\p #2}}
\newcommand{\dbar}{\overline \p}
\newcommand{\eqnref}[1]{(\ref {#1})}
\newcommand{\na}{\nabla}
\newcommand{\ep}{\epsilon}
\newcommand{\vp}{\varphi}
\newcommand{\fo}{\forall}
\newcommand{\Scal}{\mathcal{S}}
\newcommand{\BScal}{\boldsymbol{\mathcal{S}}}
\newcommand{\BDcal}{\boldsymbol{\mathcal{D}}}
\newcommand{\Dcal}{\mathcal{D}}
\newcommand{\Kcal}{\mathcal{K}}
\newcommand{\BKcal}{\boldsymbol{\mathcal{K}}}
\newcommand{\K}{\boldsymbol{\mathbb{K}}}
\newcommand{\Ecal}{\mathcal{E}}
\newcommand{\Ncal}{\mathcal{N}}
\newcommand{\BNcal}{\boldsymbol{\mathcal{N}}}
\newcommand{\Abar}{\overline A}
\newcommand{\Rcal}{\mathcal{R}}
\newcommand{\Lcal}{\mathcal {L}}
\newcommand{\Tcal}{\mathcal {T}}
\newcommand{\Gcal}{\mathcal {G}
}
\newcommand{\Cbar}{\overline C}
\newcommand{\Ebar}{\overline E}
\newcommand{\RR}{\mathbb{R}}
\newcommand{\CC}{\mathbb{C}}
\newcommand{\NN}{\mathbb{N}}
\newcommand{\Z}{\mathbb{Z}}

\title{ Small Perturbations of an Interface for Elastostatic Problems }

\date{}

\author{ Jihene Lagha \thanks{  Universit\'e de Tunis El Manar,
    Facult\'e des sciences de Tunis, UR11ES24 Optimisation,
    Mod\'elisation et Aide \`{a} la D\'ecision, 2092 Tunis,  Tunisie
    (lagha.jihene@yahoo.fr). } \and Faouzi Triki\thanks{
Laboratoire Jean Kuntzmann, Universit\'e
Grenoble-Alpes \& CNRS, 700 Avenue Centrale,
38401 Saint-Martin-d'H\`eres, France (faouzi.triki@imag.fr); FT is partly supported by
LabEx PERSYVAL-Lab (ANR-11-LABX- 0025-01). }
\and Habib Zribi \thanks{Department of Mathematics, College of Sciences, University of Hafr Al Batin, P.o. 1803, Hafr Al Batin
31991, Saudi Arabia (zribi.habib@yahoo.fr).}}

\maketitle

\begin{abstract}
 We consider the Lam\'e system for an elastic medium consisting of
an  inclusion embedded in a homogeneous background medium.
Based on  the field expansion method (FE) and layer potential techniques,
we rigorously  derived the asymptotic expansion of the perturbed
displacement field due to  small perturbations in the interface of
the inclusion.
We extend these techniques to determine
a relationship between traction-displacement measurements and the  shape of the object and  derive an asymptotic
expansion for the perturbation in the elastic moments tensors (EMTs) due to the presence of
small changes in the interface of the inclusion.
\end{abstract}

\noindent {\footnotesize {\bf Mathematics subject classification
(MSC2000):} 35B30, 35C20, 31B10}


\section{Introduction and statement of the main results}
Consider a homogeneous  isotropic elastic inclusion $D$
embedded in  the background region  $\RR^2$, which is
occupied by a homogeneous isotropic elastic material. The boundary $\p
D$ of the inclusion is assumed to be of class $\mathcal{C}^2$. In this case,  $\p D$ can  be
parametrized by a vector-valued function $t\rightarrow X(t)$,  that is,
$\p D:=\{x=X(t), t\in [a,b] \mbox{ with } a<b\}$, where $X$ is
a $\mathcal{C}^2$-function satisfying $|X'(t)|=1$ for all $t\in [a,b]$,
and $X(a) = X(b)$.

 Let  $(\lambda_0,\mu_0)$ denote the background Lam\'e constants, that
are the elastic parameters in the absence of any inclusions.
 Assume  that the Lam\'e constants in the inclusion $D$ are given by
 $(\lambda_1, \mu_1)$ where $(\lambda_1, \mu_1)\not= (\lambda_0,
 \mu_0)$. We further assume  that
$\mu_j >0,$  $\lambda_j +\mu_j>0$ for $ j=0,1,$
$(\lambda_0-\lambda_1)(\mu_0-\mu_1)\geq 0$. As in \cite{LY}, we needed
the last assumption  in order
to guarantee the well-posedeness of boundary
integral equation representation of the displacement field (see for
instance Theorem 2 in \cite{ES}).

Let $\mathbb{C}_0$ and $\mathbb{C}_1$ be the  elasticity tensors
 for $\RR^2\backslash \overline {D}$
and $D$, respectively, which are  given by
\begin{align*}
(\mathbb{{C}}_m)_{ijkl}=\lambda_m \delta_{ij}\delta_{kl}+
\mu_{m}(\delta_{ik}\delta_{jl}+\delta_{il}\delta_{jk}) \q
\mbox{for }i,j,k,l=1,2,\q m=0,1.
\end{align*}
There is another way of expressing the isotropic elastic tensor
 which will be
used later. Let  $\mathbb{I}$ be the identity  $4$-tensor and $\I$
be the identity
$2$-tensor (the $2 \times 2$ identity matrix). Then $\mathbb{{C}}_m$
can be rewritten as
\begin{align}\label{tensor-C}
\ds \mathbb{C}_m=\lambda_m \I\otimes\I+2 \mu_m \mathbb{I}, \q\q m=0,1.
\end{align}
 Then, the
elasticity tensor for $\RR^2$ in the presence of the inclusion $D$
is then given by
\begin{align*}
\ds \mathbb{C}=\mathbb{C}_0\chi_{\RR^2\backslash \overline{D}}+\mathbb{C}_1\chi_{D},
\end{align*}
where $\chi_{D}$ is the indicator function of $D$.\\

In this paper, we  consider the following transmission  problem
\begin{equation}\label{Main-Pb}
\left\{
  \begin{array}{lll}
   \ds \nabla \cdot \big(\mathbb{C} \widehat{\nabla}\u\big)=0 &
\mbox{ in } \RR^2,\\
    \nm \ds \u(x)-\H(x)=O(|x|^{-1})& \mbox{ as } |x|\rightarrow \infty,
  \end{array}
\right.
\end{equation}
where   $\H$ is a vector-valued function satisfying
$ \nabla \cdot \big(\mathbb{C}_0 \widehat{\nabla}\H\big)=0$ in  $\RR^2$, and
$\widehat{\nabla}\u=\frac{1}{2}\big(\nabla \u+(\nabla \u)^{T}\big)$ is
the
symmetric  strain tensor.
Here and throughout the paper $\M^T$
denotes the transpose of the matrix $\M$.\\

The elastostatic operator corresponding to the Lam\'e constants
$(\lambda_0, \mu_0)$ is defined by
\begin{align}\label{elastostatic-system}
\ds \mathcal{L}_{\lambda_0,\mu_0}\u:= \mu_0\Delta
\u+(\lambda_0+\mu_0) \nabla \nabla \cdot \u,
\end{align}
and the corresponding conormal derivative $ \displaystyle{\pd{\u}{\nu}}$ on $\p D$ is defined to be
\begin{align}\label{conormal-derivative}
\ds \pd{\u}{\nu}:= \lambda_0(\nabla \cdot \u)\n +\mu_0 \big(\nabla \u+(\nabla \u) ^{T}\big)\n,
\end{align}
where $\n$ is the outward unit normal to $\p D$.

Similarly, we denote  by $\displaystyle \mathcal{L}_{\lambda_1,\mu_1}$
and $\displaystyle \pd{\u}{\widetilde{\nu}}$
 the Lam\'e operator  and the conormal derivative,  respectively, associated to
the Lam\'e constants
$(\lambda_1,\mu_1)$.\\

The problem \eqref{Main-Pb}  is equivalent to the
following problem (see for instance \cite{book, book2, AKNT})
\begin{equation}\label{equation-u}
\left\{
  \begin{array}{ll}
   \ds  \mathcal{L}_{\lambda_0,\mu_0}\u=0 & \mbox{  in }\RR^2\backslash \overline{D},  \\
    \nm \ds \mathcal{L}_{\lambda_1,\mu_1}\u=0  &  \mbox{  in }{D},  \\
   \nm \ds  \u|_{-}=\u|_{+}  &\mbox{ on  } \p D, \\
   \nm\ds \pd{\u}{\widetilde{\nu}}\Big|_{-}=\pd{\u}{{\nu}}\Big|_{+} &\mbox{ on  } \p D,  \\
    \nm \ds   \u(x)-\H(x)=O(|x|^{-1})& \mbox{ as } |x|\rightarrow \infty.
  \end{array}
\right.
\end{equation}
The quantities $ \u|_{\pm}$ on $\p D$ denote the limits from outside and inside of $D$, respectively.
We will also sometimes use $\u^e$ for $\u|_{+}$ and $\u^i$ for
$\u|_{-}.$\\

Let now $D_\ep$ be an $\ep-$perturbation of $D, i.e.$,  there is $h\in \mathcal{C}^1(\p D)$ such that $\p D_\ep$
is given by
\begin{equation}
\ds \p D_{\ep}:=\Big\{\tilde {x}:\tilde{x}=x+\ep h(x)\n(x), x\in \p D\Big\}.
\end{equation}
Let $\u_{\ep}$ be the displacement field  in the presence of $D_\ep$. Then $\u_{\ep}$ is the solution to
\begin{equation}\label{equation-u-ep}
\left\{
  \begin{array}{ll}
   \ds  \mathcal{L}_{\lambda_0,\mu_0}\u_\ep=0 & \mbox{  in }\RR^2\backslash \overline{D}_\ep,  \\
    \nm \ds \mathcal{L}_{\lambda_1,\mu_1}\u_\ep=0 & \mbox{  in }{D_\ep},  \\
   \nm \ds  \u_\ep|_{-}=\u_\ep|_{+}&\mbox{ on  } \p D_\ep, \\
   \nm\ds \pd{\u_\ep}{\widetilde{\nu}}\Big|_{-}=\pd{\u_\ep}{{\nu}}\Big|_{+} &\mbox{ on  } \p D_\ep,  \\
    \nm \ds   \u_\ep(x)-\H(x)=O(|x|^{-1})& \mbox{ as } |x|\rightarrow \infty.
  \end{array}
\right.
\end{equation}

The first  main result of this paper  is the following
 derivation of the leading-order term in the asymptotic expansion of
 $(\u_{\ep}-\u)|_{\Om}$ as $\ep\rightarrow 0$, where $\Om$ is  a bounded
region outside the inclusion $D$, and away from $\p D$.
\begin{thm}\label{Main-theorem} Let  $\u$ and $\u_\ep$ be the solutions to
\eqref{equation-u} and \eqref{equation-u-ep}, respectively. Let $\Om$ be a bounded
region outside the inclusion $D$, and away from $\p D$. For $x \in \Om$, the following pointwise   asymptotic expansion  holds:
\begin{equation}\label{Main-Asymptotic}
\u_\ep(x)=\u(x)+\ep \u_1(x)+O(\ep^2),
\end{equation}
where the remainder $O(\ep^2)$ depends only on $\lambda_0,\lambda_1,\mu_0,\mu_1$,
the $\mathcal{C}^2$-norm of $X$, the $\mathcal{C}^1$-norm of $h$,  $dist (\Om, \p D)$,
and  $\u_1$ is the unique solution to
\begin{equation}\label{equation-u-1}
\left\{
  \begin{array}{ll}
   \ds  \mathcal{L}_{\lambda_0,\mu_0}\u_{1}=0 & \mbox{  in }\RR^2\backslash \overline{D},  \\
    \nm \ds \mathcal{L}_{\lambda_1,\mu_1}\u_{1}=0& \mbox{  in }{D},  \\
   \nm \ds  \u_1|_{-}-\u_1|_{+}= h(\mathbb{K}_{0,1}\widehat{\nabla}\u^{i})\n& \mbox{ on  } \p D,\\
   \nm\ds \pd{\u_1}{\widetilde{\nu}}\Big |_{-}-\pd{\u_1}{{\nu}}\Big|_{+}=
   \frac{\p }{\p \ta}\Big(h\big([\mathbb{C}_{1}-\mathbb{M}_{0,1}]\widehat{\nabla}\u^{i}\big)\ta\Big) & \mbox{ on  } \p D,  \\
    \nm \ds   \u_1(x)=O(|x|^{-1})& \mbox{ as } |x|\rightarrow \infty,
  \end{array}
\right.
\end{equation}
with $\ta$ is the tangential vector to $\p D$,
\begin{align}
\ds  \mathbb{M}_{0,1}&:=\frac{\lambda_0 (\lambda_1+2\mu_1)}{\lambda_0+2\mu_0} \I \otimes \I+2 \mu_1 \mathbb{I}+ \frac{4(\mu_0-\mu_1)(\lambda_0+\mu_0)}{\lambda_0+2\mu_0} \I\otimes (\ta \otimes \ta),\label{Tensor-M}\\
\nm \ds \mathbb{K}_{0,1}&:=\frac{\mu_0(\lambda_1-\lambda_0)+2(\mu_0-\mu_1) (\lambda_0+\mu_0)}{\mu_0(\lambda_0+2\mu_0)} \I \otimes \I+2\big(\frac{\mu_1}{\mu_0}-1\big) \mathbb{I}\nonumber\\
\nm\ds &\q+ \frac{2(\mu_1-\mu_0)(\lambda_0+\mu_0)}{\mu_0(\lambda_0+2\mu_0)} \I\otimes (\ta \otimes\ta ).\label{tensor-B}
\end{align}
\end{thm}

 Our asymptotic
expansion is also valid in the case of an elastic inclusion
with  high contrast parameters,
for more details on the behavior of the  leading and  first  order
terms $\u$ and $\u_1$ in the asymptotic expansion of the
displacement field $\u_\ep$, we refer the reader to \cite[Chapter
2]{ABGKLW}. \\

We should  notice that similar
asymptotic results  have been obtained in the context  of
interface problems in
elastostatics \cite{ AKLZ1, KZ1, KZ2, Zribi1}, the authors derive asymptotic expansions
for boundary displacement field in both cases
 of isotropic and anisotropic thin elastic inclusions and
 perturbations in the eigenvalues and elastic moments tensors
 (EMTs) caused by small perturbations of the shape of an elastic inclusion,
 the approach they use, based on energy estimates, variational  approach,  and
 fine regularity estimates for solutions of elliptic systems with
 discontinuous coefficients obtained by Li and Nirenberg \cite{LN}.
 Unfortunately,  this method does not seem to work in our case.
\\

As a consequence of the results of Theorem \ref{Main-theorem},  we obtain
the following  relationship between traction-displacement measurements
and the deformation $h$. The scalar product in $\RR^2$, will be denoted
by the dot, and sometimes to ease the notation, by $\la ,\ra$.

\begin{thm} \label{second-theorem}  Let $S$
  be a Lipschitz closed curve enclosing $D$, and  away from $\p D$. Let  $\u$ and  $\u_\ep$  be the solutions
 to \eqref{equation-u} and  \eqref{equation-u-ep},  respectively, and  $\v$ be the solution of the following system:
\begin{equation}\label{v}
\left\{
  \begin{array}{ll}
   \ds  \mathcal{L}_{\lambda_0,\mu_0}\v=0 & \mbox{  in }\RR^2\backslash \overline{D},  \\
    \nm \ds \mathcal{L}_{\lambda_1,\mu_1}\v=0  &  \mbox{  in }{D},  \\
   \nm \ds  \v|_{-}=\v|_{+}   &\mbox{ on  } \p D, \\
   \nm\ds \pd{\v}{\widetilde{\nu}}\Big|_{-}=\pd{\v}{{\nu}}\Big|_{+} &\mbox{ on  } \p D,  \\
    \nm \ds   \v (x)-\F(x)=O(|x|^{-1})&\mbox{ as } |x|\rightarrow \infty.
  \end{array}
\right.
\end{equation}

Then, the following asymptotic expansion holds:
\begin{align}\label{asymptotic-traction-displacement}
\ds &\int_{S}\big(\u_\ep-\u\big)\cdot \pd{ \F}{\nu}d\sigma-\int_{S}\big(\pd{\u_\ep}{\nu}-
\pd{ \u}{\nu}\big) \cdot \F d\sigma\nonumber\\
 \nm \ds &  =\ep\int_{\p D}h\bigg(\big(\big[\mathbb{M}_{0,1}-\mathbb{C}_{1}\big]
 \widehat{\nabla}\u^i\big)\ta\cdot \widehat{\nabla} \v^i \ta -\big(\mathbb{K}_{0,1}
  \widehat{\nabla}\u^i\big)\n\cdot (\mathbb{C}_1 \widehat{\nabla} \v^i )\n \bigg)d\sigma+O(\ep^2),
\end{align}
where  the remainder $O(\ep^2)$ depends only on $\lambda_0,\lambda_1,\mu_0,\mu_1$,
the $\mathcal{C}^2$-norm of $X$, the $\mathcal{C}^1$-norm of $h$, and $dist (S,  \p D)$.
\end{thm}
%%%%%%%%%%%%%%%%%%%%%%%%%%%%%%%%%%%%%%%%%ù
%%%%%%%%%%%%%%%%%%%%%%%%%%%%%%%%%%%%%%%

%%%%%%%%%%%%%%%%%%%%%%%%%%%%%%%%%%%%%%%%%ù
%%%%%%%%%%%%%%%%%%%%%%%%%%%%%%%%%%%%%%%

The asymptotic expansion  in \eqref{asymptotic-traction-displacement} can be
used to design new algorithms in the  identification of
 the shape of an elastic inclusion based on traction-displacement
 measurements (see for instance
\cite{AEEKL, AGKLS, AKLZ1, AKLZ2, KKL, LLZ, Zribi1}).\\

The concept of EMTs   has been studied particularly in the context
of imaging of small elastic inclusions \cite{book,   AKNT}. Recall that
EMTs $M_{\alpha \beta }^j:=(m_{\alpha \beta 1}^j,m_{\alpha \beta 2}^j
)$ for
 $\alpha ,\beta \in \NN^2$ and $j=1,2$,
associated to the inclusion $D$ with Lam\'e constants
 $(\lambda_1, \mu_1)$, and  the background medium with
Lam\'e constants $(\lambda_0, \mu_0)$
can be  described in the following manner: consider $\H$ to be a vector-valued function satisfying
$ \mathcal{L}_{\lambda_0,\mu_0}\H=0  \mbox{  in }\RR^2$. Then, the
displacement field $\u$ solution to \eqref{Main-Pb}, resulting from the perturbation of  $\H$  due to
the presence of $D$,  has the  following expansion \cite[Theorem 10.2]{book2}
\begin{equation}\label{definition-classic-EMTs}
\u(x)=\H(x)+\sum_{j=1}^2\sum_{|\alpha|\geq 1}\sum_{|\beta|\geq 1}\frac{1}{\alpha!\beta!}\p^{\alpha}\H_{j}(0) \p^{\beta}\G(x) M_{\alpha \beta }^j \q   \forall x \mbox{ with }  |x|>R,
\end{equation}
where $D \subset B_{R}(0)$ and $\G$ is the fundamental solution to $
\mathcal{L}_{\lambda_0,\mu_0}$.  An alternative definition  of EMTs
will be given in Section 6.\\

The asymptotic expansion of the EMTs has been first obtained in
\cite[Theorem 3.1]{LY} with a remainder of the order of  $O(\ep^{1+\gamma})$ with
$0<\gamma<1 $. The authors  have used an  approach  based on  that method proposed in \cite{AEEKL}.
In  this paper we give an alternative method
 to prove the asymptotic behavior of EMTs resulting from small
perturbations of the shape of an elastic inclusion with $\mathcal{C}^2$-boundary. Its main particularity is the fact that it is based on
integral equations and layer potentials
rather than variational techniques, avoiding the use (and the adaptation to our
context) of the nontrivial regularity results of Li and  Nirenberg \cite{LN}. Our approach  gives a better estimate of the
remainder  (of order $O(\ep^2)$).

\begin{thm}\label{Third-theorem} Let $(a_{j}^\alpha)$ and $(b_{k}^\beta)$
be fixed constants such that $\displaystyle \H(x)=\sum_{j=1}^2\sum_{\alpha\in \NN^2} a^\alpha_j x^\alpha e_j$
 and $\displaystyle \F(x)=\sum_{k=1}^2\sum_{\beta\in \NN^2} b^\beta_k x^\beta e_k$ are
  satisfy $\nabla \cdot \big(\mathbb{C}_0 \widehat{\nabla} \cdot \big)=0$
in $\RR^2$. Let $\u$ and $\v$ be the solutions to \eqref{equation-u} and \eqref{v},
respectively. Then, the following asymptotic expansion holds:
\begin{align}\label{asymptotic-formula-EMTs}
 \ds& \sum_{\alpha \beta j k}a_{j}^{\alpha} b_{k}^{\beta} m_{\alpha \beta k}^j(D_\ep)-
 \sum_{\alpha \beta j k}a_{j}^{\alpha} b_{k}^{\beta} m_{\alpha \beta k}^j(D) \nonumber\\
 \nm \ds & \q \q =\ep\int_{\p D}h\bigg(\big([\mathbb{C}_1-\mathbb{M}_{0,1}]
 \widehat{\nabla}\u^i\big)\ta\cdot \widehat{\nabla} \v^i \ta+\big(\mathbb{K}_{0,1}
 \widehat{\nabla}\u^i\big)\n\cdot (\mathbb{C}_1 \widehat{\nabla} \v^i )\n \bigg)d\sigma+O(\ep^2),
\end{align}
where  the remainder $O(\ep^2)$ depends only on
$\lambda_0,\lambda_1,\mu_0,\mu_1$, the $\mathcal{C}^2$-norm of $X$,
and the $\mathcal{C}^1$-norm of $h$.
\end{thm}

Based on the  asymptotic expansion in \eqref{asymptotic-formula-EMTs},
we can conceive  numerical algorithms  in the spirit of \cite{LY}
to  recover fine shape details from the higher order EMTs.
\\

The  approach and techniques developed in this paper
can be  generalized  to higher dimension  interface problems
and extended to other PDE systems, such as, Stokes and  Maxwell.\\

This paper is organized as follows. In Section 2, we review  some
preliminary results related to small perturbations of a
$\mathcal{C}^2$-interface,  differentiation of tensors, and introduce
a representation of the Lam\'e system in local coordinates. In Section 3, we formally derive the
 asymptotic expansion of the displacement by
 using the field expansion  method (Theorem \ref{Main-theorem}). In Section 4,
we derive the asymptotic expansions of layer potentials. In Section 5,
based on layer potentials techniques, we first justify the formal
expansions, and then find the  relationship between traction-displacement
measurements and the deformation $h$ (Theorem \ref{Main-theorem} \&
Theorem \ref{second-theorem}).
In Section 6, we  rigorously derive the asymptotic formula
 for the perturbation of the EMTs (Theorem \ref{Third-theorem}).
Finally, in the appendix, we provide
some useful integral representations of quantities
related to layer potentials.
%%%%%%%%%%%%%%%%%%%%%%%%%%%%%%%%%
%%%%%%%%%%%%%%%%%%%%%%%%%%%%%%%%%
%%%%%%%%%%%%%%%%%%%%%%%%%%%%%%%%%
%%%%%%%%%%%%%%%%%%%%%%%%%%%%%%%%%
%%%%%%%%%%%%%%%%%%%%%%%%%%%%%%%%%
%%%%%%%%%%%%%%%%%%%%%%%%%%%%%%%%%
%%%%%%%%%%%%%%%%%%%%%%%%%%%%%%%%%
%%%%%%%%%%%%%%%%%%%%%%%%%%%%%%%%%
%%%%%%%%%%%%%%%%%%%%%%%%%%%%%%%%%
%%%%%%%%%%%%%%%%%%%%%%%%%%%%%%%%%
%%%%%%%%%%%%%%%%%%%%%%%%%%%%%%%%%
%%%%%%%%%%%%%%%%%%%%%%%%%%%%%%%%%
%%%%%%%%%%%%%%%%%%%%%%%%%%%%%%%%%
%%%%%%%%%%%%%%%%%%%%%%%%%%%%%%%%%
%%%%%%%%%%%%%%%%%%%%%%%%%%%%%%%%%
%%%%%%%%%%%%%%%%%%%%%%%%%%%%%%%%%
%%%%%%%%%%%%%%%%%%%%%%%%%%%%%%%%%
%%%%%%%%%%%%%%%%%%%%%%%%%%%%%%%%%
%%%%%%%%%%%%%%%%%%%%%%%%%%%%%%%%%
%%%%%%%%%%%%%%%%%%%%%%%%%%%%%%%%%
%%%%%%%%%%%%%%%%%%%%%%%%%%%%%%%%%
%%%%%%%%%%%%%%%%%%%%%%%%%%%%%%%%%
%%%%%%%%%%%%%%%%%%%%%%%%%%%%%%%%%
%%%%%%%%%%%%%%%%%%%%%%%%%%%%%%%%%
%%%%%%%%%%%%%%%%%%%%%%%%%%%%%%%%%
%%%%%%%%%%%%%%%%%%%%%%%%%%%%%%%%%
%%%%%%%%%%%%%%%%%%%%%%%%%%%%%%%%%
\section{Definitions and preliminary results}

\subsection{Small perturbation of a $\mathcal{C}^2$-interface }

Let $a, b \in \RR,$ with $a<b$, and let $X(t): [a,b]\to \RR^2$ be
the arclength parametrization of $\p D$, namely, $X$ is a
$\mathcal{C}^2$-function satisfying $|X'(t)|= 1$ for all $t \in
[a,b]$, $X(a)= X(b)$, and
 $$
 \partial D:=\{x=X(t),  t\in [a,b]\}.
 $$
We assume that $X$ is a positive arclength, $i.e.$, it rotates in the anticlockwise direction. Then the outward unit normal at  $x\in \p D$, $\n(x)$,  is  given by
$\n(x)=R_{-\frac{\pi}{2}}X'(t)$, where $R_{-\frac{\pi}{2}}$ is the
rotation by $-{\pi}/{2}$,  the tangential vector at $x$, $\ta(x) =
X'(t)$, and $X'(t)\perp X''(t)$. Set the curvature $\kappa(x)$ to be
defined by
 $$
 X''(t)=\kappa (x) \n(x).
 $$
We will sometimes use $h(t)$ for $h(X(t))$ and $h'(t)$ for the
tangential derivative of $h(x)$.

Then, $\xe= \tilde{X}(t)=X(t)+\ep h(t)\n(x)=X(t)+\ep
h(t)R_{-\frac{\pi}{2}}X'(t)$ is a parametrization of ${\p D_\ep}$. By
${\n} ({\xe})$ we denote the outward unit
normal to $\p D_\ep$ at $\tilde{x}$. It is proved in
\cite{AKLZ1} that
\begin{align}\label{asymp-n}
\ds {\n}
({\xe})&=\frac{R_{-\frac{\pi}{2}}\tilde{X}'(t)}{|\tilde{X}'(t)|}
=\frac{\Big(1-\ep h(t)\kappa(x)\Big)\n(x)-\ep
h'(t)X'(t)}{\sqrt{\Big(1-\ep
h(t)\kappa(x)\Big)^2+\ep^2{h'(t)}^2}}:=\frac{\etaa(x)}{|\etaa(x)|},
\end{align}
and hence ${\n}({\xe})$ can be expanded uniformly as
 \begin{equation*}
{\n}({\xe}) = \sum_{m=0}^{\infty} \ep^m
 \n_{m}(x), \quad x \in \p D,
 \end{equation*}
where the vector-valued functions $\n_{m}$ are uniformly bounded
regardless of $m$. In particular,
\begin{equation}\label{n0n1}
 \n_{0}(x)=\n(x),\quad  \n_{1}(x) = -h'(t) \ta(x), \quad  x\in \p D.
\end{equation}
Likewise, denote by $d{ \sigma}_{\ep}(\tilde{x})$
the  length element to $\p D_\ep$ at $ \tilde{x}$ which has an
uniformly expansion \cite{AKLZ1}
 \begin{equation} \label{sigexp}
 d{ \sigma}_{\ep}(\tilde{x}) = |
 {\tilde{X}}^\prime(t)|dt=\sqrt{(1-\ep\kappa(t)h(t))^2+\ep^2 h'^2(t)}dt = \sum_{m=0}^\infty
 \ep^m \sigma_{m}(x)
 d\sigma(x),\quad x\in \p D,
 \end{equation}
where $\sigma_{m}$ are  functions bounded regardless of $m$, with
 \begin{equation} \label{sigexp01}
 \sigma_{0}(x)=1, \quad \sigma_{1}(x)=- \kappa(x)h(x),\quad x \in \p D.
 \end{equation}
%%%%%%%%%%%%%%%%%%%%%%%%%%%%%%%%%
%%%%%%%%%%%%%%%%%%%%%%%%%%%%%%%%%
\subsection{Differentiation of tensors}
In this subsection, we will use the Einstein convention for the summation notation. Let $(\e_1,\e_2)$ be an   orthonormal base of $\RR^2$. Let  $\phi$ be a differentiable scalar function. Then
\begin{align}\label{nabla-scalar}
\ds \nabla \phi =\pd{\phi}{x_i}\e_i.
\end{align}
Let  $\u=\u_i \e_i$ be a differentiable vector-valued function.  Then
\begin{equation}\label{nabla-vector}
\ds \nabla \u=\pd{\u}{x_j}\otimes \e_j=\pd{(\u_i \e_i)}{x_j}\otimes \e_j=\pd{\u_i}{x_j}\e_i\otimes \e_j.
\end{equation}
Let $\M=\M_{ij}\e_i\otimes \e_j$ be a differentiable matrix-valued function. Then
\begin{equation}\label{divergence-tensor}
\nabla \cdot \M=\pd{\M}{x_k} \e_k=\pd{(\M_{ij} \e_i\otimes \e_j)}{x_k} \e_k=\pd{\M_{ij}}{x_k} (\e_i\otimes \e_j)\e_k=\pd{\M_{ij}}{x_k} \e_i (\e_j\cdot \e_k)=\pd{\M_{ij}}{x_j} \e_i.
\end{equation}
Also, we have
\begin{equation}\label{nabla-tensor}
 \nabla \M=\pd{\M}{x_k} \otimes \e_k=\pd{(\M_{ij} \e_i\otimes \e_j)}{x_k}\otimes \e_k=\pd{\M_{ij}}{x_k} \e_i\otimes \e_j \otimes \e_k.
\end{equation}
By  $tr(\A)$ we mean the trace of the matrix $\A$. Let  $\v$ be a differentiable vector-valued  function. We have the following properties
\begin{align}
\ds & \nabla (\phi ~\u)=\phi \nabla \u+\u\otimes \nabla \phi, \label{8}\\
\nm\ds & \nabla (\phi ~\M)=\phi \nabla \M+\M\otimes \nabla \phi, \label{9}\\
\nm\ds & \nabla(\u\cdot \v)=(\nabla \u)^{T}\v+(\nabla \v)^{T}\u,\label{10}\\
\nm\ds & \nabla \cdot (\u\otimes \v)= \nabla \u ~\v+\nabla \cdot \v~ \u,\label{11}\\
\nm\ds & \nabla \cdot (\phi ~\M)= \M\nabla \phi +\phi \nabla \cdot \M,\label{12}\\
\nm\ds & \nabla \cdot (\M~\u)= \u\cdot \nabla\cdot (\M^T) + tr(\M \nabla \u),\label{13}
\end{align}

%%%%%%%%%%%%%%%%%%%%%%%%%%%%%%%%%
%%%%%%%%%%%%%%%%%%%%%%%%%%%%%%%%%%%%%%%%%
%%%%%%%%%%%%%%%%%%%%%%%%%%%%%%%%%%%%%%%%%%
%%%%%%%%%%%%%%%%%%%%%%%%%%%%%%%%%%%%%%%
%%%%%%%%%%%%%%%%%%%%%%%%%%%%%%%%%%%%%%%
%%%%%%%%%%%%%%%%%%%%%%%%%%%%%%%%%%%%%%%%%%
%%%%%%%%%%%%%%%%%%%%%%%%%%%%%%%%%%%%%
%%%%%%%%%%%%%%%%%%%%%%%%%%%%%%%%%%%%%%%%%
%%%%%%%%%%%%%%%%%%%%%%%%%%%%%%%%%%%%%%%%%%
%%%%%%%%%%%%%%%%%%%%%%%%%%%%%%%%%%%%%%%
%%%%%%%%%%%%%%%%%%%%%%%%%%%%%%%%%%%%%%%
%%%%%%%%%%%%%%%%%%%%%%%%%%%%%%%%%%%%%%%%%%
%%%%%%%%%%%%%%%%%%%%%%%%%%%%%%%%%%%%%
%%%%%%%%%%%%%%%%%%%%%%%%%%%%%%%%%
%%%%%%%%%%%%%%%%%%%%%%%%%%%%%%%%%
%%%%%%%%%%%%%%%%%%%%%%%%%%%%%%%%%
%%%%%%%%%%%%%%%%%%%%%%%%%%%%%%%%%
%%%%%%%%%%%%%%%%%%%%%%%%%%%%%%%%%
 \subsection{Lam\'e system in local coordinates}
 We begin with a review of some  basic properties of tensor products. Let $\A$ and $\B$ be two  matrices, and let  $\u,\v,$ and $\w$  be 3 vectors. We have
\begin{align}
\nm \ds &(\u\otimes \v) \w=(\v\cdot \w)\u\label{2},\\
\nm \ds &(\u\otimes \v)^{T} =\v\otimes \u \label{6},\\
\nm \ds &(\u\otimes \v\otimes \w)^{T} =\v\otimes \w\otimes \u. \label{7}\\
\nm\ds &\big(\A\otimes (\u\otimes \u)\big) \B=\big((\u\otimes \u): \B \big)\A=\la \B\u, \u\ra \A \label{identity-tensor}.
\end{align}

Let $\w$ be a twice differentiable vector-valued function on $\p D$
and   $(\n,\ta)$ be the  orthonormal base   at each point $x \in \p
D$. Then, the gradient of $\w$ in local coordinates is given by
\begin{align}\label{local-gradient}
\ds \nabla \w= \pd{\w}{\n}\otimes\n+\pd{\w }{\ta}\otimes \ta.
\end{align}
We obtain from  \eqref{nabla-tensor} and \eqref{local-gradient} that
\begin{align*}
\ds \nabla \nabla \w=&\pdl{\w}{\n^2} \otimes \n\otimes \n+\frac{\p^2 \w }{\p \n \p \ta } \otimes \ta \otimes \n
+\frac{\p^2 \w}{\p \ta  \p \n} \otimes \n\otimes \ta+\frac{\p^2 \w }{ \p \ta^2 } \otimes \ta\otimes \ta.\nonumber
\end{align*}
Taking the divergence of \eqref{local-gradient}, we get  from \eqref{divergence-tensor} that
\begin{align}\label{local-Laplacian}
\ds \Delta \w=\nabla\cdot \nabla \w= \pdl{\w}{\n^2}+\pdl{\w}{\ta^2}=\nabla \nabla \w\,\n\,\n+\nabla \nabla \w\,\ta\, \ta.
\end{align}
By using \eqref{nabla-scalar}, we find
\begin{align}\label{local-nabla-divergence}
\ds \nabla \nabla\cdot \w=&\nabla \Big (\pd{\la \w,\n\ra }{\n}+\pd{\la \w,\ta\ra }{\ta}\Big)\nonumber\\
\nm\ds =&\frac{\p^2 \la \w,\n\ra }{ \p \n^2} \n+\frac{\p^2  \la \w,\ta\ra }{\p \n \p \ta} \n
+\frac{\p^2  \la \w,\n\ra }{\p \ta  \p \n} \ta+\frac{\p^2  \la \w,\ta\ra }{\p \ta^2} \ta.
\end{align}
From\eqref{6}  and \eqref{local-gradient}, we deduce  that
\begin{align}\label{local-gradient-transpose}
\ds (\nabla \w)^{T}&=\n\otimes\pd{\w}{\n}+\ta\otimes\pd{\w }{\ta}\nonumber \\
\nm \ds &=\frac{\p \la \w, \n\ra}{\p \n} \n\otimes \n+\frac{\p \la \w, \ta\ra}{\p \n} \n\otimes \ta+
\frac{\p \la \w, \n\ra}{\p \ta} \ta\otimes \n+\frac{\p \la \w, \ta\ra}{\p \ta} \ta\otimes \ta,
\end{align}
and then it follows   from \eqref{divergence-tensor} and \eqref{2} that
\begin{align}\label{local-nabla-transpose}
\ds \nabla \cdot(\nabla \w)^{T}=\frac{\p^2 \la \w, \n\ra }{ \p \n^2} \n+
\frac{\p^2  \la \w, \ta\ra }{\p \ta \p \n} \n+\frac{\p^2  \la \w, \n\ra }{\p \n\p \ta} \ta+\frac{\p^2 \la \w, \ta\ra }{\p \ta^2} \ta.
\end{align}
It is known that $\nabla \cdot(\nabla \w)^{T}=\nabla \nabla\cdot \w$,
which  combined with \eqref{local-nabla-divergence}, and
\eqref{local-nabla-transpose} imply
\begin{align}\label{remark}
\ds \frac{\p^2 \la \w  ,  \ta \ra }{\p \ta \p \n}=\frac{\p^2 \la \w , \ta \ra }{\p \n\p \ta},\q\q \frac{\p^2
\la \w,  \n \ra }{\p \n\p \ta}=\frac{\p^2 \la \w, \n \ra }{\p \ta\p \n} .
\end{align}
Using  \eqref{nabla-tensor} and \eqref{local-gradient-transpose}, we get
\begin{align*}
\ds \nabla (\nabla \w)^{T}=&\n\otimes\pdl{\w}{\n^2}\otimes\n+\ta\otimes\frac{\p^2\w}{\p \n\p \ta}\otimes\n
+\n\otimes\frac{\p^2\w}{\p \ta\p\n}\otimes \ta+
\ta\otimes\pdl{\w}{\ta^2}\otimes \ta,\nonumber
\end{align*}
which gives
\begin{align}\label{local-nabla-nablaTranspose-nn}
\ds \nabla (\nabla \w)^{T}\, \n\, \n =&\pdl{\la \w , \n \ra }{\n^2} \n+\frac{\p^2 \la \w,  \n \ra }{\p \n \p \ta } \ta,\q\q
\nabla (\nabla \w)^{T}\,\ta\, \ta =\frac{\p^2 \la \w , \ta \ra }{\p \ta \p \n } \n+\pdl{\la \w , \ta \ra }{\ta^2} \ta.
\end{align}
Therefore, by \eqref{local-nabla-divergence}, \eqref{remark}, and \eqref{local-nabla-nablaTranspose-nn},  we obtain
 \begin{align}\label{relation-nabla-div-nabla-nabla-transpose}
\ds \nabla \nabla\cdot \w&=\la \nabla \nabla \cdot \w, \n \ra \n+\la \nabla \nabla \cdot \w , \ta\ra \ta=\nabla (\nabla \w)^{T}\, \n\, \n+\nabla (\nabla \w)^{T}\, \ta\, \ta.
\end{align}
%%%%%%%%%%%%%%%%%%%%%%%%%%%%%%%%%%%%%%%%%
%%%%%%%%%%%%%%%%%%%%%%%%%%%%%%%%%%%%%%%%%%

Let $\bphi (x)$ and $\phi(x)$ be respectively a vector and a scalar functions,  which belong to $ \mathcal
{C}^1([a,b])$ for $x=X(\cdot)\in \p D$. By ${d}/{dt}$,  we denote the tangential derivative in the
direction of $\ta(x)=X'(t)$.  We have
\begin{align*}
\ds \frac{d}{dt}\big(\bphi(x)\big)=\nabla \bphi(x) X'(t)=\frac{\p \bphi}{\p
\ta}(x),\q\q\q \frac{d}{dt}\big(\phi(x)\big)=\nabla \phi(x) \cdot X'(t)=\frac{\p \phi}{\p
\ta}(x).
\end{align*}

The following lemma holds.
 \begin{lem}\label{restriction-Lame}
The  restriction of the Lam\'e system $\Lcal_{\lambda_0,\mu_0}$ in $D$ to
a neighborhood  of $\p D$  can be expressed as follows:
\begin{align}\label{local-Lame}
\ds \Lcal_{\lambda_0,\mu_0}\bphi(x)=& \mu_0 \pdl{\bphi}{\n^2}(x)+\lambda_0 \nabla\nabla\cdot \bphi(x)\cdot \n(x)\n(x)+\mu_0 \nabla(\nabla \bphi)^T(x)\n(x)\n(x)\nonumber\\
\nm\ds&-\kappa(x)\pd{ \bphi}{\nu}(x)+\frac{d}{dt}\Big(\big(\mathbb{C}_0\widehat{\nabla} \bphi(x) \big)\ta(x)\Big),\q x\in \p D.
\end{align}
\end{lem}
\proof
According to \eqref{local-Laplacian} and \eqref{relation-nabla-div-nabla-nabla-transpose}. For $x\in \p D$, we have
\begin{align*}
\ds \Lcal_{\lambda_0,\mu_0}\bphi(x)=&\mu_0 \Delta  \bphi(x)+(\lambda_0+\mu_0)\nabla \nabla\cdot \bphi(x)\\
\nm\ds =& \mu_0 \nabla \nabla \bphi(x)\n(x)\n(x)+\lambda_0 \nabla\nabla\cdot \bphi(x)\cdot \n(x)\n(x)+\mu_0 \nabla(\nabla \bphi)^T(x)\n(x)\n(x)\\
\nm\ds&+ \mu_0 \nabla \nabla \bphi(x) \ta(x)\ta(x)+\lambda_0 \nabla\nabla\cdot \bphi(x)\cdot \ta(x)\ta(x)+\mu_0 \nabla(\nabla \bphi)^T(x)\ta(x)\ta(x).
\end{align*}
Since
\begin{align*}
\ds & \mu_0 \nabla \nabla \bphi(x) \ta(x)\ta(x)+\lambda_0 \nabla\nabla\cdot \bphi(x)\cdot \ta(x)\ta(x)+\mu_0 \nabla(\nabla \bphi)^T(x)\ta(x)\ta(x)\\
\nm\ds  &\q\q = \frac{d}{dt}\Big(\mu_0  \nabla \bphi(x) +\lambda_0 \nabla\cdot \bphi(x)+\mu_0 (\nabla \bphi)^T(x)\Big)\ta(x)\\
\nm\ds&\q\q = \frac{d}{dt}\Big(\mathbb{C}_0\widehat{\nabla} \bphi(x)  \Big)\ta(x)\\
\nm\ds &\q \q=-\kappa(x)\big(\mathbb{C}_0\widehat{\nabla} \bphi(x)  \big)\n(x)+\frac{d}{dt}\Big(\big(\mathbb{C}_0\widehat{\nabla} \bphi (x)\big)\ta(x)\Big),\q x\in \p D,
\end{align*}
then \eqref{local-Lame} holds. This completes
the proof.
%%%%%%%%%%%%%%%%%%%%%%%%%%%%%%%%%%%%%%%%%ù
%%%%%%%%%%%%%%%%%%%%%%%%%%%%%%%%%%%%%%%
%%%%%%%%%%%%%%%%%%%%%%%%%%%%%%%%%%%%%
%%%%%%%%%%%%%%%%%%%%%%%%%%%%%%%%%%ùùùùù
%%%%%%%%%%%%%%%%%%%%%%%%%%%%%%%%%%ù
%%%%%%%%%%%%%%%%%%%%%%%%%%%%%%%%%
%%%%%%%%%%%%%%%%%%%%%%%%%%%%%%%%ù
%%%%%%%%%%%%%%%%%%%%%%%%%%%%%%%%
%%%%%%%%%%%%%%%%%%%%%%%%%%%%%%%%%%ù
%%%%%%%%%%%%%%%%%%%%%%%%%%%%%%%%%%
%%%%%%%%%%%%%%%%%%%%%%%%%%%%%%%%%%%%%%%%
%%%%%%%%%%%%%%%%%%%%%%%%%%%%%%%%%%%%%
%%%%%%%%%%%%%%%%%%%%%%%%%%%%%%%%%%%%%%%%%ù
%%%%%%%%%%%%%%%%%%%%%%%%%%%%%%%%%%%%%%%
%%%%%%%%%%%%%%%%%%%%%%%%%%%%%%%%%%%%%
%%%%%%%%%%%%%%%%%%%%%%%%%%%%%%%%%%ùùùùù
%%%%%%%%%%%%%%%%%%%%%%%%%%%%%%%%%%ù
%%%%%%%%%%%%%%%%%%%%%%%%%%%%%%%%%
%%%%%%%%%%%%%%%%%%%%%%%%%%%%%%%%ù
%%%%%%%%%%%%%%%%%%%%%%%%%%%%%%%%
%%%%%%%%%%%%%%%%%%%%%%%%%%%%%%%%%%ù
%%%%%%%%%%%%%%%%%%%%%%%%%%%%%%%%%%
%%%%%%%%%%%%%%%%%%%%%%%%%%%%%%%%%%%%%%%
\section{Formal derivations: the  FE  method}\label{FE-method}
The following observations  are  useful.
\begin{prop} \label{observation-Important} Let  $\u$ be   the solution  to \eqref{Main-Pb}. Then the following identities hold:
\begin{align}
 \ds \big(\mathbb{C}_0 \widehat{\nabla} \u^e\big)\ta &=\big(\mathbb{M}_{0,1} \widehat{\nabla} \u^i\big)\ta,\label{Identity-1}\\
 \nm\ds \big(\mathbb{C}_1 \widehat{\nabla} \u^i\big)\ta &=\big(\mathbb{M}_{1,0} \widehat{\nabla} \u^e\big)\ta,\label{Identity-2}\\
\nm\ds \nabla \u^e \n -\nabla \u^i \n&=\big(\mathbb{K}_{0,1} \widehat{\nabla} \u^i\big)\n=-\big(\mathbb{K}_{1,0} \widehat{\nabla} \u^e\big)\n,\label{identity-3}
\end{align}
where the  $4$-tensors $\mathbb{M}_{l,k}$ and $\mathbb{K}_{l,k}\,\mbox{ for } l,k=0,1,$  are defined by:
\begin{align*}
\ds \mathbb{M}_{l,k}&:=\frac{\lambda_l (\lambda_k+2\mu_k)}{\lambda_l+2\mu_l} \I \otimes \I+2 \mu_k \mathbb{I}+ \frac{4(\mu_l-\mu_k)(\lambda_l+\mu_l)}{\lambda_l+2\mu_l} \I\otimes (\ta \otimes \ta),\\
\nm \ds \mathbb{K}_{l,k}&:=\frac{\mu_l(\lambda_k-\lambda_l)+2(\mu_l-\mu_k) (\lambda_l+\mu_l)}{\mu_l(\lambda_l+2\mu_l)} \I \otimes \I+2\big(\frac{\mu_k}{\mu_l}-1\big) \mathbb{I}\\
\nm\ds &\q + \frac{2(\mu_k-\mu_l)(\lambda_l+\mu_l)}{\mu_l(\lambda_l+2\mu_l)} \I\otimes (\ta \otimes \ta ).
\end{align*}
\end{prop}
\proof The solution  $\u$ of \eqref{Main-Pb}
 satisfies the following transmission
conditions along the interface $\p D$:
\begin{align}
\ds \ds  \u^i&=\u^e,\\
\nm  \ds \nabla \u^i \ta  &= \nabla \u^e \ta, \label{eq01} \\
  \nm \ds \la \widehat{\nabla} \u^i \ta, \ta\ra  &=   \la \widehat{\nabla} \u^e \ta, \ta\ra,\label{eq02} \\
\nm \ds \lambda_1 \nabla \cdot \u^i+2 \mu_1  \la \widehat{\nabla} \u^i \n, \n\ra&=
\lambda_0 \nabla \cdot \u^e+2 \mu_0  \la \widehat{\nabla} \u^e \n, \n\ra,\label{eq03}\\
\nm \ds \mu_1  \la \widehat{\nabla} \u^i \n, \ta\ra&=\mu_0  \la \widehat{\nabla} \u^e \n, \ta\ra.\label{eq04}
\end{align}
Recalling  that
\begin{align}\label{Tr-w}
\ds \nabla \cdot \u^e=\widehat{\nabla} \u^{e}: \I= tr(\widehat{\nabla} \u^{e})=\la\widehat{\nabla} \u^{e}\n,\n \ra+ \la\widehat{\nabla} \u^{e}\ta ,\ta \ra.
\end{align}
From \eqref{eq02}, \eqref{eq03}, and \eqref{Tr-w}, one can easily see that
\begin{align}\label{nabla-cdot-w}
\ds \nabla \cdot \u^e= \frac{\lambda_1+2\mu_1}{\lambda_0+2\mu_0}\nabla \cdot \u^i+\frac{2(\mu_0-\mu_1)}{\lambda_0+2\mu_0}\la \widehat{\nabla} \u^i \ta ,\ta \ra.
\end{align}
We have
\begin{align*}
\ds \nabla \u^{e} \n&=\la\nabla \u^{e} \n, \n\ra \n +\la\nabla \u^{e} \n, \ta \ra \ta \\
\nm\ds &=\la\widehat{\nabla} \u^{e} \n, \n\ra \n +2\la\widehat{\nabla} \u^{e} \n, \ta \ra \ta -\la(\nabla \u^{e})^{T} \n, \ta \ra \ta \\
\nm\ds &=\la\widehat{\nabla} \u^{e} \n, \n\ra \n +2\la\widehat{\nabla} \u^{e} \n, \ta \ra \ta -\la \nabla \u^{e} \ta , \n\ra \ta.
\end{align*}
Using \eqref{Tr-w}, we obtain
\begin{align*}
\ds \nabla \u^{e} \n=(\nabla \cdot \u^{e}) \n -\la\widehat{\nabla} \u^{e} \ta , \ta \ra \n +2\la\widehat{\nabla} \u^{e} \n, \ta \ra \ta -\la \nabla \u^{e} \ta , \n\ra \ta.
\end{align*}
In a similar way, we write
\begin{align*}
\ds \nabla \u^{i} \n&=(\nabla \cdot \u^{i})\n -\la\widehat{\nabla} \u^{i} \ta , \ta \ra \n +2\la\widehat{\nabla} \u^{i} \n, \ta \ra \ta -\la \nabla \u^{i} \ta , \n\ra \ta.
\end{align*}
It then follows from \eqref{identity-tensor}, \eqref{eq01}, \eqref{eq02},  \eqref{eq04}, \eqref{Tr-w}, and \eqref{nabla-cdot-w}, that
\begin{align*}
\ds \nabla \u^{ e} \n-\nabla \u^{i}\n&=(\nabla \cdot \u^{e}-\nabla \cdot \u^{i})\n+2\la\widehat{\nabla} \u^{e} \n, \ta \ra \ta -2\la\widehat{\nabla} \u^{i} \n, \ta \ra \ta \\
\nm \ds &=\Big(\frac{\lambda_1+2\mu_1}{\lambda_0+2\mu_0}-1\Big)(\nabla \cdot \u^{i})\n+2\big(\frac{\mu_1}{\mu_0}-1\big)\la\widehat{\nabla} \u^{i} \n, \ta \ra \ta\\
\nm\ds &\q+\frac{2(\mu_0-\mu_1)}{\lambda_0+2\mu_0} \la \widehat{\nabla} \u^i \ta , \ta \ra \n\\
\nm \ds &=\frac{\mu_0(\lambda_1-\lambda_0)+2(\mu_0-\mu_1) (\lambda_0+\mu_0)}{\mu_0(\lambda_0+2\mu_0)}(\nabla \cdot \u^{i})\n+2\big(\frac{\mu_1}{\mu_0}-1\big)\widehat{\nabla} \u^{i} \n \\
\nm\ds &\q+ \frac{2(\mu_1-\mu_0)(\lambda_0+\mu_0)}{\mu_0(\lambda_0+2\mu_0)} \la \widehat{\nabla} \u^i \ta , \ta \ra \n\\
\nm \ds &=\big(\mathbb{K}_{0,1} \widehat{\nabla} \u^i\big)\n\q\mbox{on }\p D.
\end{align*}
We obtain from \eqref{identity-tensor}, \eqref{eq02}, \eqref{eq03}, \eqref{eq04}, and \eqref{nabla-cdot-w},  that
\begin{align*}
\ds   \big(\mathbb{C}_0 \widehat{\nabla} \u^e \big)\ta &= \lambda_0 (\nabla \cdot \u^e )\ta +2\mu_0 (\widehat{\nabla} \u^e)\ta\\
\nm \ds &=\frac{\lambda_0 (\lambda_1+2\mu_1)}{\lambda_0+2\mu_0}(\nabla \cdot \u^i)\ta+\frac{2\lambda_0(\mu_0-\mu_1)}{\lambda_0+2\mu_0}\la \widehat{\nabla} \u^i \ta, \ta\ra \ta\\
\nm \ds &\q  + 2\mu_0\la \widehat{\nabla} \u^i \ta, \ta\ra \ta+2\mu_1\la \widehat{\nabla} \u^i \ta, \n\ra \n\\
\nm \ds &=\frac{\lambda_0 (\lambda_1+2\mu_1)}{\lambda_0+2\mu_0}(\nabla \cdot \u^i)\ta+\frac{2\lambda_0(\mu_0-\mu_1)}{\lambda_0+2\mu_0}\la \widehat{\nabla} \u^i \ta, \ta\ra \ta\\
\nm \ds &\q+2\mu_1 \widehat{\nabla }\u^i \ta + 2\mu_0\la \widehat{\nabla} \u^i \ta, \ta\ra \ta- 2\mu_1\la \widehat{\nabla} \u^i \ta, \ta\ra \ta\\
\nm \ds &= \frac{\lambda_0 (\lambda_1+2\mu_1)}{\lambda_0+2\mu_0} (\nabla \cdot \u^i)\ta+2\mu_1 \widehat{\nabla }\u^i \ta +\frac{4(\mu_0-\mu_1)(\lambda_0+\mu_0)}{\lambda_0+2\mu_0} \la \widehat{\nabla} \u^i \ta, \ta\ra \ta\\
\nm\ds &=\big(\mathbb{M}_{0,1} \widehat{\nabla} \u^i \big)\ta\q\mbox{on }\p D.
\end{align*}
The identities $(\mathbb{C}_1 \widehat{\nabla}
\u^i\big)\ta=\big(\mathbb{M}_{1,0} \widehat{\nabla} \u^e\big)\ta$ and
$\nabla \u^e \n -\nabla \u^i \n=-\big(\mathbb{K}_{1,0}
\widehat{\nabla} \u^e\big)\n$ can be obtained in exactly the
 same manner as above.  The
proof of the proposition is then achieved.\\
%%%%%%%%%%%%%%%%%%%%%%%%%%%%%%%%%%%%%%%%%ù
%%%%%%%%%%%%%%%%%%%%%%%%%%%%%%%%%%%%%%%
%%%%%%%%%%%%%%%%%%%%%%%%%%%%%%%%%%%%%
%%%%%%%%%%%%%%%%%%%%%%%%%%%%%%%%%%ùùùùù
%%%%%%%%%%%%%%%%%%%%%%%%%%%%%%%%%%ù
%%%%%%%%%%%%%%%%%%%%%%%%%%%%%%%%%
%%%%%%%%%%%%%%%%%%%%%%%%%%%%%%%%ù
%%%%%%%%%%%%%%%%%%%%%%%%%%%%%%%%
%%%%%%%%%%%%%%%%%%%%%%%%%%%%%%%%%%ù
%%%%%%%%%%%%%%%%%%%%%%%%%%%%%%%%%%
%%%%%%%%%%%%%%%%%%%%%%%%%%%%%%%%%%%%%%%%
%%%%%%%%%%%%%%%%%%%%%%%%%%%%%%%%%%%%%
%%%%%%%%%%%%%%%%%%%%%%%%%%%%%%%%%%%%%%%%%ù
%%%%%%%%%%%%%%%%%%%%%%%%%%%%%%%%%%%%%%%
%%%%%%%%%%%%%%%%%%%%%%%%%%%%%%%%%%%%%
%%%%%%%%%%%%%%%%%%%%%%%%%%%%%%%%%%ùùùùù
%%%%%%%%%%%%%%%%%%%%%%%%%%%%%%%%%%ù
%%%%%%%%%%%%%%%%%%%%%%%%%%%%%%%%%
%%%%%%%%%%%%%%%%%%%%%%%%%%%%%%%%ù
%%%%%%%%%%%%%%%%%%%%%%%%%%%%%%%%
%%%%%%%%%%%%%%%%%%%%%%%%%%%%%%%%%%ù
%%%%%%%%%%%%%%%%%%%%%%%%%%%%%%%%%%
%%%%%%%%%%%%%%%%%%%%%%%%%%%%%%%%%%%%%%%

We now derive, based on the FE method \cite{CGHIR}, formally the asymptotic expansion of $\u_\ep$, solution to \eqref{equation-u-ep}, as $\ep$ goes to zero.
We start by expanding  $\u_\ep$  in powers of $\ep$, $i.e.$
\begin{align*}
\ds \u_\ep(x)=\u_0(x)+\ep \u_1(x)+O(\ep^2),\q x \in \Om,
\end{align*}
where $\u_n, n=0,1,$ are  well defined in $\RR^2 \backslash  \p D$,
and satisfy
\begin{equation*}
\left\{
  \begin{array}{ll}
   \ds  \mathcal{L}_{\lambda_0,\mu_0}\u_n=0 & \mbox{  in }\RR^2\backslash \overline{D},  \\
    \nm \ds \mathcal{L}_{\lambda_1,\mu_1}\u_n=0  &  \mbox{  in }{D}, \\
    \nm \ds   \u_n(x)-\H(x)\delta_{0n}=O(|x|^{-1})& \mbox{ as } |x|\rightarrow \infty.
  \end{array}
\right.
\end{equation*}
Here $\delta_{0n}$ is the Kronecker symbol.

Let $\tilde x= x+\ep h(x)\n(x)\in \p D_\ep$  for $x\in \p D$.
The conormal derivative $\displaystyle \pd{\u^e_{\ep}}{\nu}(\xe)$ on $\p D_\ep$  is given by
\begin{align}\label{formal0}
\ds \pd{\u^e_{\ep}}{\nu}(\xe)=\lambda_0 \nabla \cdot \u^e_{\ep}(\xe){{\n}}(\tilde x)
 +\mu_0 \Big(\nabla \u^e_{\ep}(\xe)+(\nabla \u^e_{\ep}) ^{T}(\xe)\Big) {{\n}}(\tilde x),
\end{align}
where ${{\n}}(\tilde{x})$ is the outward unit normal to $\p D_\ep$ at $˜\tilde x$ defined  by  \eqref{asymp-n}. By the Taylor expansion, we write
\begin{align}\label{formal1}
\ds \nabla \cdot \u_{\ep}^{e}(\xe)&= \nabla \cdot\u_{0}^{e}\big(x+\ep h(x)\n(x)\big)+\ep  \nabla \cdot \u_{1}^{e}\big(x+\ep h(x)\n(x)\big)+O(\ep^2)\nonumber\\
\nm \ds &=\nabla \cdot \u_{0}^{e}( x)+\ep h(x) \nabla  \nabla \cdot \u_{0}^{e}( x) \cdot  \n(x)+\ep \nabla \cdot \u_{1}^{e}( x)+O(\ep^2),\q x \in \p D.
\end{align}
In a similar way, we get
\begin{align}\label{formal2}
\ds \nabla \u_{\ep}^{e}(\tilde {x})+(\nabla \u_{\ep}^{e}) ^{T}(\tilde {x})
= &\big[\nabla \u_{0}^{e}(x)+(\nabla \u_{0}^{e}) ^{T}(x)\big]+\ep \big[\nabla \u_{1}^{e}(x)+(\nabla \u_{1}^{e}) ^{T}(x)\big]\nonumber\\
\nm \ds &+ \ep h(x) \big[\nabla \nabla \u_{0}^{e}(x)\n(x)+\nabla(\nabla \u_{0}^{e}) ^{T}(x)\n(x)\big]+O(\ep^2), \q x \in \p D.
\end{align}
It then follows from  \eqref{n0n1}, \eqref{formal0}, \eqref{formal1} and \eqref{formal2} that
\begin{align}\label{expan-conormal}
\ds \pd{\u_{\ep}^{e}}{\nu}(\xe)=&\pd{\u_{0}^{e}}{\nu}(x)+\ep \pd{\u_{1}^{e}}{\nu}(x)- \ep h'(t) \big(\mathbb{C}_0 \widehat{\nabla} \u_0^e(x)\big)\ta(x)\nonumber\\
\nm \ds &+\ep  h(x)\Bigg(\lambda_0 \nabla  \nabla \cdot \u_{0}^{e}( x) \cdot  \n(x)\n(x)+\mu_0\nabla \nabla \u_{0}^{e}(x)\n(x)\n(x)\nonumber\\
\nm\ds& \q\q\q\q\q\q\q\q \q \q+\mu_0\nabla(\nabla \u_{0}^{e}) ^{T}(x)\n(x)\n(x)\Bigg)+O(\ep^2), \q x \in \p D.
\end{align}
Since $\u_0^e$ satisfies $\displaystyle \Lcal_{\lambda_0,\mu_0}\u_0^e=0 $  in $\RR^2 \backslash \overline{D}$, then, by \eqref{local-Lame}, we obtain
\begin{align*}
\ds \mu_0[\nabla \nabla \u_{0}^{e}]\,\n\,\n+\lambda_0 [\nabla \nabla \cdot \u_0^e]\cdot \n\,\n+\mu_0 [\nabla(\nabla \u_{0}^{e}) ^{T}]\,\n\,\n=\kappa\pd{ \u_0^e}{\nu}-\frac{\p }{\p \ta}\Big(\big(\mathbb{C}_0\widehat{\nabla} \u_0^e \big)\ta\Big)\q  \mbox{on } \p D,
\end{align*}
and hence, we derive  from \eqref{expan-conormal} the following formal asymptotic expansion
\begin{align}\label{asymptotic-conormal-e}
\ds \pd{\u_{\ep}^{e}}{\nu}(\tilde x)=&\pd{\u_{0}^{e}}{\nu}(x)+\ep \pd{\u_{1}^{e}}{\nu}(x)+\ep\kappa(x)h(x)\pd{\u_{0}^{e}}{\nu}(x)
-\ep\frac{d}{dt}\Big(h(x)\big[\mathbb{C}_0\widehat{\nabla}\u^{e}_0(x)\big]\ta(x)\Big)\nonumber\\
\nm\ds &+O(\ep^2), \q x\in \p D.
\end{align}
Similarly to \eqref{asymptotic-conormal-e} , we have
 \begin{align}\label{asymptotic-conormal-i}
\ds \pd{\u_{\ep}^{i}}{\widetilde{\nu}}(\tilde x)=&\pd{\u_{0}^{i}}{\widetilde{\nu}}(x)+\ep \pd{\u_{1}^{i}}{\widetilde{\nu}}(x)+\ep\kappa(x)h(x)\pd{\u_{0}^{i}}{\widetilde{\nu}}(x)
-\ep\frac{d}{dt}\Big(h(x)\big[\mathbb{C}_1\widehat{\nabla}\u^{i}_0(x)\big]\ta(x)\Big)\nonumber\\
\nm\ds &+O(\ep^2), \q x\in \p D.
\end{align}
By using  $ \displaystyle \pd{\u_\ep^i}{\widetilde{\nu}}=\pd{\u_\ep^e}{{\nu}}$ on  $\p D_\ep$,  we deduce
from  \eqref{asymptotic-conormal-e} and
 \eqref{asymptotic-conormal-i} that
\begin{align}
\ds \pd{\u_{0}^{i}}{\widetilde{\nu}}=\pd{\u_{0}^{e}}{\nu}& \q \mbox{ on } \p D,\nonumber\\
\nm \ds \pd{\u_{1}^{i}}{\widetilde{\nu}}-\pd{\u_{1}^{e}}{\nu}&=
\frac{\p }{\p \ta }\Big(h\big(\mathbb{C}_1\widehat{\nabla}\u^{i}_0\big)\ta \Big)-\frac{\p }{\p \ta }\Big(h\big(\mathbb{C}_0\widehat{\nabla}\u^{e}_0\big)\ta \Big)\q \mbox{ on } \p D.\label{conormal-u-1}
\end{align}

For $\tilde x=x+\ep h(x) \n(x)\in \p  D_\ep$. We have the following Taylor expansion
\begin{align*}%\label{asymptotic-u-e}
\ds \u_\ep^{e}(\tilde x)&=\u_0^{e}(\tilde x)+\ep \u_1^{e}(\tilde x)+O(\ep^2)\nonumber \\
\nm \ds &=\u^{e}_0(x)+\ep h(x) \nabla \u^{e}_0(x) \n(x)+\ep \u^{e}_1( x)+O(\ep^2), \q x \in \p D.
\end{align*}
Likewise, we obtain
\begin{align*}%\label{asymptotic-u-i}
\ds \u_\ep^{i}(\tilde x)&=\u^{i}_0(x)+\ep h(x) \nabla \u^{i}_0(x) \n(x)+\ep \u^{i}_1( x)+O(\ep^2), \q x \in \p D.
\end{align*}
The transmission condition $\u_\ep^{i}=\u_\ep^{e}$  on $\p D_\ep$,  immediately yields
\begin{align*}
\ds\u^{i}_0=\u^{e}_0\q \mbox{ on }  \p D,
\end{align*}
and
\begin{align}\label{006}
\ds \u^{i}_1-\u^{e}_1=h\big( \nabla \u^{e}_0 \n-\nabla \u^{i}_0 \n\big)\q \mbox{ on }    \p D.
\end{align}
Note that $\u_0=\u$ which is the solution to \eqref{equation-u}.  It then follows from \eqref{conormal-u-1},  \eqref{006}, and Lemma \ref{observation-Important} that
\begin{align}\label{condition-u1-pD}
\ds \u_{1}^{ i}- \u_{1}^{ e}&= h\big(\mathbb{K}_{0,1}\widehat{\nabla}\u^{i}\big)\n \q \mbox{ on } \p D,\\
\nm \ds \pd{\u_{1}^{i}}{\widetilde{\nu}}-\pd{\u_{1}^{e}}{\nu}&=
\frac{\p }{\p \ta }\Big(h\big([\mathbb{C}_1-\mathbb{M}_{0,1}]\widehat{\nabla}\u^{i}\big)\ta \Big) \q \mbox{ on } \p D.
\end{align}
Thus  we formally obtain  Theorem \ref{Main-theorem}, as desired. For a proof, see Subsection \ref{proof}.
%%%%%%%%%%%%%%%%%%%%%%%%%%%%%%%%%%%%%%%%%ù
%%%%%%%%%%%%%%%%%%%%%%%%%%%%%%%%%%%%%%%
%%%%%%%%%%%%%%%%%%%%%%%%%%%%%%%%%%%%%
%%%%%%%%%%%%%%%%%%%%%%%%%%%%%%%%%%ùùùùù
%%%%%%%%%%%%%%%%%%%%%%%%%%%%%%%%%%ù
%%%%%%%%%%%%%%%%%%%%%%%%%%%%%%%%%
%%%%%%%%%%%%%%%%%%%%%%%%%%%%%%%%ù
%%%%%%%%%%%%%%%%%%%%%%%%%%%%%%%%
%%%%%%%%%%%%%%%%%%%%%%%%%%%%%%%%%%ù
%%%%%%%%%%%%%%%%%%%%%%%%%%%%%%%%%%
%%%%%%%%%%%%%%%%%%%%%%%%%%%%%%%%%%%%%%%%
%%%%%%%%%%%%%%%%%%%%%%%%%%%%%%%%%%%%%
%%%%%%%%%%%%%%%%%%%%%%%%%%%%%%%%%%%%%%%%%ù
%%%%%%%%%%%%%%%%%%%%%%%%%%%%%%%%%%%%%%%
%%%%%%%%%%%%%%%%%%%%%%%%%%%%%%%%%%%%%
%%%%%%%%%%%%%%%%%%%%%%%%%%%%%%%%%%ùùùùù
%%%%%%%%%%%%%%%%%%%%%%%%%%%%%%%%%%ù
%%%%%%%%%%%%%%%%%%%%%%%%%%%%%%%%%
%%%%%%%%%%%%%%%%%%%%%%%%%%%%%%%%ù
%%%%%%%%%%%%%%%%%%%%%%%%%%%%%%%%
%%%%%%%%%%%%%%%%%%%%%%%%%%%%%%%%%%ù
%%%%%%%%%%%%%%%%%%%%%%%%%%%%%%%%%%
%%%%%%%%%%%%%%%%%%%%%%%%%%%%%%%%%%%%%%%%
\section{Asymptotic formulae of layer potentials}
\subsection{Layer potentials}
Let us  review some well-known properties of the layer potentials on a
Lipschitz domain for the
elastostatics.  \\

Let
\begin{align*}
\ds \Psi:=\Big \{\bpsi: \p_i \bpsi_j+\p_j \bpsi_i=0,\q 1\leq i,j\leq 2\Big \}.
\end{align*}
or equivalently,
\begin{align*}
\ds \Psi=\mbox{span}\Bigg\{\theta_1(x):={1 \brack 0},\, \theta_2(x):={0 \brack 1},\,\theta_3(x):={x_2 \brack -x_1}\Bigg\}.
\end{align*}
Introduce the space
\begin{align*}
\ds \displaystyle L^2_{\Psi}(\p D):=\Big \{\f \in L^2(\p D): \int_{\p D} \f\cdot \bpsi ~d\sigma=0 \mbox{ for all } \bpsi \in \Psi \Big\}.
\end{align*}
In particular, since $\Psi$ contains constant functions, we get
\begin{align*}
 \ds \int_{\p D} \f d\sigma=0
\end{align*}
for any $\f \in L^2_{\Psi}(\p D)$. The following fact is useful later.
\begin{align}\label{L-Psi}
\ds \mbox{If} \q \w \in W^{1,\frac{3}{2}}(D)  \q\mbox{satisfies}\q \Lcal_{\lambda_0,\mu_0} \w=0  \mbox{ in } D,\q\mbox{then}\q \pd{\w}{\nu}\Big|_{\p D}\in L^2_{\Psi}(\p D).
\end{align}

The  Kelvin matrix of fundamental  solution $\G$ for the  Lam\'e
system $\Lcal_{\lambda_0,\mu_0}$ in $\RR^2$,  is known to be
\begin{align}\label{Kelvin}
\ds \G(x)=\frac{A}{2\pi}\log |x|\I-\frac{B}{2\pi} \frac{x\otimes x }{|x|^2},\q x \neq 0,
\end{align}
where
\begin{align*}
\ds A=\frac{1}{2}\Big(\frac{1}{\mu_0}+\frac{1}{2\mu_0+\lambda_0}\Big) \q\q \mbox{and}\q \q B=\frac{1}{2}\Big(\frac{1}{\mu_0}-\frac{1}{2\mu_0+\lambda_0}\Big).
\end{align*}
The single and double layer potentials of the density function  $\bphi$ on $L^2(\p D)$
associated with the Lam\'e parameters $(\lambda_0, \mu_0)$  are defined by
 \begin{align}
\ds \BScal_{D} [\bphi](x)
 =& \int_{\partial D} \G(x-y) \bphi(y)d\sigma(y), \quad x \in
 \RR^2,\label{single-layer}\\
\nm \ds \BDcal_{D} [\bphi] (x)
=&\int_{\p D}\bigg(\lambda_0 \nabla_y \cdot \G(x-y)\otimes \n(y)\nonumber\\
\nm\ds &\q\q \q+ \mu_0 \Big(\big[\nabla_y \G(x-y)\n(y)\big]^{T}+(\nabla_y \G)^{T}(x-y)\n(y)\Big)\bigg)\bphi(y)d\sigma(y)\nonumber\\
\nm\ds :=&\int_{\p D}\KK(x-y)\bphi(y)d\sigma(y), \quad x \in \RR^2 \setminus \p
D. \label{double-layer}
\end{align}

The followings are  well-known properties of  the single and double layer potentials due to Dahlberg, Keing, and Verchota \cite{DKV}. Let $D$ be a Lipschitz  bounded domain in $\RR^2$. Then, we have
\begin{align}
\ds \pd{\BScal_D [\bphi] }{\nu}\Big |_{\pm}(x) & = \Big (\pm \frac{1}{2} \I
+ \BKcal_D^* \Big ) [\bphi] (x) \quad
\mbox{a.e. } x \in \p D \label{nuS}, \\
 \nm \ds \BDcal_D [\bphi] \big|_{\pm} (x) & = \Big(\mp \frac{1}{2}
\I + \BKcal_D \Big) [\bphi](x) \quad \mbox{a.e. } x \in \p D,
\label{doublejump-h}
\end{align}
where $\BKcal_D$ is defined  by
\begin{align*}
\ds \BKcal_D [\bphi](x)={p. v.}\int_{ \p D}\KK(x-y) \bphi(y) d\sigma(y) \quad \mbox{a.e. } x \in \p D,
\end{align*}
 and  $\BKcal_D^*$ is the adjoint operator of $\BKcal_D $, that is,
\begin{align}\label{KD-star}
\ds \BKcal_D^*  [\bphi](x)=&{p.v.}\int_{\p D}\KK^{T}(x-y)\bphi(y)d\sigma(y)\nonumber\\
\nm \ds=&{p. v.}\int_{\p D}  \bigg(\lambda_0 \n(x)\otimes\nabla_x \cdot \G(x-y)\nonumber\\
 \nm \ds &\q + \mu_0 \Big(\nabla_x \G(x-y)\n(x)+\big[(\nabla_x \G)^{T}(x-y)\n(x)\big]^{T}\Big)\bigg)\bphi(y)d\sigma(y)\quad \mbox{a.e. } x \in \p D,
\end{align}
with
\begin{align}\label{hat-K}
\ds\KK^{T}(x-y)=& \frac{1}{2\pi}\frac{(A-B)}{(A+B)} ~\frac{\la x-y,\n(x)\ra }{|x-y|^2}\I+\frac{1}{2\pi}\frac{(A-B)}{(A+B)} \frac{(x-y)\otimes \n(x)-\n(x)\otimes (x-y)}{|x-y|^2}\nonumber\\
\nm\ds &+\frac{2}{\pi}\frac{B}{(A+B)}\frac{\la x-y,\n(x)\ra }{|x-y|^2}\frac{(x-y)\otimes(x-y) }{|x-y|^2}\q \mbox{ for } x,y\in \p D, \q x\neq y.
\end{align}
Here ${p.v.}$ denotes the Cauchy principal value. The operators $\BKcal_D $ and $\BKcal_D^*$ are singular integral operators and bounded on $L^2(\p D)$.

Even though the derivation of the kernel $\KK^{T}(x-y)$   is easy,
we give its proof for the reader's convenience.
Denote by $\x:=x-y$, one can easily see  from \eqref{11} and  \eqref{12} that
\begin{align*}
\ds \displaystyle \nabla_x \cdot \G(\x)=\frac{A-B}{2\pi}\frac{\x}{|\x|^2},
\end{align*}
and hence
\begin{align}\label{part1}
\n(x)\otimes \nabla_x \cdot \G(\x)= \frac{(A-B)}{2\pi}\frac{\n(x)\otimes \x }{|\x|^2}.
\end{align}
It follows from \eqref{nabla-tensor} and \eqref{7} that
\begin{align}\label{nabla-x-times-x}
\ds \nabla_x(\x\otimes \x )=\nabla_x(\x_i\x_je_i\otimes e_j)=\frac{\p(\x_i\x_j)}{\p x_k} e_i\otimes e_j\otimes e_k&=(\x_j\delta_{ik}+\x_i \delta_{jk})e_i\otimes e_j\otimes e_k\nonumber\\
\nm\ds &=\x_je_k\otimes e_j\otimes e_k+\x_i e_i\otimes e_k\otimes e_k\nonumber\\
\nm\ds&=(\I\otimes \x)^{T}+ (\x\otimes \I).
\end{align}
Here we used the Einstein convention for the summation notation.

From \eqref{8}, \eqref{9},  \eqref{10}, and \eqref{nabla-x-times-x}, we get
\begin{align}\label{transpose-gradient-gamma }
\ds \nabla_x \G(\x)=\frac{A}{2\pi}\frac{\I\otimes \x}{|\x|^2}+\frac{B}{\pi} \frac{\x\otimes \x \otimes \x}{|\x|^4}-\frac{B}{2\pi}\frac{(\I\otimes \x)^{T}+ (\x\otimes \I)}{|\x|^2},
\end{align}
and thus
\begin{align*}
\ds \nabla_x \G(\x)\n(x)=\frac{A}{2\pi}\frac{\la  \x, \n(x)\ra}{|\x|^2}\I+\frac{B}{\pi} \frac{ \la  \x, \n(x)\ra}{|\x|^4}(\x\otimes \x)-\frac{B}{2\pi}\frac{\x\otimes \n(x)+ \n(x)\otimes \x}{|\x|^2}.
\end{align*}
Using \eqref{7}, the transpose of $\nabla_x \G(\x)$ is given by
\begin{align*}
\ds (\nabla_x \G)^{T}(\x)=\frac{A}{2\pi}\frac{(\I\otimes \x)^{T}}{|\x|^2}+\frac{B}{\pi} \frac{\x\otimes \x \otimes \x}{|\x|^4}-\frac{B}{2\pi}\frac{(\x\otimes \I)+(\I\otimes \x)}{|\x|^2},
\end{align*}
and hence  we obtain
\begin{align*}
\ds (\nabla_x \G)^{T}(\x)\n(x)=\frac{A}{2\pi}\frac{\n(x)\otimes \x}{|\x|^2}+\frac{B}{\pi} \frac{ \la \x, \n(x)\ra}{|\x|^4}\x\otimes \x-\frac{B}{2\pi}\frac{\x\otimes \n(x)}{|\x|^2}-\frac{B}{2\pi}\frac{\la\x, \n(x)\ra}{|\x|^2}\I.
\end{align*}
Therefore
\begin{align}\label{kernel-adjoint}
\ds \nabla_x \G(\x)\n(x)+[(\nabla_x \G)^{T}(\x)\n(x)]^{T}=&\frac{(A-B)}{2\pi}\frac{\la \x,\n(x)\ra }{|\x|^2}\I+\frac{(A-B)}{2\pi}\frac{\x \otimes \n(x) }{|\x|^2}\nonumber\\
\nm\ds &-\frac{B}{\pi}\frac{\n(x)\otimes\x}{|\x|^2}+\frac{2 B}{\pi} \frac{\la \x,\n(x)\ra}{|\x|^4} (\x\otimes \x).
\end{align}
We finally  get  $\KK^{T}(x-y)$ in \eqref{hat-K}  from \eqref{part1} and \eqref{kernel-adjoint}, as desired.

Let $\BDcal_{D}^{\sharp}$  be the standard  double layer potential which is  defined for any $\bphi \in L^2(\p D)$ by
\begin{align}\label{dcal-sharp}
\ds \BDcal_{D}^{\sharp}[\bphi](x)=\int_{\p D} \pd{\G(x-y)}{\n(y)}\bphi(y)d\sigma(y), \q x\in \RR^2\backslash {\p D}.
\end{align}
One can easily see that
\begin{align}\label{dcal-sharp-presentation}
\ds  \pd{\G(x-y)}{\n(y)}\bphi(y)=&\bigg[-\frac{A}{2\pi}\frac{\la \x, \n(y) \ra}{|\x|^2}\I-\frac{B}{\pi}\frac{\la \x, \n(y) \ra}{|\x|^4} (\x\otimes \x)+\frac{B}{2\pi}\frac{\x\otimes \n(y)+\n(y)\otimes \x} {|\x|^2}\bigg]\bphi(y)\nonumber\\
\nm\ds=&-\frac{A}{2\pi} \frac{\la \x,\n(y) \ra}{|\x|^2}\bphi(y)-\frac{B}{\pi}
\frac{\la\x,\n(y)\ra\la \x,\bphi(y)\ra}{|\x|^4}\x\nonumber\\
\nm\ds &+\frac{B}{2\pi}\frac{\la \n(y),\bphi(y)\ra \x +\la\x,\bphi(y)
\ra\n(y)}{|\x|^2}\nonumber\\
\nm\ds :=& \Lambda_1(x,y)+\Lambda_2(x,y)+\Lambda_3(x,y)\q \mbox{for }x\neq y.
\end{align}
For $i=1,2,3$, it follows from \eqref{8}-\eqref{13}  that
\begin{align*}
\ds \Lcal_{\lambda_0,\mu_0}\big(\Lambda_i(\cdot,y)\big)(x)=C_i
\Bigg(\frac{\la\bphi(y),\n(y)\ra}{|\x|^4}\x&+ \frac{\la\x,\n(y)\ra}{|\x|^4}\bphi(y)+\frac{\la\x,\bphi(y)\ra}{|\x|^4}\n(y)\\
&\q\q\q\q \q -4\frac{\la\x,\n(y)\ra\la\x,\bphi(y)\ra}{|\x|^6}\x \Bigg)\q \mbox{for }x\neq y,
\end{align*}
with
\begin{align*}
\ds C_1=\frac{(\lambda_0+\mu_0)A}{\pi},\q\q C_2=-\frac{2\mu_0B}{\pi},\q\q C_3=-\frac{(\lambda_0+\mu_0)B}{\pi}.
\end{align*}
Since $C_1+C_2+C_3=0$, then  $\BDcal_D^{\sharp}(\bphi)$ satisfies
\begin{align}\label{D-sharp-equation}
\ds \Lcal_{\lambda_{0},\mu_0}\big(\BDcal_D^{\sharp}[\bphi]\big)=0 \q \mbox{in  } \RR^2\backslash{\p D}.
\end{align}
The following proposition holds.

\begin{prop}\label{D-sharp} Let $D$ be a bounded Lipschitz domain in
  $\RR^2$. For $\bphi \in L^2(\p D)$, we have
\begin{align}
 \ds \BDcal^{\sharp}_D [\bphi] \big|_{\pm} (x) & = \Big(\mp \frac{1}{2\mu_0}~\I\pm B ~\n \otimes \n
 + \BKcal^{\sharp}_D \Big) \bphi(x) \quad \mbox{a.e. } x \in \p D,\label{D-Sharp-n}\\
\nm \ds \pd{\BScal_D [\bphi]}{\n} \Big |_{\pm}(x) & = \Big (\pm \frac{1}{2\mu_0}~\I\mp B ~\n \otimes \n
+ \big(\BKcal_D^\sharp\big)^* \Big ) \bphi (x) \quad
\mbox{a.e. } x \in \p D,\label{S-Sharp-n}
\end{align}
where  $\BKcal^\sharp_D$ is defined  by
\begin{align*}
\ds \BKcal^\sharp_D [\bphi](x)=\mbox{p.v.} \int_{ \p D}\pd{}{\n(y)} \G(x-y) \bphi(y) d\sigma(y) \quad \mbox{a.e. } x \in \p D,
\end{align*}
and  $\big(\BKcal^\sharp_D\big)^*$ is the adjoint operator of $\BKcal^\sharp_D $,
 that is,
\begin{align}\label{KD-star-sharp}
\ds \big(\BKcal^\sharp_D\big)^*  [\bphi](x)=  p.v. \int_{ \p D}\pd{}{\n(x)} \G(x-y) \bphi(y) d\sigma(y)\quad \mbox{a.e. } x \in \p D.
\end{align}
The operators $\BKcal^\sharp_D$ and $(\BKcal^\sharp_D)^*$ are singular integral operators and bounded on $L^2(\p D)$.
\end{prop}
\proof  Standard arguments yield the trace formulas \cite{FJR}
\begin{align}\label{Important-limit}
 \ds \p_i\big(\BScal_{D}[\bphi]\big)_j(x)\big|_{\pm}&=\pm  \Big\{\frac{1}{2\mu_0}\n_i(x)\bphi_j(x)-B \la \bphi,  \n \ra \n_{i}(x)\n_{j}(x)\Big\}\nonumber\\
 \nm\ds \q&\q + p.v. \int_{\p D}\p_{i}\G_{jk}(x-y)\bphi_k(y)d\sigma(y),\q x \in \p D,
\end{align}
namely,
\begin{align}\label{limit-relation}
\ds \nabla \BScal_D[\bphi](x)\big|_{\pm}=&\pm \Big\{\frac{1}{2\mu_0}\bphi(x)\otimes\n(x)-B \la  \n(x), \bphi (x)\ra \n(x)\otimes \n(x)\Big\}\nonumber\\
\nm\ds &+ p.v.\int_{ \p D}\nabla_x \big[\G(x-y) \bphi(y)\big] d\sigma(y),\q x \in \p D.
\end{align}
Clearly the jump relation of the normal derivative of the single layer potential in \eqref{S-Sharp-n} follows from \eqref{2} and \eqref{limit-relation}. The jump formula in \eqref{D-Sharp-n} can be proved by using standard arguments from the proof  of the theorem $3.28$ in \cite{ Folland76}. The operators $\BKcal^\sharp_D $ and $\big(\BKcal^\sharp_D\big)^*$ are bounded on  $L^2(\p D)$ by the theorem of Coifman-McIntosh-Meyer \cite{CMM82}.

The operators  $\BDcal^{\sharp}_D $ and ${\p \BScal_D}/{\p \n}$ can be viewed as  unfamiliar layer potentials for  the system of elastostatics.

Note that we  will drop the $p.v.$ in the below; this is because $\p D$ is $\mathcal{C}^2$ and throughout this paper  we will denote  by $\displaystyle \widetilde{\BScal}_D, \widetilde{\BDcal}_D, \,\widetilde{\BKcal}^*_D,\, \widetilde{\BDcal}^{\sharp}_D,$ and $(\widetilde{\BKcal}^{\sharp}_D)^*$  the   layer potentials  corresponding to the Lam\'e constants $(\lambda_1,\mu_1)$.

 Let us note  simple, but important relations.
\begin{lem}\label{relations}
\begin{enumerate}
 \item  If $\f\in W^{1,2}(D)$ and $\Lcal_{\lambda_0,\mu_0}\f=0$ in $D$, then for all $\g \in W^{1,2}(D)$,
\begin{align}\label{Important-relation}
 \ds \displaystyle  \int_{\p D} \g\cdot \pd{\f}{{\nu}} ~d\sigma=\int_{D} \lambda_0 (\nabla \cdot \f)(\nabla \cdot \g)+\frac{\mu_0}{2}\big(\nabla \f+(\nabla \f) ^{T}\big): \big(\nabla \g+(\nabla \g) ^{T}\big) d\sigma.
\end{align}
 \item If $\f \in W^{1,2}(\RR^2\backslash \overline{D})$ and $\Lcal_{\lambda_0,\mu_0}\f=0$ in $\RR^2\backslash \overline{D}$, $\f(x)=O(|x|^{-1})$  as $|x|\rightarrow \infty$. Then for all $\g \in W^{1,2}(\RR^2\backslash \overline{D})$,  $\g(x)=O(|x|^{-1})$  as $|x|\rightarrow \infty$, we have
\begin{align}\label{Important-relation-2}
\ds  - \displaystyle  \int_{\p D} \g\cdot \pd{\f}{{\nu}} ~d\sigma=\int_{\RR^2\backslash \overline{D}} \lambda_0 (\nabla \cdot \f)(\nabla \cdot \g)+\frac{\mu_0}{2}\big(\nabla \f+(\nabla \f) ^{T}\big): \big(\nabla \g+(\nabla \g) ^{T}\big) d\sigma.
\end{align}
Here,  for $2\times 2$ matrices $\M$ and $\N$,  $\M: \N= \displaystyle \sum_{ij}\M_{ij}\N_{ij}$.
\end{enumerate}
\end{lem}

\subsection{Asymptotic formula of $\BKcal_{D_\ep}^*$}
Let $\xe$, $\ye \in \p D_\ep$, that is,
\begin{align}\label{xtile-ytilde}
\ds \xe= x+\ep h(x)\n(x),\q\q\q  \ye= y+\ep h(y)\n(y),\q x,y\in\p D.
\end{align}
Denote by
\begin{align*}
\ds E(x,y):=h(x)\n(x)-h(y)\n(y).
\end{align*}
It follows from \eqref{xtile-ytilde} that
\begin{align}\label{absolute-value-2}
\ds |\xe-\ye|^2=|x-y|^2 \Big(1+2\ep F(x,y)+\ep^2 G(x,y)\Big),
\end{align}
where
\begin{align*}
\ds F(x,y)=\frac{\la x-y,E(x,y)\ra}{|x-y|^2}, \q\q G(x,y)=\frac{|E(x,y)|^2}{|x-y|^2}.
\end{align*}
Since $\p D$ is of class $\mathcal{C}^2$, then
\begin{align*}
\ds \frac{\la x-y,\n(x)\ra}{|x-y|^2}, \frac{\la x-y,\n(y)\ra}{|x-y|^2}\leq C \q  \mbox{ for } x,y \in\p D.
\end{align*}
We have $h\nu \in \mathcal{C}^1(\p D)$. Then, one can easily see that
\begin{align*}
\ds |F(x,y)|+|G(x,y)|^{\frac{1}{2}}\leq C \|X\|_{\mathcal{C}^2}\|h\|_{\mathcal{C}^1} \q  \mbox{ for } x,y \in\p D.
\end{align*}
We denote by $|\cdot|_{\infty}$ the matrix infinity  norm. For $x,y \in \p D$, we have
\begin{align*}
\ds \bigg|\frac{(x-y)\otimes (x-y)}{|x-y|^2}  \bigg| _{\infty}\leq 1,
\end{align*}
and
\begin{align*}
\ds \bigg|\frac{E(x,y)\otimes (x-y)}{|x-y|^2}\bigg|_{\infty}, \bigg|\frac{(x-y)\otimes E(x,y)}{|x-y|^2}\bigg|_{\infty}, \bigg|\frac{E(x,y)\otimes E(x,y)}{|x-y|^2}\bigg|^{\frac{1}{2}}_{\infty}\leq C \|X\|_{\mathcal{C}^2}\|h\|_{\mathcal{C}^1}.
\end{align*}

For $\widetilde \bphi \in L^2(\p D_\ep)$, the operator $\BKcal_{D_{\ep}}^*$ is defined by
\begin{align*}
\ds \BKcal_{D_\ep}^*  [\widetilde\bphi](\xe)=&\int_{\p D_{\ep}}\KK^{T}(\xe-\ye)\widetilde{\bphi}(\ye)d{\sigma}_{\ep}(\ye), \nonumber
\end{align*}
where
\begin{align*}
\ds\KK^{T}(\xe-\ye)=& \frac{1}{2\pi}\frac{(A-B)}{A+B} ~\frac{\la \xe-\ye,\n(\xe)\ra }{|\xe-\ye|^2}\I+\frac{1}{2\pi}\frac{(A-B)}{(A+B)} \frac{(\xe-\ye)\otimes \n(\xe)-\n(\xe)\otimes (\xe-\ye)}{|\xe-\ye|^2}\nonumber\\
\nm\ds &+\frac{2}{\pi}\frac{B}{A+B}\frac{\la \xe-\ye,\n(\xe)\ra }{|\xe-\ye|^2}\frac{(\xe-\ye)\otimes(\xe-\ye) }{|\xe-\ye|^2}\q \mbox{ for } \xe,\ye \in \p D_\ep , ~ \xe\neq \ye.
\end{align*}
It follows from \eqref{asymp-n}, \eqref{sigexp}, and \eqref{absolute-value-2} that
\begin{align*}
\ds& \frac{(\xe-\ye)\otimes \n(\xe)-\n(\xe)\otimes(\xe-\ye)}{|\xe-\ye|^2}d\sigma_\ep(\ye)\\
\nm\ds& =\frac{(\xe-\ye)\otimes \etaa(x)-\etaa(x)\otimes(\xe-\ye)}{|x-y|^2}\\
\nm\ds&\q\q\q\q\q\q  \times \frac{1}{1+2\ep F(x,y)+\ep^2G(x,y)} \frac{\sqrt{\Big(1-\ep h(y)\kappa(y)\Big)^2+\ep^2h'(s)^2}}{\sqrt{\Big(1-\ep h(x)\kappa(x)\Big)^2+\ep^2h'(t)^2}}d\sigma(y).
\end{align*}
We have
\begin{align}\label{term1}
\ds&\frac{1}{1+2\ep F(x,y)+\ep^2G(x,y)}\times  \frac{\sqrt{\Big(1-\ep h(y)\kappa(y)\Big)^2+\ep^2h'(s)^2}}{\sqrt{\Big(1-\ep h(x)\kappa(x)\Big)^2+\ep^2h'(t)^2}}d\sigma(y)\nonumber\\
\nm\ds&=\Big[1-2\ep \frac{\la x-y,h(x)\n(x)-h(y)\n(y)\ra}{|x-y|^2}+\ep\Big(\kappa(x)h(x)-\kappa(y)h(y)\Big)\Big] d\sigma(y)+O(\ep^2),
\end{align}
where the remainder $O(\ep^2)$ depends only  on the
$\mathcal{C}^2$-norm of $X$ and $\mathcal{C}^1$-norm of $h$.\\

According to \eqref{asymp-n} and  \eqref{xtile-ytilde}, we write
\begin{align}\label{asymptotic-tensor-product}
\ds \frac{(\xe-\ye)\otimes \etaa(x)- \etaa(x)\otimes(\xe-\ye)}{|x-y|^2}&=\big(1-\ep\kappa(x)h(x)\big)\frac{(x-y)\otimes \n(x)-\n(x)\otimes(x-y)}{|x-y|^2}\nonumber\\
\nm \ds &\q + \ep \frac{E(x,y)\otimes \n(x)-\n(x)\otimes E(x,y)}{|x-y|^2}\nonumber\\
\nm\ds  &\q - \ep h'(t)\frac{(x-y)\otimes \ta(x)-\ta(x)\otimes (x-y)}{|x-y|^2}+O(\ep^2).
\end{align}
Therefore, by \eqref{term1} and  \eqref{asymptotic-tensor-product}, we get
\begin{align}\label{asymp-term1}
\ds & \frac{(\xe-\ye)\otimes \n(\xe)-\n(\xe)\otimes(\xe-\ye)}{|\xe-\ye|^2}d\sigma_\ep(\ye)\nonumber\\
\nm\ds &=\frac{\x \otimes \n(x)-\n(x)\otimes \x}{|\x|^2}d\sigma(y)+\ep\Big(\kappa(x)h(x)-\kappa(y)h(y)\Big) \frac{\x \otimes \n(x)-\n(x)\otimes \x}{|\x|^2}d\sigma(y)\nonumber\\
\nm\ds&\q-\ep \kappa(x) h(x) \frac{\x \otimes \n(x)-\n(x)\otimes \x}{|\x|^2}d\sigma(y)-\ep h'(t) \frac{\x \otimes \ta(x)-\ta(x)\otimes \x}{|\x|^2}d\sigma(y)\nonumber\\
\nm\ds&\q+ \ep h(y)\bigg[2 \frac{\la \x ,\n(y)\ra}{|\x|^4}\Big(\x \otimes \n(x)-\n(x)\otimes \x\Big)-\frac{\n(y)\otimes \n(x)-\n(x)\otimes \n(y)}{|\x|^2}\bigg]d\sigma(y)\nonumber\\
\nm\ds &\q - 2 \ep h(x) \frac{\la \x ,\n(x)\ra}{|\x|^4}\Big(\x \otimes \n(x)-\n(x)\otimes \x\Big)d\sigma(y)+O(\ep^2).
\end{align}
It is proved in \cite{AKLZ1} that
\begin{align}\label{asymp-term3}
\ds \frac{\la \tilde{x}-\tilde{y},\n(\tilde x)\ra }{|\tilde{x}-\tilde{y}|^2}d\sigma_\ep(\ye)
=&\frac{\la \x,\n(x)\ra }{|\x|^2}d\sigma(y)+\ep \Big(\kappa(x)h(x)-\kappa(y)h(y)\Big)\frac{\la \x,\n(x)\ra }{|\x|^2}d\sigma(y)\nonumber\\
\nm\ds &-\ep \kappa(x)h(x)\frac{\la \x,\n(x)\ra }{|\x|^2}d\sigma(y)-\ep h'(t)\frac{\la \x,\ta(x)\ra }{|\x|^2}d\sigma(y)\nonumber\\
\nm\ds &+\ep h(y)\bigg(2 \frac{\la \x,\n(y)\ra\la \x,\n(x)\ra}{|\x|^4}-\frac{\la \n(x),\n(y)\ra }{|\x|^2}\bigg)d\sigma(y)\nonumber\\
\nm\ds&+ \ep h(x)\bigg(-2 \Big(\frac{\la \x,\n(x)\ra}{|\x|^2}\Big)^2+\frac{1}{|\x|^2}\bigg)d\sigma(y)+O(\ep^2).
\end{align}
Using \eqref{xtile-ytilde} and \eqref{absolute-value-2}, we obtain
\begin{align}\label{asymp-term4}
\ds \frac{(\xe-\ye)\otimes(\xe-\ye)}{|\xe-\ye|^2}=&\frac{(\xe-\ye)\otimes(\xe-\ye)}{|x-y|^2}\times \frac{1}{1+2\ep F(x,y)+\ep^2G(x,y)}\nonumber\\
\nm\ds=&\bigg(1-2\ep\frac{\la x-y, h(x)\n(x)-h(y)\n(y)\ra}{|x-y|^2}\bigg) \frac{(x-y)\otimes(x-y)}{|x-y|^2}\nonumber\\
\nm\ds &+ \ep \frac{(x-y)\otimes \Big(h(x)\n(x)-h(y)\n(y)\Big)}{|x-y|^2}\nonumber\\
\nm\ds &+ \ep \frac{\Big(h(x)\n(x)-h(y)\n(y)\Big)\otimes(x-y) }{|x-y|^2}
+O(\ep^2).
\end{align}
It follows from  \eqref{asymp-term3} and \eqref{asymp-term4} that
\begin{align}\label{asymp-term5}
\ds& \frac{\la \xe-\ye,\n(\tilde x)\ra }{|\tilde{x}-\tilde{y}|^2}\frac{(\xe-\ye)\otimes(\xe-\ye)}{|\xe-\ye|^2}d\sigma_\ep(\ye)\nonumber\\
\nm\ds&\q  =\Bigg[ \frac{\la\x ,\n( x)\ra }{|\x|^4}(\x\otimes \x )+\ep\Big(\kappa(x)h(x)-\kappa(y)h(y)\Big)\frac{\la\x ,\n( x)\ra }{|\x|^4}(\x\otimes \x )\nonumber\\
\nm\ds&\q\q -\ep \kappa(x)h(x)\frac{\la\x ,\n( x)\ra }{|\x|^4}(\x\otimes \x )-\ep h'(t)\frac{\la\x ,\ta( x)\ra }{|\x|^4}(\x\otimes \x )\nonumber\\
\nm\ds &\q\q + \ep h(x)\Bigg(-4 \frac{(\la \x ,\n( x)\ra)^2 }{|\x |^6}(\x \otimes \x)+\frac{(\x \otimes \x)}{|\x|^4}
+\frac{\la \x ,\n( x)\ra }{|\x|^4}\Big(\n(x)\otimes\x +\x \otimes\n(x)\Big) \Bigg)\nonumber\\
\nm\ds &\q\q + \ep h(y)\Bigg(4 \frac{\la \x ,\n( y)\ra \la \x ,\n( x)\ra }{|\x |^6}(\x\otimes \x )-\frac{\la \n(x),\n(y)\ra}{|\x|^4}(\x\otimes \x)\nonumber\\
\nm\ds &\q\q\q\q\q\q\q\q\q\q\q\q\q-\frac{\la \x ,\n( x)\ra }{|\x |^4}\Big(\n(y)\otimes\x +\x \otimes\n(y)\Big)\Bigg)\Bigg]d\sigma(y)+O(\ep^2).
\end{align}
From \eqref{asymp-term1}, \eqref{asymp-term3},  and \eqref{asymp-term5}, we write
\begin{align*}
\ds \KK^{T}(\xe-\ye)=\KK^{T}(x-y)+\ep  \KK_1(x-y)+O(\ep^2) \q \mbox{for } x, y \in \p D, x\neq y.
\end{align*}
Introduce the integral operator $\Kcal_D^{(1)}$, defined for any $\bphi \in L^2(\p D)$ by
\begin{align}\label{KD-1-verion0}
\ds \BKcal_{D}^{(1)}[\bphi](x):=\int_{\p D}  \KK_1(x-y)\bphi(y)d\sigma(y),\q x\in \p D.
\end{align}
The operator $\BKcal_{D}^{(1)}$ is bounded on $L^2(\p D)$. In fact, this is an
immediate consequence of the celebrated theorem of Coifman-McIntosh-Meyer \cite{CMM82}.

Let $\Phi_\ep$ be the diffeomorphism from $\p D$ onto $\p D_\ep$  given by $\Phi_\ep(x) = x+\ep h(t)\n(x)$,
where $x = X(t) \in \p D$. The following theorem holds.

\begin{thm}  \label{theorem-KD-1} There exists $C>0$ depending only on
  $\lambda_0,\lambda_1,\mu_0,\mu_1$, $\|X\|_{\mathcal{C}^2}$, and
  $\|h\|_{\mathcal{C}^1}$ such that for any $\widetilde\bphi \in
  L^2(\p D_\ep)$, we have
\begin{align}\label{theorem-KD-star}
\bigg\|\BKcal^*_{D_{\ep}}[\widetilde \bphi] \circ \Phi_{\ep}- \BKcal^*_{D}[\bphi]-\ep\BKcal^{(1)}_{D}[\bphi]
 \bigg\|_{L^2(\p D)}\leq C \ep^2 \big\|\bphi \big\|_{L^2(\p D)},
\end{align}
where $\bphi=\widetilde {\bphi} \circ \Phi_\ep$ and $\BKcal_D^{(1)}$ is defined in \eqref{KD-1-verion0}.
\end{thm}

The following theorem is of particular importance to us in order
to establish our asymptotic
expansions.
\begin{thm} \label{conormal-D-sharp-prop} Let $\bphi \in \mathcal{C}^{1,k}(\p D),$ for some $0<k<1$. Then
\begin{align}\label{conormal-D-sharp-identity}
\ds \frac{\p\BDcal_{D}^{\sharp}[\bphi]}{\p \nu}\Big|_{+}-\frac{\p\BDcal_{D}^{\sharp}[\bphi]}{\p \nu}\Big|_{-}&=\frac{\p}{\p  \ta }\Big(\la \bphi, \ta \ra \n +\frac{\lambda_0}{2\mu_0+\lambda_0}\la \bphi, \n \ra \ta \Big) \q \mbox{ on } \p D.
\end{align}
\end{thm}
\proof
For a function $\w$ defined on $\RR^2\backslash \p D$, we denote
\begin{align*}
\ds \w(x)|_{\pm}=\lim_{t\neq0, t\rightarrow 0^{\pm}}\w(x_t) \q \mbox{ for } x\in \p D, \q x_t:=x +t\n(x).
\end{align*}
Let $\widetilde \bphi \in L^2(\p D_\ep)$ and $\bphi=\widetilde{\bphi}
\circ \Phi_\ep$. Following the same arguments as in the case of
$\BKcal_{D_\ep}$ (taking $h=1$) and using the integral representations in the appendix, we can prove that
\begin{align*}
\ds &\Big(\pd{\BScal_{D_\ep}}{{\nu}}[\widetilde\bphi] \circ \Phi_{\ep}-\pd{\BScal_{D}}{{\nu}}[ \bphi]\Big)(x_t)\\
\nm \ds &=\ep \bigg(\kappa(x)\frac{\p \BScal_{D}[\bphi]}{\p\nu}(x_t)- \frac{\p \BScal_{D}[\kappa \bphi]}{\p\nu}(x_t)+\frac{\p\BDcal_{D}^{\sharp}[\bphi]}{\p \nu}(x_t)
 -\kappa(x)\frac{\p \BScal_{D}[\bphi]}{\p\nu}(x_t)\\
 \nm\ds &\q\q+\lambda_0 \nabla \nabla\cdot \BScal_D[\bphi](x_t)\cdot \n(x_t)\n(x_t)
+\mu_0 \nabla \Big(\nabla \BScal_D[\bphi](x_t)
+\big(\nabla \BScal_D[\bphi](x_t)\big)^{T}\Big)\n(x_t)\n(x_t)\bigg)\\
\nm\ds&\q +O(\ep^2).
\end{align*}
If $\bphi \in \mathcal{C}^{1,k}(\p D),$ then $\BScal_D[\bphi]$ is  $\mathcal{C}^{2,k}$ and $\BDcal^{\sharp}_D[\bphi]$ is $\mathcal{C}^{1,k}$  on  $\overline{D}$ and $\RR^2\backslash D$. Thus
\begin{align*}
&\ds \Big(\pd{\BScal_{D_\ep}}{{\nu}}[\widetilde \bphi] \circ \Phi_{\ep}-\pd{\BScal_{D}}{{\nu}}[ \bphi]\Big)\bigg|_{\pm}\\
\nm\ds&
\q=\ep \bigg(\kappa\frac{\p \BScal_{D}[\bphi]}{\p\nu}- \frac{\p \BScal_{D}[\kappa \bphi]}{\p\nu}+\frac{\p\BDcal_{D}^{\sharp}[\bphi]}{\p \nu} -\kappa\frac{\p \BScal_{D}[\bphi]}{\p\nu}+\lambda_0 \nabla \nabla\cdot \BScal_D[\bphi]\cdot \n\n\nonumber\\
\nm\ds &\q\q\q\q\q  +\mu_0 \nabla \Big(\nabla \BScal_D[\bphi]
+\big(\nabla \BScal_D[\bphi]\big)^{T}\Big)\n\n\bigg)\bigg|_{\pm}+O( |t|^k\ep)+O(\ep^2) \q \mbox{on }\p D.
\end{align*}
Since $\Lcal_{\lambda_0,\mu_0}\BScal_D[\cdot]=0$ in  $ \RR^2\backslash \p D$, it  follows from   the representation of the Lam\'e system on $\p D$ in \eqref{local-Lame} that
\begin{align*}
\ds &\frac{\p}{\p \ta}\Big(\big(\mathbb{C}_0 \widehat{\nabla}{\BScal_D[\bphi]}\big)\ta\Big)\Big|_{\pm}\\
\nm \ds &\,  =\bigg(\kappa\frac{\p \BScal_{D}[\bphi]}{\p\nu}-\lambda_0 \nabla \nabla\cdot \BScal_D[\bphi]\cdot \n\n
-\mu_0 \nabla \Big(\nabla \BScal_D[\bphi]+\big(\nabla \BScal_D[\bphi]\big)^{T}\Big)\n\n\bigg)\bigg|_{\pm}\q \mbox{on } \p D,
\end{align*}
and hence
\begin{align*}
\ds \Big(\pd{\BScal_{D_\ep}}{{\nu}}[\widetilde \bphi] \circ \Phi_{\ep}-\pd{\BScal_{D}}{{\nu}}[ \bphi]\Big)\bigg|_{\pm}
&=\ep \bigg(\kappa\frac{\p \BScal_{D}[\bphi]}{\p\nu}- \frac{\p \BScal_{D}[\kappa \bphi]}{\p\nu}+\frac{\p\BDcal_{D}^{\sharp}[\bphi]}{\p \nu}\\
\nm\ds& \q\q\,\,\,-\frac{\p }{\p \ta }\Big(\big(\mathbb{C}_0 \widehat{\nabla}{\BScal_D[\bphi]}\big)\ta\Big)\bigg)\bigg|_{\pm}  +O( |t|^k\ep)+O(\ep^2)\,\,\,\mbox{on } \p D.
\end{align*}
According to \eqref{nuS}, we have
\begin{align*}
\ds \bigg(\pd{\BScal_{D_\ep}[\widetilde \bphi]}{{\nu}} \circ \Phi_{\ep}-\pd{\BScal_{D}[ \bphi]}{{\nu}}\bigg)\bigg|_{+}
=\bigg(\pd{\BScal_{D_\ep}[\widetilde \bphi]}{{\nu}} \circ \Phi_{\ep}-\pd{\BScal_{D}[ \bphi]}{{\nu}}\bigg)\bigg|_{-}\q \mbox{on }\p D,
\end{align*}
which gives
\begin{align}\label{jump}
\ds &\kappa\frac{\p \BScal_{D}[\bphi]}{\p\nu}\Big|_{+}- \frac{\p \BScal_{D}[\kappa \bphi]}{\p\nu}\Big|_{+}+\frac{\p\BDcal_{D}^{\sharp}[\bphi]}{\p \nu}\Big|_{+}-\frac{\p }{\p \ta }\Big(\big(\mathbb{C}_0 \widehat{\nabla}{\BScal_D[\bphi]}\big)\ta \Big)\Big|_{+}\nonumber\\
\nm\ds &=\kappa\frac{\p \BScal_{D}[\bphi]}{\p\nu}\Big|_{-}- \frac{\p \BScal_{D}[\kappa \bphi]}{\p\nu}\Big|_{-}+\frac{\p\BDcal_{D}^{\sharp}[\bphi]}{\p \nu}\Big|_{-}-\frac{\p }{\p \ta}\Big(\big(\mathbb{C}_0 \widehat{\nabla}{\BScal_D[\bphi]}\big)\ta \Big)\Big|_{-}\q \mbox{on } \p D.
\end{align}
By  \eqref{nuS} again, we have
\begin{align*}
\ds &\bigg(\kappa\frac{\p \BScal_{D}[\bphi]}{\p\nu}- \frac{\p \BScal_{D}[\kappa \bphi]}{\p\nu}\bigg)\bigg|_{+}=\bigg(\kappa\frac{\p \BScal_{D}[\bphi]}{\p\nu}- \frac{\p \BScal_{D}[\kappa \bphi]}{\p\nu}\bigg)\bigg|_{-}\q \mbox{on } \p D.
\end{align*}
It then  follows from  \eqref{jump} that
\begin{align}
\ds \frac{\p\BDcal_{D}^{\sharp}[\bphi]}{\p \nu}\Big|_{+}-\frac{\p }{\p \ta }\Big(\big(\mathbb{C}_0 \widehat{\nabla}{\BScal_D[\bphi]}\big)\ta \Big)\Big|_{+}=\frac{\p\BDcal_{D}^{\sharp}[\bphi]}{\p \nu}\Big|_{-}-\frac{\p }{\p \ta }\Big(\big(\mathbb{C}_0 \widehat{\nabla}{\BScal_D[\bphi]}\big)\ta\Big)\Big|_{-}\q \mbox{on } \p D, \label{jump-continuous}
\end{align}
that is,
$
\frac{\p\BDcal_{D}^{\sharp}[\bphi]}{\p \nu}-\frac{\p}{\p \ta }\Big(\big(\mathbb{C}_0 \widehat{\nabla}{\BScal_D[\bphi]}\big)\ta \Big) \mbox{ is continuous on } \p D,
$
but
$
\frac{\p\BDcal_{D}^{\sharp}[\bphi]}{\p \nu}$ and\\
$\frac{\p }{\p \ta}\Big(\big(\mathbb{C}_0 \widehat{\nabla}{\BScal_D[\bphi]}\big)\ta \Big)$  are discontinuous on  $\p D$, and we have the following relationship
\begin{align*}
\ds \frac{\p\BDcal_{D}^{\sharp}[\bphi]}{\p \nu}\Big|_{+}-\frac{\p\BDcal_{D}^{\sharp}[\bphi]}{\p \nu}\Big|_{-}=\frac{\p }{\p\ta }\Big(\big(\mathbb{C}_0 \widehat{\nabla}{\BScal_D[\bphi]}\big)\ta \big|_{+}- \big(\mathbb{C}_0 \widehat{\nabla}{\BScal_D[\bphi]}\big)\ta \big|_{-}\Big)\q \mbox{on } \p D.
\end{align*}
It follows from  \eqref{Important-limit} that
$$
\big(\mathbb{C}_0 \widehat{\nabla}{\BScal_D[\bphi]}\big)\ta \big|_{+}- \big(\mathbb{C}_0 \widehat{\nabla}{\BScal_D[\bphi]}\big)\ta \big|_{-}= \la \bphi, \ta \ra \n +\frac{\lambda_0}{2\mu_0+\lambda_0}\la \bphi, \n \ra \ta \q \mbox{ on } \p D.
$$
Thus   \eqref{conormal-D-sharp-prop} is proved, as desired. This
finishes the proof of the theorem.\\

As a direct consequence of \eqref{nuS}, \eqref{jump-continuous}, and
the expansions in the appendix,  the integral representation of $\BKcal^{(1)}_{D}$ in \eqref{KD-1-verion0}, can be rewritten as
\begin{align}\label{KD-1-final}
\ds \BKcal_{D}^{(1)}[\bphi](x) =&\bigg(\kappa h(x) \frac{\p \BScal_{D}[\bphi]}{\p\nu}(x)- \frac{\p \BScal_{D}[\kappa h\bphi]}{\p\nu}(x)\bigg)\bigg|_{\pm}\nonumber\\
\nm\ds &+\bigg(\frac{\p\BDcal_{D}^{\sharp}[h\bphi]}{\p \nu}(x)
-\frac{d}{dt}\Big(h(x)\big(\mathbb{C}_0 \widehat{\nabla}{\BScal_D[\bphi](x)}\big)\ta(x)\Big)\bigg)\bigg|_{\pm},\q  x\in \p D.
\end{align}

\subsection{Asymptotic expansion of $\BScal_{D_{\ep}}$ }
For $\widetilde \bphi \in L^2(\p D_\ep)$, we have
\begin{align*}
\ds \BScal_{D_\ep}[\widetilde \bphi](\xe) &= \int_{\p D_\ep}\bigg(\frac{A}{2\pi}\log|\xe-\ye| -\frac{B}{2\pi} \frac{(\xe-\ye)\otimes(\xe-\ye) }{|\xe-\ye|^2}\bigg)\widetilde{\bphi}(\ye) d\sigma_\ep(\ye), \q \tilde x \in \p D_\ep.
\end{align*}
It follows from \eqref{sigexp} and \eqref{absolute-value-2} that
\begin{align}\label{asymp-term6}
\ds &\log |\xe-\ye|d\sigma_\ep(\tilde y )\nonumber\\
\nm\ds &=\frac{1}{2}\log\Big(|x-y|^2\big(1+2\ep F(x,y)+\ep^2G(x,y)\big)\Big)d\sigma_\ep(\tilde y )\nonumber\\
\nm\ds&= \Big(\log |x-y|+\ep F(x,y)+O(\ep^2)\Big)\times \Big(d\sigma(y)-\ep \kappa(y)h(y)d\sigma(y)+O(\ep^2)\Big)\nonumber\\
\nm\ds&= \Bigg[\log |x-y|+\ep \bigg(-\kappa(y)h(y) \log|x-y|+h(x)\frac{\la x-y,\n(x)\ra}{|x-y|^2}
-h(y)\frac{\la x-y,\n(y)\ra}{|x-y|^2}\bigg)\Bigg]d\sigma(y)\nonumber\\
\nm \ds &\q +O(\ep^2)\Big(\log |x-y|+1\Big).
\end{align}
According to \eqref{sigexp}, \eqref{asymp-term4},  and  \eqref{asymp-term6}, we obtain
\begin{align}\label{kernel-SD}
 \ds&\G(\xe-\ye)d\sigma(\ye)\nonumber\\
\nm\ds &=\Bigg[\G(\x)-\ep \kappa(y) h(y)\G(\x)\nonumber\\
\nm\ds&\q+\ep h(x)\Bigg(\frac{A}{2\pi}\frac{\la \x, \n(x) \ra}{|\x|^2}\I+\frac{B}{\pi}\frac{\la \x, \n(x) \ra}{|\x|^2} \frac{\x\otimes \x} {|\x|^2}-\frac{B}{2\pi}\frac{\x\otimes \n(x)+\n(x)\otimes \x} {|\x|^2}\Bigg)\nonumber\\
\nm\ds&\q+\ep h(y)\Bigg(-\frac{A}{2\pi}\frac{\la \x, \n(y) \ra}{|\x|^2}\I-\frac{B}{\pi}\frac{\la \x, \n(y) \ra}{|\x|^2} \frac{\x\otimes \x} {|\x|^2}+\frac{B}{2\pi}\frac{\x\otimes \n(y)+\n(y)\otimes \x} {|\x|^2}\Bigg)\Bigg]d\sigma(y)\nonumber\\
\nm\ds&\q+O(\ep^2)\Big(\log |\x|+1).
\end{align}
Introduce an integral operator $\BScal_D^{(1)}$, defined for any $\bphi \in L^2(\p D)$ by
\begin{align}\label{SD-1}
\BScal_{D}^{(1)}[\bphi](x)=-\BScal_{D}[\kappa h \bphi](x)+\Big(h(x)\frac{\p\BScal_D [\bphi]}{\p\n}(x)+\BDcal_{D}^{\sharp}[h\bphi](x)\Big)\Big|_{\pm}, \q x\in \p D.
\end{align}
The operators $\BScal_{D}^{(1)}$ and $\pd{\BScal_{D}^{(1)}}{ \ta }$  are bounded on $L^2(\p D)$ by the theorem of Coifman, McIntosh, and Meyer \cite{CMM82}. Therefore, we get from \eqref{kernel-SD}
\begin{align}\label{single-L2}
\Big\|\BScal_{D_{\ep}}[\widetilde \bphi] \circ \Phi_{\ep}- \BScal_{D}[\bphi]-\ep\BScal^{(1)}_{D}[\bphi]
 \Big\|_{L^2(\p D)}\leq C \ep^2 \big\|\bphi \big\|_{L^2(\p D)},
\end{align}
where $\bphi=\widetilde \bphi \circ \Phi_\ep$.

We have
\begin{align*}
\ds \pd{\BScal_{D_{\ep}}[\widetilde \bphi] }{\ta}(\xe)=\int_{\p D} \nabla \G\big(\xe-\Phi_\ep(y)\big)R_{\frac{\pi}{2}}\etaa(x) \bphi (y)\times \frac{\sqrt{\Big(1-\ep h(y)\kappa(y)\Big)^2+\ep^2h'(s)^2}}{\sqrt{\Big(1-\ep h(x)\kappa(x)\Big)^2+\ep^2h'(t)^2}}d\sigma(y),
\end{align*}
where $\nabla \G$ and $\etaa$ are defined in \eqref{transpose-gradient-gamma } and \eqref{asymp-n}, respectively. Following the same argument as in the case of  $\BKcal_{D_\ep}^*$, we can prove that
\begin{align}\label{single-tangente-L2}
\bigg\|\pd{\BScal_{D_{\ep}}[\widetilde \bphi]}{\ta}\circ \Phi_{\ep}- \pd{\BScal_{D}[\bphi]}{\ta}-\ep\pd{\BScal^{(1)}_{D}[\bphi]}{\ta}
 \bigg\|_{L^2(\p D)}\leq C \ep^2 \big\|\bphi \big\|_{L^2(\p D)}.
\end{align}
Throughout this paper $W_1^2(\p D)$ denotes the first $L^2$-Sobolev of
space of order $1$ on $\p D$.  From \eqref{single-L2} and
\eqref{single-tangente-L2}, we obtain
the following theorem.

\begin{thm} \label{theorem-SD-1} There exists $C>0$ depending only on $\lambda_0,\lambda_1,\mu_0,\mu_1$,
$\|X\|_{\mathcal{C}^2}$, and $\|h\|_{\mathcal{C}^1}$ such that for any $\widetilde \bphi \in L^2(\p D_\ep)$,
\begin{align}\label{theorem-SD}
\bigg\|\BScal_{D_{\ep}}[\widetilde \bphi] \circ \Phi_{\ep}- \BScal_{D}[\bphi]-\ep\BScal^{(1)}_{D}[\bphi]
 \bigg\|_{W_1^2(\p D)}\leq C \ep^2 \big\|\bphi \big\|_{L^2(\p D)},
\end{align}
where $\bphi=\widetilde\bphi \circ \Phi_\ep$ and $\BScal_D^{(1)}$ is defined in \eqref{SD-1}.
\end{thm}

\section{Asymptotic of the displacement field}

The following
solvability result done by Escauriaza and Seo \cite{ES}.
\begin{thm} \label{ES-theorem} Suppose that $(\lambda_0-\lambda_1)(\mu_0-\mu_1)\geq 0$ and  $0<\lambda_1,\mu_1<+\infty.$ For any given $(\F,\GG)\in W^{2}_1(\p D)\times L^2(\p D),$ there exists a unique pair $(\f,\g)\in L^{2}(\p D)\times L^2(\p D)$ such that
\begin{align}\label{ES-system}
\left\{
  \begin{array}{lll}
  \ds  \widetilde{\BScal}_{D}[\f]\big |_{-}- \BScal_{D}[\g]\big |_{+}
   =\F &\mbox{ on }\p D, \\
    \nm \ds \big(-\frac{1}{2}\I+\widetilde{\BKcal}^*_{D}\big)[\f]- \big(\frac{1}{2}\I+{\BKcal}^*_{D}\big)[\g]
   =\GG &\mbox{ on }\p D,   \\
  \end{array}
\right.
\end{align}
and there exists a constant $C>0$ depending only on $\lambda_0,\mu_0,\lambda_1,\mu_1,$ and the Lipschitz character of $D$ such that
\begin{align}\label{ES-bounds}
\ds \|\f\|_{L^2(\p D)}+\|\g\|_{L^2(\p D)}\leq C\Big(\|\F\|_{W^{2}_1(\p D)}+\|\GG\|_{L^2(\p D)}\Big).
\end{align}
Moreover, if $\GG \in L^2_{\Psi}(\p D)$, then $\g \in L^2_{\Psi}(\p D)$.
\end{thm}
The following proposition is of particular importance to us.
\begin{prop} Suppose that $(\lambda_0-\lambda_1)(\mu_0-\mu_1)\geq 0$ and $0<\lambda_1,\mu_1<+\infty.$ For any given $(\F,\GG)\in W^{2}_{1}(\p D)\times L^2(\p D)$,
there exists a unique pair $(\f,\g)\in L^2(\p D)\times L^2(\p D)$ such that
\begin{equation}\label{solvability}
\left\{
  \begin{array}{lll}
   \ds  \Big(\widetilde{\BScal}_{D}+\ep \widetilde{\BScal}_D^{(1)}\Big)[\f]- \Big(\BScal_{D}+\ep \BScal_D^{(1)}\Big)[\g]
   =\F\q \mbox{on }\p D, \\
    \nm \ds \Big(-\frac{1}{2}\I+\widetilde{\BKcal}_{D}^*+\ep \widetilde{\BKcal}_D^{(1)}\Big)[\f]- \Big(\frac{1}{2}\I+{\BKcal}_{D}^*+\ep {\BKcal}_D^{(1)}\Big)[\g]
   =\GG\q \mbox{on }\p D.  \\
  \end{array}
\right.
\end{equation}
Furthermore, there exists a constant $C>0$ depending only on  $\lambda_0, \mu_0, \lambda_1, \mu_1$,
and the Lipschitz character of $D$ such that
\begin{align}\label{estimate}
\|\f\|_{L^2(\p D)}+\|\g\|_{L^2(\p D)}\leq C \Big(\|\F\|_{W_{1}^2(\p D)}+\|\GG\|_{L^2(\p D)}\Big).
\end{align}
\end{prop}
\proof Let $\mathcal{X}:=L^2(\p D)\times L^2(\p D)$ and $\mathcal{Y}:=W_{1}^2(\p D)\times L^2(\p D)$. For $n=0,1$, define the operator $\Tcal_n: \mathcal{X} \rightarrow \mathcal{Y}$ by
\begin{align*}
\ds \Tcal_0(\f,\g):=\bigg(\widetilde{\BScal}_D[\f]\big|_{-}-\BScal_D [\g]\big|_{+}, \Big(-\frac{1}{2}\I+\widetilde{\BKcal}_{D}^*\Big)[\f] - \Big(\frac{1}{2}\I+{\BKcal}_{D}^*\Big)[\g]\bigg),
\end{align*}
and
\begin{align*}
\ds \Tcal_1(\f,\g):=\bigg(\widetilde{\BScal}^{(1)}_D[\f]\big|_{-}-\BScal^{(1)}_D [\g]\big|_{+}, \widetilde{\BKcal}_D^{(1)}[\f] -\BKcal_D^{(1)}[\g] \bigg).
\end{align*}
The operator $\Tcal_1$ is bounded on $ \mathcal{X}$ because it is a
linear combination of bounded  integral operators. According to
Theorem  \ref{ES-theorem}, the operator $\Tcal_0$ is invertible. For $\ep$ small enough, it follows from Theorem $1.16$, section $4$  of \cite{Kato}, that  the operator $\Tcal_0+\ep \Tcal_1$ is invertible. This completes the proof of solvability of \eqref{solvability}. The estimate \eqref{estimate} is a consequence of solvability and the closed graph  theorem.

\subsection{Representation of solutions}

For more details on the following representation formulae, we refer to \cite{book, book2, AKNT}. The solution $\u_\ep$ to \eqref{equation-u-ep} can be represented as
\begin{equation}\label{Representation-u-ep}
\u_\ep(x)=
\left\{
  \begin{array}{lll}
   \ds  \H(x)+\BScal_{D_\ep}[\bvarphi_\ep](x), \q x\in \RR^2 \backslash \overline{D}_\ep, \\
    \nm \ds\widetilde{ \BScal}_{D_\ep}[\bpsi_\ep](x),\q x\in D_\ep,  \\
  \end{array}
\right.
\end{equation}
where the pair $(\bpsi_\ep,\bvarphi_\ep)$ is the unique solution in $L^2(\p D_\ep)\times L^2_{\Psi}(\p D_\ep)$ of
\begin{equation}\label{phi-psi-ep}
\left\{
  \begin{array}{lll}
   \ds  \widetilde{\BScal}_{D_\ep}[\bpsi_\ep]\big|_{-}- \BScal_{D_\ep}[\bvarphi_\ep]\big |_{+}
   =\H \q \mbox{on }\p D_\ep, \\
    \nm \ds \big(-\frac{1}{2}\I+\widetilde{\BKcal}^*_{D_\ep}\big)[\bpsi_\ep]- \big(\frac{1}{2}\I+{\BKcal}^*_{D_\ep}\big)[\bvarphi_\ep]
   =\pd{\H}{\nu}\q \mbox{on }\p D_\ep.  \\
  \end{array}
\right.
\end{equation}
Similarly, the solution to \eqref{equation-u} has the following representation
\begin{equation}\label{Representation-u}
\u(x)=
\left\{
  \begin{array}{lll}
   \ds  \H(x)+\BScal_{D}[\bvarphi](x), \q x\in  \RR^2 \backslash \overline{D}, \\
    \nm \ds\widetilde{ \BScal}_{D}[\bpsi](x),\q x\in D,  \\
  \end{array}
\right.
\end{equation}
where the pair $(\bpsi,\bvarphi)$ is the unique solution in $L^2(\p D)\times L^2_{\Psi}(\p D)$ of
\begin{equation}\label{phi-psi}
\left\{
  \begin{array}{lll}
   \ds  \widetilde{\BScal}_{D}[\bpsi]\big |_{-}- \BScal_{D}[\bvarphi]\big |_{+}
   =\H\q \mbox{on }\p D, \\
    \nm \ds \big(-\frac{1}{2}\I+\widetilde{\BKcal}^*_{D}\big)[\bpsi]- \big(\frac{1}{2}\I+{\BKcal}^*_{D}\big)[\bvarphi]
   =\pd{\H}{\nu}\q \mbox{on }\p D.  \\
  \end{array}
\right.
\end{equation}
Let $\Om$ be a bounded
region outside the inclusion $D$, and away from $\p D$. It then follows from \eqref{Representation-u-ep} and \eqref{Representation-u} that
\begin{align}\label{main-equation}
\ds \u_{\ep}(x)-\u(x)= \BScal_{D_\ep}[\bvarphi_\ep](x)-\BScal_{D}[\bvarphi](x),\q x\in \Om.
\end{align}
In order to prove the  asymptotic expansion for $\displaystyle (\u_\ep-\u)|_{ \Om}$ as $\ep$ tends to $0$, we next
investigate the asymptotic behavior of $\BScal_{D_\ep}[\bvarphi_\ep]$ as $\ep \rightarrow 0.$
%%%%%%%%%%%%%%%%%%%%%%%%%%%%ù
%%%%%%%%%%%%%%%%%%%%%%%%%%%%%ù
%%%%%%%%%%%%%%%%%%%%%%%%%%%%%%%ù
%%%%%%%%%%%%%%%%%%%%%%%%%%%%%%%%%%ù
\subsection{Proof of the theorem \ref{Main-theorem}}\label{proof}
For $\xe=x+\ep h(x)\n(x)\in \p D_\ep$. We have the following Taylor expansion
\begin{align}\label{H-ep}
\ds \H\big(x+\ep h(x)\n(x)\big)&=\H(x)+\ep h(x) \pd{\H}{\n}(x)+O(\ep^2), \q x\in \p D,
\end{align}
where the remainder $O(\ep^2)$ depends only on $\|h\|_{\mathcal{C}^0(\p D)}$ and $\|X\|_{\mathcal{C}^1(\p D)}$.

Similarly, by the Taylor expansion, \eqref{n0n1}, and  \eqref{local-Lame}, we obtain that
\begin{align}\label{conormal-H-ep}
\ds \pd{\H}{\nu}(\xe)=&\lambda_0 \nabla \cdot \H(\xe){\n}(\xe)+\mu_0\Big(\nabla \H(\xe)+(\nabla \H)^{T}(\xe)\Big){\n}(\xe)\nonumber\\
\ds =&\lambda_0 \nabla \cdot \H(x)\n(x)+\mu_0\Big(\nabla \H(x)+(\nabla \H)^{T}(x)\Big)\n(x)\nonumber\\
\ds &+\ep h(x)\Big[\lambda_0 \nabla \nabla \cdot \H(x)\cdot \n(x)+\mu_0 \nabla\nabla \H(x)\n(x)+\mu_0 \nabla(\nabla \H)^T(x)\n(x)\Big]\n(x)\nonumber\\
 \ds &-\ep h'(t)\Big[\lambda_0 \nabla \cdot \H(x)\ta(x)+\mu_0\Big(\nabla \H(x)+(\nabla \H)^{T}(x)\Big)\ta (x)\Big]+O(\ep^2)\nonumber\\
\ds =&\pd{\H}{\nu}(x)+ \ep \kappa(x)h(x)\pd{\H}{\nu}(x)-\ep\frac{d}{dt}\Big(h(x)\big(\mathbb{C}_0\widehat{\nabla} \H\big)\ta (x)\Big)+O(\ep^2), \q x\in \p D.
\end{align}

Now, we introduce  $(\bpsi^{(1)}, \bvarphi^{(1)})$ as a solution to  the following system
\begin{equation}\label{phi-psi-ep-system}
\left\{
  \begin{array}{lll}
   \ds  \widetilde{\BScal}_{D}[\bpsi^{(1)}]\big|_{-}- \BScal_{D}[\bvarphi^{(1)}]\big|_{+}
   =h\pd{\H}{ \n}-\big(\widetilde{\BScal}_D^{(1)}[\bpsi]-\BScal_D^{(1)}[\bvarphi]\big)\q \mbox{on }\p D, \\
    \nm \ds \Big(-\frac{1}{2}\I+\widetilde{\BKcal}^*_{D}\Big)[\bpsi^{(1)}]- \Big(\frac{1}{2}\I+{\BKcal}^*_{D}\Big)[\bvarphi^{(1)}]
   =\kappa h \pd{ \H}{ \nu}-\frac{\p}{\p \ta}\Big(h\big(\mathbb{C}_0 \widehat{\nabla} \H \big)\ta\Big)\\
   \nm\ds \q\q\q\q\q\q\q\q\q\q\q\q\q \q\q\q\q\q\q-\big(\widetilde{\BKcal}_D^{(1)}[\bpsi]-\BKcal_D^{(1)}[\bvarphi]\big)\q \mbox{on }\p D,
  \end{array}
\right.
\end{equation}
 where $(\bpsi,\bvarphi)$ is the solution to \eqref{phi-psi}. One can easily check the existence and uniqueness of $(\bpsi^{(1)}, \bvarphi^{(1)})$  by using the theorem \ref{ES-theorem}.

It follows from \eqref{phi-psi-ep}, \eqref{phi-psi-ep-system},  and Theorems \ref{theorem-KD-1} and \ref{theorem-SD-1} that
\begin{equation}\label{phi-ep-psi-ep}
\left\{
  \begin{array}{lll}
   \ds  \Big(\widetilde{\BScal}_D+\ep \widetilde{\BScal}^{(1)}_{D}\Big)\big [\widetilde{\bpsi} -\bpsi-\ep \bpsi^{(1)}\big ]\big|_{-}- \Big({\BScal}_D+\ep {\BScal}^{(1)}_{D}\Big)\big [\widetilde{\bvarphi} -\bvarphi-\ep \bvarphi^{(1)}\big ]\big|_{+}\\
  \nm\ds \q\q\q\q\q\q\q\q\q\q\q\q\q\q\q\q\q\q\q  =\H\circ \Phi_\ep-\H-\ep h \pd{\H}{\n}+O_1(\ep^2)\q \mbox{on }\p D, \\
    \nm \ds \Big(-\frac{1}{2}\I+\widetilde{\BKcal}_{D}^*+\ep \widetilde{\BKcal}_D^{(1)}\Big)\big [\widetilde{\bpsi} -\bpsi-\ep \bpsi^{(1)}\big ]- \Big(\frac{1}{2}\I+{\BKcal}_{D}^*+\ep {\BKcal}_D^{(1)}\Big)\big [\widetilde{\bvarphi} -\bvarphi-\ep \bvarphi^{(1)}\big ]\\
   \nm\ds \q\q\q\q\q\q\q\q   =\pd{\H}{\nu}\circ \Phi_\ep-\pd{ \H}{ \nu}-\ep \kappa h \pd{ \H}{ \nu}+\ep\frac{\p}{\p \ta}\Big(h\big(\mathbb{C}_0 \widehat{\nabla} \H \big)\ta\Big)+O_2(\ep^2)\q \mbox{on }\p D,
  \end{array}
\right.
\end{equation}
with  $\widetilde{\bvarphi}:=\bvarphi_\ep \circ \Phi_\ep$,   $\widetilde{\bpsi}:=\bpsi_\ep \circ \Phi_\ep$,  and $\|O_1(\ep^2)\|_{W_1^2(\p D)}, \|O_2(\ep^2)\|_{L^2(\p D)}\leq C\ep^2$, where the constant $C$ depends only on  $\lambda_0,\lambda_1,\mu_0,\mu_1$, the $\mathcal{C}^2$-norm of $X$, and the $\mathcal{C}^1$-norm of $h$.

The following lemma follows immediately from  \eqref{H-ep}, \eqref{conormal-H-ep}, \eqref{phi-ep-psi-ep},  and the estimate   in \eqref{estimate}.
\begin{lem}\label{Lemma} For $\ep$ small enough, there exists $C$ depending only on $\lambda_0,\lambda_1,\mu_0,\mu_1$, the $\mathcal{C}^2$-norm of $X$, and the $\mathcal{C}^1$-norm of $h$ such that
\begin{equation}\label{estimate-important}
\ds \Big\|\bpsi_\ep \circ \Phi_{\ep}-\bpsi-\ep \bpsi^{(1)} \Big\|_{L^2(\p D)}+\Big\| \bvarphi_\ep \circ \Phi_{\ep}-\bvarphi-\ep \bvarphi^{(1)} \Big\|_{L^2(\p D)}
\leq C \ep^2,
\end{equation}
where $(\bpsi_\ep,\bvarphi_\ep),(\bpsi,\bvarphi),$ and $(\bpsi^{(1)},\bvarphi^{(1)})$  are the solutions to \eqref{phi-psi-ep}, \eqref{phi-psi}, and  \eqref{phi-psi-ep-system}, respectively.
\end{lem}
%%%%%%%%%%%%%%%%%%%%%%%%%%%%%%%
%%%%%%%%%%%%%%%%%%%%%%%%%%%%%%%%%%%%%%ù
%%%%%%%%%%%%%%%%%%%%%%%%%%%%%%%%%%%%%
%%%%%%%%%%%%%%%%%%%%%%%%%%%%%%%%%%%%%%%%%
%%%%%%%%%%%%%%%%%%%%%%%%%%%%%%%%%%%%%%%%%%%ùù
%%%%%%%%%%%%%%%%%%%%%%%%%%%%%%ù
%%%%%%%%%%%%%%%%%%%%%%%%%%%%%%%%
%%%%%%%%%%%%%%%%%%%%%%%%%%%%%ù
%%%%%%%%%%%%%%%%%%%%%%%%%ùù

Recall that the domain $D$ is separated apart from  $ \Om$, then
\begin{align*}
\ds \sup_{x\in \Om, y\in \p D}\Big|\p^i \G(x-y)\Big|\leq C,\q i\in \NN^2,
\end{align*}
for some constant $C>0$ depending on $dist(D,\Om)$.  After the change of variables $\ye=\Phi_\ep(y)$,
we get from \eqref{sigexp}, \eqref{estimate-important}, and the Taylor
expansion of $\G(x-\ye)$ for $y\in \p D$, and $x\in  \Om$ fixed that
\begin{align}\label{main-equation-expansion}
\ds \BScal_{D_\ep}[\bvarphi_\ep](x)=&\int_{\p D_\ep}\G(x-\tilde y) \bvarphi_\ep(\tilde y)d\sigma(\tilde y)\nonumber\\
\nm\ds =&\int_{\p D}\Big(\G(x-y)+\ep h(y)\nabla \G(x-y)\n(y)\Big)\Big( \bvarphi(y)+\ep \bvarphi^{(1)}(y)\Big)\nonumber\\
\nm\ds &\q\q\q\q\q\q\q\q\q\q\q\q\q\q \q\q \q\q \times \Big(1-\ep\kappa(y)h(y)\Big)d\sigma (y) +O(\ep^2)\nonumber\\
\nm\ds =& \BScal_{D}[\bvarphi](x)+\ep\Big(\BScal_{D}[\bvarphi^{(1)}](x)-\BScal_{D}[\kappa h\bvarphi](x)+\BDcal^{\sharp}_{D}[h\bvarphi](x)\Big)+O(\ep^2).
\end{align}
The following theorem  follows immediately from \eqref{main-equation} and \eqref{main-equation-expansion}.
\begin{thm}\label{Thm-asymp-version0} Let $\ep$ be small enough. The following pointwise  expansion  holds  for $ x\in\Om$
\begin{align}\label{asymp-version0}
\ds \u_\ep(x)=\u(x)+\ep\Big(\BScal_{D}[\bvarphi^{(1)}](x)-\BScal_{D}[\kappa h\bvarphi](x)+\BDcal^{\sharp}_{D}[h\bvarphi](x)\Big)+O(\ep^2),
\end{align}
where $\bvarphi$ and $\bvarphi^{(1)}$ are defined  by \eqref{phi-psi} and \eqref{phi-psi-ep-system}, respectively.
The remainder $O(\ep^2)$ depends only  on $\lambda_0,\lambda_1,\mu_0,\mu_1$, the $\mathcal{C}^2$-norm of $X$,
the $\mathcal{C}^1$-norm of $h$, and $dist(\Om,D).$
\end{thm}

We now prove the following  representation theorem  for the solution
of the transmission problem \eqref{equation-u-1} which will be very
helpful in the proof of  theorem \ref{Main-theorem}.
\begin{thm} \label{Representation-u1} The solution $\u_1$ of \eqref{equation-u-1} is represented by
\begin{equation}\label{u1-second}
\u_1(x)=
\left\{
  \begin{array}{lll}
   \ds  \BScal_{D}[\bvarphi^{(1)}](x)-\BScal_{D}[\kappa h\bvarphi ](x)+\BDcal^\sharp_{D}[h\bvarphi](x),\q  x\in\RR^2 \backslash \overline{D}, \\
    \nm \ds   \widetilde{\BScal}_{D}[\bpsi^{(1)}](x)-\widetilde{\BScal}_{D}[\kappa h\bpsi](x)+\widetilde{\BDcal}^\sharp_{D}[h\bpsi](x),\q  x\in D,  \\
  \end{array}
\right.
\end{equation}
where  $(\bpsi, \bvarphi)$ and $(\bpsi^{(1)},\bvarphi^{(1)}) $ are defined by  \eqref{phi-psi} and \eqref{phi-psi-ep-system}, respectively.
\end{thm}
\proof  One can easily see that
\begin{align*}
\nm \ds \Lcal_{\lambda_0,\mu_0}\u_1=0 \q \mbox{in }\RR^2\backslash \overline{D},\q\q\q\q \Lcal_{\lambda_1,\mu_1}\u_1=0\q \mbox{in }D.
\end{align*}
It follows from \eqref{identity-3}, \eqref{SD-1}, \eqref{Representation-u}, and \eqref{phi-psi-ep-system} that
\begin{align*}
\ds \ds \u_1^i-\u_1^e=&\widetilde{\BScal}_D[\bpsi^{(1)}]-\BScal_D[\bvarphi^{(1)}]+\BScal_{D}[\kappa h\bvarphi ]-\widetilde{\BScal}_{D}[\kappa h\bpsi ]+\widetilde{\BDcal}^\sharp_{D}[h\bpsi]\big|_{-}-\BDcal^\sharp_{D}[h\bvarphi]\big|_{+}\\
\nm\ds =&h \pd{\H}{\n}+\BScal_D^{(1)}[\bvarphi]-\widetilde{\BScal}_D^{(1)}[\bpsi]+\BScal_{D}[\kappa h\bvarphi ]-\widetilde{\BScal}_{D}[\kappa h\bpsi ]+\widetilde{\BDcal}^\sharp_{D}[h\bpsi]\big|_{-}-\BDcal^\sharp_{D}[h\bvarphi]\big|_{+}\\
\nm\ds =& h\Big( \pd{\H}{\n}+\pd{\BScal_D[\bvarphi]}{\n}\Big|_{+}-\pd{\widetilde{\BScal}_D[\bpsi]}{\n}\Big|_{-}\Big)\\
\nm\ds =& h\big(\nabla \u^e\n-\nabla\u^i \n\big)\\
\nm\ds =&h \big(\mathbb{K}_{0,1}\widehat{\nabla} \u^i\big)\n \q \mbox{on } \p D.
\end{align*}
Using \eqref{phi-psi-ep-system}, we get
\begin{align*}
\ds \ds \pd{\u_1}{\widetilde{\nu}}\Big|_{-}-\pd{\u_1}{\nu}\Big|_{+}=&\pd{\widetilde{\BScal}_D[\bpsi^{(1)}]}{\widetilde{\nu}}\Big|_{-}-
\pd{\BScal_D[\bvarphi^{(1)}]}{\nu}\Big|_{+} +\pd{\BScal_{D}[\kappa h\bvarphi ]}{\nu}\Big|_{+}-\pd{\widetilde{\BScal}_{D}[\kappa h\bpsi ]}{\widetilde{\nu}}\Big|_{-}\\
\nm\ds&+\pd{\widetilde{\BDcal}^\sharp_{D}[h\bpsi]}{\widetilde{\nu}}\Big|_{-}-\pd{\BDcal^\sharp_{D}[h\bvarphi]}{\nu}\Big|_{+}\\
\nm\ds =&\kappa h \pd{\H}{\nu}-\frac{\p}{\p \ta}\Big(h\big(\mathbb{C}_0\widehat{\nabla}\H)\big)\ta\Big)-\widetilde{\BKcal}_D^{(1)}[\bpsi]+\BKcal_D^{(1)}[\bvarphi]\\
\nm\ds &+\pd{\BScal_{D}[\tau h\bvarphi ]}{\nu}\Big|_{+}-\pd{\widetilde{\BScal}_{D}[\kappa h\bpsi ]}{\widetilde{\nu}}\Big|_{-}+\pd{\widetilde{\BDcal}^\sharp_{D}[h\bpsi]}{\widetilde{\nu}}\Big|_{-}-\pd{\BDcal^\sharp_{D}[h\bvarphi]}{\nu}\Big|_{+}.
\end{align*}
According to \eqref{Identity-1}, \eqref{KD-1-final}, \eqref{Representation-u}, and \eqref{conormal-H-ep} we obtain
\begin{align*}
\ds  \pd{\u_1}{\widetilde\nu}\Big|_{-}-\pd{\u_1}{{\nu}}\Big|_{+}=&\frac{\p}{\p \ta}\big(h (\mathbb{C}_1\widehat{\nabla} \u^i )\ta\big) -\frac{\p}{\p \ta}\big(h (\mathbb{C}_0\widehat{\nabla} \u^e )\ta\big)=\frac{\p}{\p \ta}\big(h [\mathbb{C}_1-\mathbb{M}_{0,1}]\widehat{\nabla} \u^i )\ta\big).
\end{align*}
Now,  let us check the  condition
\begin{align}\label{phi1}
\ds \BScal_{D}[\bvarphi^{(1)}-\kappa h\bvarphi ](x)\rightarrow 0\q \mbox{ as } |x|\rightarrow \infty.
\end{align}
To do this,  we rewrite the system of equations \eqref{phi-psi-ep-system}
\begin{align}\label{phi-psi-ep-system-second}
\left\{
  \begin{array}{lll}
   \ds  \widetilde{\BScal}_{D}\big[\bpsi^{(1)}-\kappa h \bpsi\big]\Big|_{-}- \BScal_{D}\big[\bvarphi^{(1)}-\kappa h \bvarphi\big]\Big|_{+}
= h \big(\mathbb{K}_{0,1}\widehat{\nabla} \u^i\big)\n-\widetilde{\BDcal}_{D}^{\sharp}[h\bpsi]\big|_{-}+\BDcal_{D}^{\sharp}[h\bvarphi]\big|_{+}\\
    \nm \ds \pd{\widetilde{\BScal}_D}{\tilde \nu}\big[\bpsi^{(1)}-\kappa h \bpsi\big]\Big |_{-}- \pd{{\BScal}_D}{ \nu}\big[\bvarphi^{(1)}-\kappa h \bvarphi\big]\Big|_{+}=\frac{\p\BDcal_{D}^{\sharp}[h\bvarphi]}{\p \nu}\Big|_{+}-\frac{\p\widetilde{\BDcal}_{D}^{\sharp}[h\bpsi]}{\p \tilde \nu}\Big|_{-}\\
\nm\ds \q\q\q\q\q\q\q\q\q\q\q \q \q \q \q \q \q\q \q +\frac{\p}{\p\ta}\Big(h\big(\mathbb{C}_1 \widehat{\nabla}\u^i\big)\ta\Big) -\frac{\p}{\p \ta}\Big(h\big(\mathbb{C}_0 \widehat{\nabla}\u^e\big)\ta\Big).
  \end{array}
\right.
\end{align}
  It is clear that
\begin{align*}
\ds \int_{\p D} \Big[\frac{\p}{\p\ta}\Big(h\big(\mathbb{C}_1 \widehat{\nabla}\u^i\big)\ta\Big) -\frac{\p}{\p \ta}\Big(h\big(\mathbb{C}_0 \widehat{\nabla}\u^e\big)\ta\Big)\Big]\cdot \theta_m ~d\sigma =0\q \mbox{for } m=1,2.
\end{align*}
We have
\begin{align*}
\ds \int_{\p D} \frac{d}{dt}\Big( h\big(\mathbb{C}_0 \widehat{\nabla}\u^e(x)\big)\ta(x)\Big) \cdot \theta_3(x) d\sigma& =-
\int_{\p D}h(x) \big(\mathbb{C}_0 \widehat{\nabla}\u^e(x)\big)\ta(x)\cdot \n(x) d\sigma\\
\nm\ds &=-\mu_0\int_{\p D}h(x) \Big(\nabla \u^e(x)+(\nabla \u^e)^{T}(x) \Big)\ta(x)\cdot \n (x) d\sigma\\
\nm\ds &=-\mu_0 \int_{\p D} h(x)\Big(\nabla \u^e(x)+(\nabla \u^e)^{T}(x) \Big)\n(x)\cdot \ta(x) d\sigma\\
\nm\ds &=-\int_{\p D} h(x)\pd{\u^e}{\nu}(x)\cdot \ta (x)d\sigma.
\end{align*}
Similarly, we get
\begin{align*}
\ds \int_{\p D} \frac{d}{dt}\Big( h(x)\big(\mathbb{C}_1 \widehat{\nabla}\u^i(x)\big)\ta(x) \Big) \cdot \theta_3(x) d\sigma
=-\int_{\p D} h(x)\pd{\u^i}{\tilde \nu}(x)\cdot \ta(x)d\sigma.
\end{align*}
Thus
\begin{align*}
\ds \int_{\p D} \Big[\frac{\p }{\p \ta }\Big(h\big(\mathbb{C}_1 \widehat{\nabla}\u^i\big)\ta \Big) -\frac{\p}{\p \ta }\Big(h\big(\mathbb{C}_0 \widehat{\nabla}\u^e\big)\ta\Big)\Big]\cdot \theta_3d\sigma =
\int_{\p D} h\Big(\pd{\u^e}{\nu}-\pd{\u^i}{\tilde \nu}\Big)\cdot \ta d\sigma=0.
\end{align*}
Consequently,
\begin{align*}
\ds \frac{\p}{\p\ta}\Big(h\big(\mathbb{C}_1 \widehat{\nabla}\u^i\big)\ta\Big) -\frac{\p}{\p \ta}\Big(h\big(\mathbb{C}_0 \widehat{\nabla}\u^e\big)\ta \Big) \in L^2_\Psi(\p D).
\end{align*}
By \eqref{L-Psi},  $\p\widetilde{\BScal}_D[\bpsi^{(1)}-\kappa h \bpsi]/{\p{\tilde \nu }}\big |_{-}$ and $\p\widetilde{\BDcal}_{D}^{\sharp}[h\bpsi]/\p \tilde \nu\big |_{-} \in L^2_\Psi(\p D)$. It then follows from \eqref{phi-psi-ep-system-second} that $\p \BScal_D[\bvarphi^{(1)}-\kappa h \bvarphi]/\p{ \nu}\big|_{+}+ \p\BDcal_{D}^{\sharp}[h\bvarphi]/\p \nu\big|_{+} \in L^2_\Psi(\p D)$. Since
\begin{align*}
\ds \bvarphi^{(1)}-\kappa h \bvarphi+\frac{\p }{\p \ta}\Big(h\la \bvarphi, \ta\ra \n +\frac{\lambda_0}{2\mu_0+\lambda_0}h \la \bvarphi, \n \ra \ta\Big)=&
\pd{{\BScal}_D}{ \nu}\big[\bvarphi^{(1)}-\kappa h \bvarphi\big]\Big|_{+}+\frac{\p\BDcal_{D}^{\sharp}[h\bvarphi]}{\p \nu}\Big|_{+}\\
\nm \ds &-\pd{{\BScal}_D}{ \nu}\big[\bvarphi^{(1)}-\kappa h \bvarphi\big]\Big|_{-}-\frac{\p\BDcal_{D}^{\sharp}[h\bvarphi]}{\p \nu}\Big|_{-},
\end{align*}
with $\p \BScal_D[\bvarphi^{(1)}-\kappa h \bvarphi]/\p{ \nu}\big|_{-}, \p\BDcal_{D}^{\sharp}[h\bvarphi]/\p \nu\big|_{-} \in L^2_\Psi(\p D)$, see \eqref{L-Psi}. Then
\begin{align*}
\ds \bvarphi^{(1)}-\kappa h \bvarphi+\frac{\p }{\p \ta}\Big(h\la \bvarphi, \ta\ra \n +\frac{\lambda_0}{2\mu_0+\lambda_0}h\la \bvarphi, \n \ra \ta\Big)\in L^2_\Psi(\p D).
\end{align*}
Therefore, we have
\begin{align*}
\ds \BScal_{D}[\bvarphi^{(1)}-\kappa h\bvarphi ](x)&=\G(x)\int_{\p D}(\bvarphi^{(1)}-\kappa h\bvarphi) d\sigma+O(|x|^{-1})\\
\nm \ds &=\G(x)\int_{\p D}\bigg(\bvarphi^{(1)}-\kappa h\bvarphi+\frac{\p }{\p \ta}\Big(h \la \bvarphi, \ta\ra \n +\frac{\lambda_0}{2\mu_0+\lambda_0}h\la \bvarphi, \n \ra \ta\Big)\bigg) d\sigma\\
\nm \ds&\q-\G(x)\int_{\p D}\frac{\p }{\p \ta}\Big(h\la \bvarphi, \ta\ra \n +\frac{\lambda_0}{2\mu_0+\lambda_0}h \la \bvarphi, \n \ra \ta\Big) d\sigma+O(|x|^{-1})\\
\nm\ds &=O(|x|^{-1})\q \mbox{ as }|x|\rightarrow \infty.
\end{align*}
Thus   $\u_1$ defined by   \eqref{u1-second} satisfies $\u_1(x)=O(|x|^{-1})$  as $|x|\rightarrow \infty.$ This completes the proof of the theorem \ref{Representation-u1}.

The  theorem \ref{Main-theorem}  immediately follows from  the integral representation of $\u_1$   in  \eqref{u1-second} and the theorem \ref{Thm-asymp-version0}.
%%%%%%%%%%%%%%%%%%%%%%%%%%%%%%%%%%%%%%%%%%
%%%%%%%%%%%%%%%%%%%%%%%%%%%%%%%%%%%%%%%
%%%%%%%%%%%%%%%%%%%%%%%%%%%%%%%%%%%%%%%
%%%%%%%%%%%%%%%%%%%%%%%%%%%%%%%%%%%%%%%%%%
%%%%%%%%%%%%%%%%%%%%%%%%%%%%%%%%%%%%%
%%%%%%%%%%%%%%%%%%%%%%%%%%%%%%%%%%%%%%%%%%
%%%%%%%%%%%%%%%%%%%%%%%%%%%%%%%%%%%%%%%
\subsection{Proof of the theorem \ref{second-theorem}}
The following corollary can be proved in exactly the same manner as
Theorem \ref{Main-theorem}.
\begin{cor}Let  $\u$ and $\u_\ep$ be the solutions to
\eqref{equation-u} and \eqref{equation-u-ep}, respectively. Let $\Om$ be a bounded
region outside the inclusion $D$, and away from $\p D$. For $x \in \Om$, the following pointwise asymptotic expansion holds:
\begin{equation}\label{Main-Asymptotic-normal}
\ds \pd{\u_\ep}{\nu}(x)=\pd{\u}{\nu}(x)+\ep \pd{\u_1}{\nu}(x)+O(\ep^2),
\end{equation}
where the remainder $O(\ep^2)$ depends only on $\lambda_0,\lambda_1,\mu_0,\mu_1$,
the $\mathcal{C}^2$-norm of $X$, the $\mathcal{C}^1$-norm of $h$,  $dist (\Om, \p D)$,
and  $\u_1$ is the unique solution of \eqref{equation-u-1}.
\end{cor}
 Let $S$ be a Lipschitz closed curve enclosing $D$  away from $\p D$. Let $\v$ be the solution to \eqref{v}. It follows from \eqref{Main-Asymptotic}, \eqref{Important-relation-2}, and \eqref{Main-Asymptotic-normal} that
\begin{align*}
\ds\int_{S}\big(\u_\ep-\u\big)\cdot\pd{ \F}{\nu}d\sigma-\int_{S}\big(\pd{\u_\ep}{\nu}-
\pd{ \u}{\nu}\big) \cdot\F d\sigma& =\ep \int_{S}\Big(\u_1\cdot \pd{ \v}{\nu}-\pd{\u_1}{\nu}
\cdot \v\Big) d\sigma+O(\ep^2).
\end{align*}
By using Lemma \ref{relations} to  the integral on the right-hand side, we get
$$
\int_{S}\Big( \pd{ \v}{\nu}\cdot\u_1-\v\cdot\pd{\u_1}{\nu}
\Big) d\sigma=
 \int_{\p D} \Big(\pd{\v^{e}}{\nu}\cdot \u^{e}_1-\v^e \cdot\pd{\u^{e}_1}{\nu}\Big)d\sigma.
$$
According to the jump conditions for $\u_1$ in \eqref{equation-u-1}, we deduce that
\begin{align}\label{eq001}
\ds\int_{S}\Big( \pd{ \v}{\nu}\cdot\u_1-\v\cdot\pd{\u_1}{\nu}
\Big) d\sigma=& \int_{\p D} \Big(\pd{\v^{i}}{\tilde{\nu}}\cdot \u^{i}_1-\v^i \cdot \pd{\u^{i}_1}{\tilde{\nu}}\Big)d\sigma\nonumber\\
\nm\ds &-\int_{\p D} h\big(\mathbb{K}_{0,1}\widehat{\nabla}\u^{i}\big)\n \cdot \big(\mathbb{C}_1\widehat{\nabla}\v^{i}\big)\n d\sigma\nonumber\\
\nm \ds &+ \int_{\p D} \frac{\p}{\p \ta }\Big(h\big([\mathbb{C}_1-\mathbb{M}_{0,1}]\widehat{\nabla}\u^{i}\big)\ta \Big) \cdot \v^{i} d\sigma .
\end{align}
It follows from \eqref{Important-relation} that
\begin{align}\label{eq002}
\int_{\p D} \Big(\pd{\v^{i}}{\tilde{\nu}}\cdot \u^{i}_1-\v^i \cdot \pd{\u^{i}_1}{\tilde{\nu}}\Big)d\sigma=0.
\end{align}
We have
\begin{align}\label{eq003}
 \ds \int_{\p D} \frac{\p}{\p \ta } \Big(h\big([\mathbb{C}_1-\mathbb{M}_{0,1}]\widehat{\nabla}\u^{i}\big)\ta \Big) \cdot \v^{i} d\sigma &=
-\int_{\p D}  h\big([\mathbb{C}_1-\mathbb{M}_{0,1}]\widehat{\nabla}\u^{i}\big)\ta  \cdot \nabla \v^{i} \ta d\sigma.
\end{align}
One can easily check that
\begin{align}\label{eq004}
 \ds \big([\mathbb{C}_1-\mathbb{M}_{0,1}]\widehat{\nabla}\u^{i}\big)\ta  \cdot \nabla \v^{i} \ta
=\big([\mathbb{C}_1-\mathbb{M}_{0,1}]\widehat{\nabla}\u^{i}\big)\ta  \cdot \widehat{\nabla} \v^{i} \ta.
\end{align}
We finally obtain  from \eqref{eq001}-\eqref{eq004} the  relationship between traction-displacement measurements and the
shape deformation $h$ \eqref{asymptotic-traction-displacement}.
%%%%%%%%%%%%%%%%%%%%%%%%%%%%%%%%%%%%%%%
%%%%%%%%%%%%%%%%%%%%%%%%%%%%%%%%%%%%%%%%%%
%%%%%%%%%%%%%%%%%%%%%%%%%%%%%%%%%%%%%
%%%%%%%%%%%%%%%%%%%%%%%%%%%%%%%%%%%%%%%%%%
%%%%%%%%%%%%%%%%%%%%%%%%%%%%%%%%%%%%%%%
%%%%%%%%%%%%%%%%%%%%%%%%%%%%%%%%%%%%%%%
%%%%%%%%%%%%%%%%%%%%%%%%%%%%%%%%%%%%%%%%%%
%%%%%%%%%%%%%%%%%%%%%%%%%%%%%%%%%%%%%
\section{Asymptotic expansion of EMTs}
We  introduce the notion of EMTs associated with $D$ and Lam\'e parameters $(\lambda_0, \mu_0)$  for the background and $(\lambda_1, \mu_1)$ for $D$
as follows (see \cite{book, book2}): For multi-index $\alpha \in \NN^2$ and $j=1,2$, let the pair $(\f_\alpha^j, \g_\alpha^j)$ in $L^2(\p D)\times L^2(\p D)$ be the unique solution to
\begin{equation}\label{phi-psi-EMTs}
\left\{
  \begin{array}{lll}
 \ds  \widetilde{\BScal}_{D}[\f_\alpha^j]\big |_{-}- \BScal_{D}[\g_\alpha^j]\big |_{+}
 =x^\alpha \e_j\q& \mbox{ on } \p D, \\
\nm \ds \pd{\widetilde{\BScal}_{D}[\f_\alpha^j]}{\widetilde{\nu}}\Big |_{-}- \pd{\BScal_{D}[\g_\alpha^j]}{\nu}\Big |_{+}
=\pd{(x^\alpha \e_j)}{\nu} \q& \mbox{ on }  \p D.   \\
  \end{array}
\right.
\end{equation}
Now for multi-index $\beta \in \NN^2$, the EMTs are defined by
\begin{align}\label{Definition-EMTs}
M_{\alpha \beta}^j=(m^j_{\alpha\beta 1}, m^j_{\alpha\beta 2}):= \int_{\p D} y^\beta \g_\alpha^j(y) d\sigma(y).
\end{align}
Let $\displaystyle \H(x)=\sum_{j=1}^2\sum_{\alpha \in \NN^2}a^\alpha_j x^\alpha e_j$ and $\displaystyle \F(x)=\sum_{k=1}^2\sum_{\beta \in \NN^2}b^\beta_k x^\beta e_k$ be tow polynomials satisfying $\nabla \cdot \big(\mathbb{C}_0 \widehat{\nabla} \cdot \big)=0$
in $\RR^2$. The EMTs $m_{\alpha \beta k}^j(D)$ associated with $D$ satisfy
\begin{align}\label{EMT}
\sum_{\alpha \beta j k}a_{j}^{\alpha} b_{k}^{\beta}m_{\alpha \beta k}^j(D) = \int_{\p D}\F(y)\bvarphi(y)d\sigma(y),
\end{align}
where $\bvarphi$  is defined in \eqref{phi-psi}.

The perturbed $m_{\alpha \beta k}^j(D_\ep)$ satisfy
\begin{align}\label{EMT-ep}
\ds \sum_{\alpha \beta j k}a_{j}^{\alpha} b_{k}^{\beta}m_{\alpha \beta k}^j(D_\ep) = \int_{\p D_\ep}\F(\ye){\bvarphi}_\ep(\ye)d\sigma_\ep(\ye),
\end{align}
where ${\bvarphi}_\ep$ is defined in \eqref{phi-psi-ep}.\\

The purpose of this  section is  to prove the  asymptotic behavior of  $ \sum_{\alpha \beta j k}a_{j}^{\alpha} b_{k}^{\beta}m_{\alpha \beta k}^j(D_\ep)$  defined in \eqref{EMT-ep} as $\ep$ tends to zero.

By Taylor expansion, we have
\begin{align*}
 \F(\ye)=\F\big(y+\ep h(y) \n(y)\big)=\F(y)+\ep h(y)\pd{\F}{\n}(y)+O(\ep^2),\q y\in \p D.
\end{align*}
It follows from Lemma \ref{Lemma} that
\begin{align*}
 {\bvarphi}_\ep(\ye)=\bvarphi_\ep\big(y+\ep h(y) \n(y)\big)=\bvarphi(y)+\ep \bvarphi^{(1)}(y)+O(\ep^2),\q y\in \p D,
\end{align*}
where  $\bvarphi^{(1)}$ is defined in  \eqref{phi-psi-ep-system}.

Recall that $d\sigma_\ep(\ye)=\big(1-\ep \kappa h(y)\big)d\sigma(y)+O(\ep^2)$   for $y\in \p D$. After the
change of variables $\ye=y+\ep h(y) \n(y)$, we get from  \eqref{EMT-ep} that
\begin{align}
 \ds \sum_{\alpha \beta j k}a_{j}^{\alpha} b_{k}^{\beta}m_{\alpha \beta k}^j(D_\ep) &= \int_{\p D}\Big(\F+\ep h\pd{\F}{\n}\Big)\cdot\Big(\bvarphi+\ep \bvarphi^{(1)}\Big)\big(1-\ep \kappa h\big)d\sigma+O(\ep^2)\nonumber\\
 \nm \ds &=\sum_{\alpha \beta j k}a_{j}^{\alpha} b_{k}^{\beta}m_{\alpha \beta k}^j(D) +\ep \int_{\p D}\F\cdot\Big(\bvarphi^{(1)}-\kappa h\bvarphi\Big)d\sigma\nonumber\\
\nm \ds &\q+\ep \int_{\p D}h\pd{\F}{\n}\cdot \bvarphi d\sigma+O(\ep^2). \label{eqEMT-1}
\end{align}
From \eqref{nuS}, we have
\begin{align*}
\int_{\p D}\F\cdot\Big(\bvarphi^{(1)}-\kappa h\bvarphi\Big)d\sigma=\int_{\p D}\F\cdot\bigg(\pd{\BScal_{D}[\bvarphi^{(1)}-\kappa h\bvarphi]}{\nu}\Big |_{+}-\pd{\BScal_{D}[\bvarphi^{(1)}-\kappa h\bvarphi]}{\nu}\Big |_{-}\bigg)d\sigma.
\end{align*}
By using \eqref{Important-relation-2} and \eqref{phi1}, we get
\begin{align*}
\ds \int_{\p D}\F\cdot\pd{\BScal_{D}[\bvarphi^{(1)}-\kappa h\bvarphi]}{\nu}\Big |_{+}d\sigma&=\int_{\p D}(\F-\v^{e})\cdot\pd{\BScal_{D}[\bvarphi^{(1)}-\kappa h\bvarphi]}{\nu}\Big |_{+}d\sigma\\
\nm\ds &\q+\int_{\p D}\v^{e}\cdot\pd{\BScal_{D}[\bvarphi^{(1)}-\kappa h\bvarphi]}{\nu}\Big |_{+}d\sigma\\
\nm \ds &=\int_{\p D}\Big(\pd{\F}{\nu}-\pd{\v^{e}}{\nu}\Big)\cdot \BScal_{D}[\bvarphi^{(1)}-\kappa h\bvarphi]d\sigma\\
\nm\ds &\q+\int_{\p D}\v^{e}\cdot\pd{\BScal_{D}[\bvarphi^{(1)}-\kappa h\bvarphi]}{\nu}\Big |_{+}d\sigma.
\end{align*}
Since, by using \eqref{Important-relation},  we get
\begin{align*}
\ds \int_{\p D}\pd{\F}{\nu}\cdot \BScal_{D}[\bvarphi^{(1)}-\kappa h\bvarphi]d\sigma-\int_{\p D}\F\cdot\pd{\BScal_{D}[\bvarphi^{(1)}-\kappa h\bvarphi]}{\nu}\Big |_{-}d\sigma=0.
\end{align*}
It  then follows from \eqref{v} and \eqref{phi-psi-ep-system-second} that
\begin{align}
\ds&\int_{\p D}\F\cdot\Big(\bvarphi^{(1)}-\kappa h\bvarphi\Big)d\sigma\nonumber\\
\nm \ds &=\int_{\p D}\v^{e}\cdot\pd{\BScal_{D}[\bvarphi^{(1)}-\kappa h\bvarphi]}{\nu}\Big |_{+}d\sigma
-\int_{\p D}\pd{\v^{e}}{\nu}\cdot \BScal_{D}[\bvarphi^{(1)}-\kappa h\bvarphi]d\sigma\nonumber\\
\nm \ds&=\int_{\p D}\v^{i}\cdot\pd{\widetilde{\BScal}_{D}[\bpsi^{(1)}-\kappa h\bpsi]}{\tilde{\nu}}\Big |_{-}d\sigma
-\int_{\p D}\pd{\v^{i}}{\tilde{\nu}}\cdot \widetilde{\BScal}_{D}[\bpsi^{(1)}-\kappa h\bpsi] d\sigma\nonumber\\
\nm \ds&\q+\int_{\p D}\v^{i}\cdot\frac{\p\widetilde{\BDcal}_{D}^{\sharp}[h\bpsi]}{\p \tilde \nu}\Big|_{-}d\sigma
-\int_{\p D}\pd{\v^{i}}{\tilde{\nu}}\cdot\widetilde{\BDcal}_{D}^{\sharp}[h\bpsi]\big|_{-} d\sigma\nonumber\\
\nm \ds &\q -\int_{\p D}\v^{e}\cdot\frac{\p \BDcal_{D}^{\sharp}[h\bvarphi]}{\p  \nu}\Big|_{+}d\sigma
+\int_{\p D}\pd{\v^{e}}{\nu}\cdot\BDcal_{D}^{\sharp}[h\bvarphi]\big|_{+} d\sigma\nonumber\\
\nm\ds &\q-\int_{\p D}\v^{i}\cdot\frac{\p}{\p \ta}\Big(h\big([\mathbb{C}_1-\mathbb{M}_{0,1}] \widehat{\nabla}\u^i\big)\ta\Big)d\sigma+\int_{\p D}\pd{\v^{i}}{\tilde\nu}\cdot\Big( h (\mathbb{K}_{0,1} \widehat{\nabla}\u^i\big)\n\Big)d\sigma.\label{eqEMT-2}
\end{align}
One can  easily  see  that
\begin{align}
\ds &-\int_{\p D}\v^i\cdot\frac{\p }{\p \ta }\Big(h\big([\mathbb{C}_1-\mathbb{M}_{0,1}] \widehat{\nabla}\u^i\big)\ta\Big)d\sigma+\int_{\p D}\pd{\v^i}{\tilde\nu}\cdot\Big( h (\mathbb{K}_{0,1} \widehat{\nabla}\u^i\big)\n\Big)d\sigma\nonumber\\
\nm \ds & =\int_{\p D}h\bigg(\big([\mathbb{C}_1-\mathbb{M}_{0,1}] \widehat{\nabla}\u^i\big)\ta \cdot \nabla \v^i \ta+(\mathbb{K}_{0,1} \widehat{\nabla}\u^i\big)\n\cdot \big(\mathbb{C}_1 \widehat{\nabla} \v^i \big)\n \bigg)d\sigma.\label{eqEMT-3}
\end{align}
We now apply \eqref{Important-relation} to obtain that
\begin{align}
\ds& \int_{\p D}\v^{i}\cdot\pd{\widetilde{\BScal}_{D}[\bpsi^{(1)}-\kappa h\bpsi]}{\tilde{\nu}}\Big |_{-}d\sigma
-\int_{\p D}\pd{\v^{i}}{\tilde{\nu}}\cdot \widetilde{\BScal}_{D}[\bpsi^{(1)}-\kappa h\bpsi] d\sigma=0, \label{eqEMT-4}\\
\nm \ds&\int_{\p D}\v^{i}\cdot\frac{\p\widetilde{\BDcal}_{D}^{\sharp}[h\bpsi]}{\p \tilde \nu}\Big|_{-}d\sigma
-\int_{\p D}\pd{\v^{i}}{\tilde{\nu}}\cdot\widetilde{\BDcal}_{D}^{\sharp}[h\bpsi]\big|_{-} d\sigma=0.\label{eqEMT-5}
\end{align}
It follows from \eqref{Important-relation}, \eqref{Important-relation-2}, \eqref{conormal-D-sharp-identity}, and the proposition \ref{D-sharp}  that
\begin{align*}
  \ds & \int_{\p D}\pd{\v^{e}}{\nu}\cdot\BDcal_{D}^{\sharp}[h\bvarphi]\big|_{+} d\sigma-\int_{\p D}\v^{e}\cdot\frac{\p \BDcal_{D}^{\sharp}[h\bvarphi]}{\p  \nu}\Big|_{+}d\sigma\\
\nm\ds & \q =\int_{\p D}\pd{\v^{e}}{\nu}\cdot\BDcal_{D}^{\sharp}[h\bvarphi]\big|_{+} d\sigma-\int_{\p D}(\v^{e}-\F)\cdot\frac{\p \BDcal_{D}^{\sharp}[h\bvarphi]}{\p  \nu}\Big|_{+}d\sigma-\int_{\p D} \F\cdot\frac{\p \BDcal_{D}^{\sharp}[h\bvarphi]}{\p  \nu}\Big|_{+}d\sigma\\
\nm \ds &\q =\int_{\p D}\pd{\F}{\nu}\cdot\BDcal_{D}^{\sharp}[h\bvarphi]\big|_{+} d\sigma-\int_{\p D} \F\cdot\frac{\p \BDcal_{D}^{\sharp}[h\bvarphi]}{\p  \nu}\Big|_{+}d\sigma\\
\nm \ds &\q =-\int_{\p D}\F \cdot \frac{\p }{\p \ta}\Big(\la h\bvarphi , \ta\ra \n +\frac{\lambda_0}{2\mu_0+\lambda_0}\la h\bvarphi, \n \ra \ta\Big)d\sigma \\
\nm \ds &\q\q +\int_{\p D}\pd{\F}{\nu}\cdot\BDcal_{D}^{\sharp}[h\bvarphi]\big|_{+} d\sigma-\int_{\p D} \F\cdot\frac{\p \BDcal_{D}^{\sharp}[h\bvarphi]}{\p  \nu}\Big|_{-}d\sigma\\
\nm \ds &\q =\int_{\p D}  \Big(\la h\bvarphi, \ta\ra \la \nabla \F \ta, \n \ra +\frac{\lambda_0}{2\mu_0+\lambda_0}\la h\bvarphi, \n \ra \la \nabla \F \ta, \ta \ra\Big)d\sigma \\
\nm \ds &\q\q +\int_{\p D} \pd{\F}{\nu}\cdot\BDcal_{D}^{\sharp}[h\bvarphi]\big|_{+}d\sigma-\int_{\p D}\pd{\F}{\nu}\cdot\BDcal_{D}^{\sharp}[h\bvarphi]\big|_{-} d\sigma\\
\nm \ds &\q = \int_{\p D}  h \Big(\la \nabla \F \ta, \n\ra \ta +\frac{\lambda_0}{2\mu_0+\lambda_0}\la \nabla \F \ta,\ta\ra\n\Big)\cdot \bvarphi \, d\sigma \\
\nm \ds &\q\q +\int_{\p D} h \bigg [\pd{\BScal_{D}[\pd{\F}{\nu}]}{\n}\Big |_{-}-\pd{\BScal_{D}[\pd{\F}{\nu}]}{\n}\Big |_{+}\bigg]\cdot \bvarphi \,d\sigma.
\end{align*}
By using  \eqref{conormal-derivative}, \eqref{S-Sharp-n}, and  the identity $\nabla \cdot \F=\la \nabla \F \n , \n \ra+\la \nabla \F \ta, \ta\ra$,  we get
\begin{align*}
\ds \pd{\BScal_{D}[\pd{\F}{\nu}]}{\n}\Big |_{-}-\pd{\BScal_{D}[\pd{\F}{\nu}]}{\n}\Big |_{+}=-\nabla \F \n-\la \nabla \F \ta, \n\ra \ta -\frac{\lambda_0}{2\mu_0+\lambda_0}\la \nabla \F \ta,\ta\ra\n,
\end{align*}
and hence
\begin{align}\label{eqEMT-6}
\ds \int_{\p D}\pd{\v^e}{\nu}\cdot\BDcal_{D}^{\sharp}[h\bvarphi]\big|_{+} d\sigma-\int_{\p D}\v^e\cdot\frac{\p \BDcal_{D}^{\sharp}[h\bvarphi]}{\p  \nu}\Big|_{+}d\sigma=-\int_{\p D} h \pd{\F}{\n}\cdot \bvarphi \,d\sigma.
\end{align}
In conclusion, we obtain from \eqref{eq004} and \eqref{eqEMT-1}-\eqref{eqEMT-6} the theorem \ref{Third-theorem}.

In the  remaining part of this section we show that the asymptotic
expansion in \eqref{asymptotic-formula-EMTs} coincides with that one obtained in \cite[Theorem 3.1]{LY}. We can easily see from Proposition \ref{observation-Important} that \eqref{asymptotic-formula-EMTs} is equivalent to
\begin{align*}
 \ds \sum_{\alpha \beta j k}a_{j}^{\alpha} b_{k}^{\beta} m_{\alpha \beta k}^j(D_\ep)&=
 \sum_{\alpha \beta j k}a_{j}^{\alpha} b_{k}^{\beta} m_{\alpha \beta k}^j(D)\nonumber  \\
 \nm \ds & \q +\ep\int_{\p D}h\bigg(\big([\mathbb{M}_{1,0}-\mathbb{C}_0] \widehat{\nabla}\u^e\big)\ta\cdot \widehat{\nabla} \v^e \ta -(\mathbb{K}_{1,0} \widehat{\nabla}\u^e\big)\n\cdot (\mathbb{C}_0 \widehat{\nabla} \v^e )\n\bigg)d\sigma\\
 \nm\ds&\q +O(\ep^2),
\end{align*}
with
\begin{align}\label{Tensor-L}
\ds \q \mathbb{M}_{1,0}-\mathbb{C}_0&=\frac{2(\lambda_1\mu_0-\lambda_0 \mu_1)}{(\lambda_1+2\mu_1)}\I\otimes\I
+\frac{4(\mu_1-\mu_0)(\lambda_1+\mu_1)}{\lambda_1+2\mu_1}\I\otimes(\ta\otimes \ta)\nonumber\\
\nm \ds&:=\eta \I\otimes\I+\delta \I\otimes(\ta\otimes \ta),
\end{align}
and
\begin{align}\label{tensor-B}
\ds-\mathbb{K}_{1,0}&=\frac{(\lambda_1-\lambda_0)\mu_1-2(\mu_1-\mu_0)(\lambda_1+\mu_1)}{\mu_1(\lambda_1+2\mu_1)}\I\otimes \I+2\big(1-\frac{\mu_0}{\mu_1}\big) \mathbb{I}\nonumber\\
\nm\ds &\q +2\frac{(\mu_1-\mu_0)(\lambda_1+\mu_1)}{\mu_1(\lambda_1+2\mu_1)}\I\otimes (\ta\otimes \ta)\nonumber\\
\nm\ds &:=\rho\I\otimes \I+\tau\mathbb{I}+\varrho\I\otimes (\ta\otimes \ta).
\end{align}
Simple computations, yield
\begin{align*}
\ds \big([\mathbb{M}_{1,0}-\mathbb{C}_0] \widehat{\nabla}\u^e\big)\ta\cdot \widehat{\nabla} \v^e \ta&=\eta (\nabla \cdot \u^e)\la \widehat{\nabla} \v^e \ta, \ta\ra+\delta\la \widehat{\nabla} \u^e \ta, \ta\ra\la \widehat{\nabla} \v^e \ta, \ta\ra,\\
\nm \ds -\big(\mathbb{K}_{1,0} \widehat{\nabla}\u^e\big)\n\cdot \big(\mathbb{C}_0\widehat{\nabla} \v^e \big)\n &=\lambda_0 \rho (\nabla \cdot \u^e) (\nabla \cdot \v^e)
+2\mu_0 \rho \nabla \cdot \u^e \la \widehat{\nabla} \v^e \n, \n\ra\\
\nm\ds&\q +\lambda_0\tau (\nabla \cdot \v^e)\la \widehat{\nabla} \u^e \n, \n\ra+2\mu_0 \tau \la \widehat{\nabla} \u^e \n, \widehat{\nabla} \v^e \n\ra\\
\nm\ds &\q +\lambda_0\varrho (\nabla \cdot \v^e)\la \widehat{\nabla} \u^e \ta, \ta\ra+2\mu_0 \varrho \la \widehat{\nabla} \u^e \ta,\ta\ra\la \widehat{\nabla} \v^e \n, \n\ra.
\end{align*}
Note that
\begin{align*}
2\mu_0 \tau \la \widehat{\nabla} \u^e \n, \widehat{\nabla} \v^e \n\ra=\mu_0 \tau  \widehat{\nabla} \u^e:\widehat{\nabla} \v^e +\mu_0 \tau (\nabla \cdot \u^e)\la \widehat{\nabla} \v^e\n,\n\ra-\mu_0 \tau (\nabla \cdot \v^e)\la \widehat{\nabla} \u^e \ta,\ta\ra.
\end{align*}
Hence
\begin{align*}
\ds \big([\mathbb{M}_{1,0}-\mathbb{C}_0] \widehat{\nabla}\u^e\big)\ta\cdot \widehat{\nabla} \v^e \ta -(\mathbb{K}_{1,0} \widehat{\nabla}\u^e\big)\n\cdot (\mathbb{C}_0 \widehat{\nabla} \v^e )\n=\mathbb{S} \widehat{\nabla}\u^e:\widehat{\nabla} \v^e ,
\end{align*}
where
\begin{align}\label{S}
\mathbb{S} \widehat{\nabla}\u^e &=\lambda_0(\rho+\tau) (\nabla \cdot \u^e)\I+(\lambda_0\varrho-\lambda_0\tau+2\mu_0\varrho-\mu_0 \tau)\la \widehat{\nabla} \u^e \ta, \ta\ra \I+\eta (\nabla \cdot \u^e) \ta\otimes \ta\nonumber\\
\nm \ds &\q+(\delta-2\mu_0 \varrho)\la \widehat{\nabla} \u^e \ta, \ta\ra \ta\otimes \ta+(2\mu_0 \rho+\mu_0 \tau)(\nabla \cdot \u^e)\n\otimes \n+\mu_0\tau \widehat{\nabla}  \u^e.
\end{align}
It is proved in \cite{LY} that
\begin{align*}
 \ds \sum_{\alpha \beta j k}a_{j}^{\alpha} b_{k}^{\beta} m_{\alpha \beta k}^j(D_\ep)&=
 \sum_{\alpha \beta j k}a_{j}^{\alpha} b_{k}^{\beta} m_{\alpha \beta k}^j(D)+ \ep\int_{\p D}h \big(\mathbb{M}\widehat{\nabla} \u^e\big) : \widehat{\nabla} \v^e +O(\ep^{1+\gamma}),
\end{align*}
for some positive $\gamma$ and
\begin{align*}
\ds  \mathbb{M}\widehat{\nabla} \u^e:=\big(\mathbb{C}_1-\mathbb{C}_0\big)\mathbb{C}_1^{-1}\Big(\big(\mathbb{K}\widehat{\nabla} \u^e \ta\big)\otimes \ta+\big(\mathbb{C}_0\widehat{\nabla} \u^e \n\big)\otimes \n\Big),
\end{align*}
where the $4-$tensor $\mathbb{K}$ is defined by
\begin{align*}
 \ds \mathbb{K}:= p \I \otimes \I+2 \mu_0 \mathbb{I}+q \I\otimes (\ta\otimes \ta),
\end{align*}
where
\begin{align*}
 \ds p:= \frac{\lambda _1(\lambda_0+2\mu_0)}{\lambda_1+2\mu_1}\q\q  \mbox{and}\q \q q:= \frac{4(\mu_1-\mu_0)(\lambda_1+\mu_1)}{\lambda_1+2\mu_1}.
\end{align*}
Denote by
\begin{align*}
\lambda:=\frac{\lambda_1-\lambda_0+\mu_1-\mu_0}{2(\lambda_1+\mu_1)}-\frac{\mu_1-\mu_0}{2\mu_1}, \q\q \mu=\frac{\mu_1-\mu_0}{2\mu_1}.
\end{align*}
It is proved in  \cite{AEEKL} that
\begin{align}\label{M}
\ds  \mathbb{M}\widehat{\nabla} \u^e=&\big[\lambda (p+\lambda_0+2\mu_0)+2\mu p-\eta\big](\nabla \cdot \u^e)\I+\lambda q \la \widehat{\nabla} \u^e \ta, \ta\ra \I+\eta(\nabla \cdot \u^e)\ta\otimes \ta\nonumber\\
\nm \ds &+2 \mu q \la \widehat{\nabla} \u^e \ta, \ta\ra \ta\otimes \ta +\big[2 \mu \lambda_0+\eta-2\mu p\big] (\nabla \cdot \u^e)\n\otimes \n+4 \mu\mu_0 \widehat{\nabla} \u^e.
\end{align}
Looking at the coefficients in  \eqref{S} and \eqref{M}, we confirm that
\begin{align*}
\ds \mathbb{M}=\mathbb{S}&=\frac{\lambda_0(\lambda_1-\lambda_0)+2 \lambda_0(\mu_1-\mu_0)}{\lambda_1+2\mu_1}\I\otimes \I+2\big(1-\frac{\mu_0}{\mu_1}\big)\frac{\mu_0 \lambda_1-\mu_1\lambda_0}{\lambda_1+2\mu_1}\I\otimes(\ta\otimes \ta)\\
\nm \ds &\q +\frac{2(\lambda_1\mu_0-\lambda_0 \mu_1)}{(\lambda_1+2\mu_1)} (\ta\otimes \ta)\otimes\I+4\big(1-\frac{\mu_0}{\mu_1}\big)\frac{(\mu_1-\mu_0)(\lambda_1+\mu_1)}{\lambda_1+2\mu_1}(\ta\otimes \ta)\otimes(\ta\otimes \ta)\\
\nm\ds &\q +2\big(\frac{\mu_0}{\mu_1}\big)\frac{(\lambda_1-\lambda_0)\mu_1-(\mu_1-\mu_0)\lambda_1}{\lambda_1+2\mu_1}(\n\otimes \n)\otimes\I+\frac{2 \mu_0(\mu_1-\mu_0)}{\mu_1} \mathbb{I}.
\end{align*}
%%%%%%%%%%%%%%%%%%%%%%%%%%%%%%%%%%%%%%%%%
%%%%%%%%%%%%%%%%%%%%%%%%%%%%%%%%%%%%%%%%%%
%%%%%%%%%%%%%%%%%%%%%%%%%%%%%%%%%%%%%%%
%%%%%%%%%%%%%%%%%%%%%%%%%%%%%%%%%%%%%%%
%%%%%%%%%%%%%%%%%%%%%%%%%%%%%%%%%%%%%%%%%%
%%%%%%%%%%%%%%%%%%%%%%%%%%%%%%%%%%%%%
%%%%%%%%%%%%%%%%%%%%%%%%%%%%%%%%%%%%%%%%%
%%%%%%%%%%%%%%%%%%%%%%%%%%%%%%%%%%%%%%%%%%
%%%%%%%%%%%%%%%%%%%%%%%%%%%%%%%%%%%%%%%
%%%%%%%%%%%%%%%%%%%%%%%%%%%%%%%%%%%%%%%
%%%%%%%%%%%%%%%%%%%%%%%%%%%%%%%%%%%%%%%%%%
%%%%%%%%%%%%%%%%%%%%%%%%%%%%%%%%%%%%%
%%%%%%%%%%%%%%%%%%%%%%%%%%%%%%%%%%%%%%%%%
%%%%%%%%%%%%%%%%%%%%%%%%%%%%%%%%%%%%%%%%%%
%%%%%%%%%%%%%%%%%%%%%%%%%%%%%%%%%%%%%%%
%%%%%%%%%%%%%%%%%%%%%%%%%%%%%%%%%%%%%%%
%%%%%%%%%%%%%%%%%%%%%%%%%%%%%%%%%%%%%%%%%%
%%%%%%%%%%%%%%%%%%%%%%%%%%%%%%%%%%%%%
\newpage
\section*{Appendix}
\subsection*{1) Derivation of the $\displaystyle\frac{\p\BScal_D[\bphi]}{\p\nu}(x)$}
We have
\begin{align*}
\ds \BScal_D[\bphi](x)&=\int_{\p D}\G(x-y)\bphi(y)d\sigma(y)\\
\nm \ds &=\int_{\p D}\Big(\frac{A}{2\pi}\log|x-y|\bphi(y)-\frac{B}{2\pi}\frac{\la x-y,\bphi(y)\ra}{|x-y|^2}(x-y)\Big)d\sigma(y).
\end{align*}
Let $\x:=x-y$. By using \eqref{8}, we get
\begin{align*}
\ds &\nabla_x\Bigg(\frac{A}{2\pi}\log|\x|\bphi(y)-\frac{B}{2\pi}\frac{\la \x,\bphi(y)\ra}{|\x|^2}\x\Bigg)\\
\nm\ds &\q \q = \frac{A}{2\pi} \frac{\bphi(y)\otimes \x}{|\x|^2}-\frac{B}{2\pi}\frac{\la \x,\bphi(y)\ra }{|\x|^2}\I
+\frac{B}{\pi}\frac{\la \x,\bphi(y)\ra}{|\x|^4}\big(\x\otimes \x\big)-\frac{B}{2\pi}\frac{\x\otimes \bphi(y)}{|\x|^2}.
\end{align*}
Therefore
\begin{align}\label{nabla-nablaT-SD}
\ds &\nabla_x\Big(\frac{A}{2\pi}\log|\x|\bphi(y)-\frac{B}{2\pi}\frac{\la \x,\bphi(y)\ra}{|\x|^2}\x\Big)+\left[\nabla_x\Big(\frac{A}{2\pi}\log|\x|\bphi(y)-\frac{B}{2\pi}\frac{\la \x,\bphi(y)\ra}{|\x|^2}\x\Big)\right]^{T}\nonumber\\
\nm\ds &~~= \frac{(A-B)}{2\pi} \frac{\bphi(y)\otimes \x}{|\x|^2}+\frac{(A-B)}{2\pi} \frac{\x\otimes \bphi(y)}{|\x|^2}-\frac{B}{\pi}\frac{\la \x,\bphi(y)\ra }{|\x|^2}\I
+\frac{2B}{\pi}\frac{\la \x,\bphi(y)\ra}{|\x|^4}\big(\x\otimes\x\big),
\end{align}
which gives
\begin{align*}
\ds &\left(\nabla_x\Big(\frac{A}{2\pi}\log|\x|\bphi(y)-\frac{B}{2\pi}\frac{\la \x,\bphi(y)\ra}{|\x|^2}\x\Big)+\left[\nabla_x\Big(\frac{A}{2\pi}\log|\x|\bphi(y)-\frac{B}{2\pi}\frac{\la \x,\bphi(y)\ra}{|\x|^2}\x\Big)\right]^{T}\right)\n(x)\\
\nm\ds&= \frac{(A-B)}{2\pi} \frac{\la \x,\n(x)\ra}{|\x|^2}\bphi(y)+\frac{(A-B)}{2\pi} \frac{\la \n(x),\bphi(y)\ra }{|\x|^2}\x-\frac{B}{\pi}\frac{\la \x,\bphi(y)\ra }{|\x|^2}\n(x)\\
&\q+\frac{2B}{\pi}\frac{\la \x,\bphi(y)\ra\la \x,\n(x)\ra}{|\x|^4}\x\\
\nm\ds&= \frac{(A-B)}{2\pi} \frac{\la \x,\n(x)\ra}{|\x|^2}\bphi(y)+\frac{(A-B)}{2\pi} \frac{\x\otimes \n(x) }{|\x|^2}\bphi(y)-\frac{B}{\pi}\frac{\n(x)\otimes \x }{|\x|^2}\bphi(y)\\
&\q+\frac{2B}{\pi}\frac{\la \x,\n(x)\ra}{|\x|^4}(\x\otimes \x)\bphi(y)\\
\nm\ds&= \Big[\frac{(A-B)}{2\pi} \frac{\la \x,\n(x)\ra}{|\x|^2}\I+\frac{(A-B)}{2\pi} \frac{\x\otimes \n(x) }{|\x|^2}-\frac{B}{\pi}\frac{\n(x)\otimes \x }{|\x|^2}
+\frac{2B}{\pi}\frac{\la \x,\n(x)\ra}{|\x|^4}(\x\otimes \x)\Big]\bphi(y)\\
\nm\ds&:=\QQ(\x) \bphi(y).
\end{align*}
Hence, we obtain
\begin{align}\label{2-partie-SD}
\ds \Big(\nabla \BScal_D[\bphi](x)+\big[\nabla \BScal_D[\bphi](x)\big]^{T}\Big)\n(x)=\int_{\p D}\QQ(x-y) \bphi(y)d\sigma(y), \q x\in  \p D.
\end{align}
It follows from  \eqref{13} that
\begin{align}\label{Divergence-SD}
\ds \nabla_x \cdot \Big(\G(\x)\bphi(y)\Big)&=\la \nabla_x\cdot \G(\x), \bphi(y)\ra.
\end{align}
Thus
\begin{align*}
\ds \nabla_x \cdot \Big(\G(\x)\bphi(y)\Big)\n(x)&=\la \nabla_x\cdot \G(\x), \bphi(y)\ra \n(x)=\Big(\n(x)\otimes \nabla_x\cdot \G(\x)\Big) \bphi(y).
\end{align*}
Since
\begin{align}\label{Divergence-Gamma}
\ds \nabla_x \cdot \G(\x)&=\nabla_x\cdot  \Big(\frac{A}{2\pi}\log|\x|\I-\frac{B}{2\pi}\frac{\x\otimes \x}{|\x|^2}\Big)\nonumber\\
\nm\ds &=\frac{A}{2\pi}\frac{\x}{|\x|^2}-\frac{B}{2\pi}\Big(-2\frac{(\x\otimes \x)\x}{|\x|^4}+\frac{\nabla \x~ \x+\nabla\cdot\x ~\x}{|\x|^2}\Big)\nonumber\\
\nm\ds &=\frac{A-B}{2\pi}\frac{\x}{|\x|^2},
\end{align}
then
\begin{align*}
\ds \nabla_x \cdot \Big(\G(\x)\bphi(y)\Big)\n(x)=\frac{A-B}{2\pi}\frac{\n(x)\otimes \x}{|\x|^2}\bphi(y):=\PP(\x)\bphi(y),
\end{align*}
and hence
\begin{align}\label{1-partie-SD}
\ds \nabla  \cdot \BScal_D [\bphi](x) \n(x)=\int_{\p D}\PP(x-y)\bphi(y)d\sigma(y).
\end{align}
It then follows from  \eqref{conormal-derivative},  \eqref{2-partie-SD} and \eqref{1-partie-SD} that
\begin{align}\label{normal-derivative-SD}
\ds \frac{\p\BScal_D[\bphi]}{\p\nu}(x)=\int_{\p D}\big( \lambda_0 \PP(x-y)+ \mu_0 \QQ(x-y)\Big)\bphi(y)d\sigma(y).
\end{align}
Note that  $\lambda_0\PP(x-y)+\mu_0\QQ(x-y)=\KK^{T}(x-y)$ for $x,y \in \p D, x\neq y$,  where $\KK^{T}$ is defined by \eqref{hat-K}.

\subsection*{2) Derivation of  the   $\frac{\p \BDcal^{\sharp}_{D}[\bphi]}{\p\nu}(x)$}

According to \eqref{dcal-sharp-presentation}, we have
\begin{align*}
\ds& \nabla_x \Big(\pd{\G(\x)}{\n(y)}\bphi(y)\Big)\\
\nm\ds&\q=-\frac{A}{2\pi}\Big[\frac{\bphi(y)\otimes \n(y)}{|\x|^2}
-2 \frac{\la \x ,\n(y)\ra}{|\x|^4}\bphi(y)\otimes \x\Big]\\
\nm\ds &\q\q -\frac{B}{\pi}\Big[\frac{\la \x ,\bphi(y)\ra}{|\x|^4}\x\otimes \n(y)+\frac{\la \x ,\n(y)\ra}{|\x|^4} \x\otimes\bphi(y)-4\frac{\la \x ,\n(y)\ra \la \x, \bphi(y)\ra}{|\x|^6}\x\otimes \x\\
\nm\ds &~~~~~~~~~~~~~~ ~~~~~~~~~~~~~~~~~~~~~~~~~~~~~~~~~~~~~~~~~~~~~~~~~~~~~~~~~~~~~~~~~~~~+\frac{\la \x ,\n(y)\ra\la \x, \bphi(y)\ra}{|\x|^4}\I\Big]\\
\nm\ds &\q\q+\frac{B}{2\pi}\Big[\frac{\la \n(y) ,\bphi(y)\ra}{|\x|^2}\I-2\frac{\la \n(y) ,\bphi(y)\ra}{|\x|^4}\x\otimes \x+\frac{\n(y)\otimes \bphi(y)}{|\x|^2}-2\frac{\la \x ,\bphi(y)\ra}{|\x|^4}\n(y)\otimes \x\Big].
\end{align*}
Then
\begin{align}\label{sharp1}
\ds& \nabla_x \Big(\pd{\G(\x)}{\n(y)}\bphi(y)\Big)\n(x)\nonumber\\
\nm\ds&\q=-\frac{A}{2\pi}\Big[\frac{\la \n(x), \n(y)\ra}{|\x|^2}\I
-2 \frac{\la \x ,\n(y)\ra\la \x ,\n(x)\ra}{|\x|^4}\I\Big]\bphi(y)\nonumber\\
\nm\ds &\q\q -\frac{B}{\pi}\Big[\frac{\la \n(x) ,\n(y)\ra}{|\x|^4}\x\otimes \x+\frac{\la \x ,\n(y)\ra}{|\x|^4} \x\otimes\n(x)-4\frac{\la \x ,\n(y)\ra \la \x, \n(x)\ra}{|\x|^6}\x\otimes \x\\
\nm\ds &~~~~~~~~~~~~~~ ~~~~~~~~~~~~~~~~~~~~~\q~~~~~~~~~~~~~~~~~~~~~~~~~~~~~~~~~~~~~~~+\frac{\la \x ,\n(y)\ra}{|\x|^4}\n(x)\otimes \x\Big]\bphi(y)\nonumber\\
\nm\ds &\q\q+\frac{B}{2\pi}\Big[\frac{\n(x)\otimes \n(y) }{|\x|^2}-2\frac{\la \x ,\n(x)\ra}{|\x|^4}\x\otimes \n(y)+\frac{\n(y)\otimes \n(x)}{|\x|^2}-2\frac{\la \x ,\n(x)\ra}{|\x|^4}\n(y)\otimes \x\Big]\bphi(y).\nonumber
\end{align}
Likewise, we get
\begin{align}\label{sharp2}
\ds& \Big[\nabla_x \Big(\pd{\G(\x)}{\n(y)}\bphi(y)\Big)\Big]^{T}\n(x)\nonumber\\
\nm\ds&\q=-\frac{A}{2\pi}\Big[\frac{ \n(y)\otimes \n(x)}{|\x|^2}
-2 \frac{\la \x ,\n(y)\ra}{|\x|^4} \x \otimes \n(x)\Big]\bphi(y)\nonumber\\
\nm\ds &\q\q -\frac{B}{\pi}\Big[\frac{\la \x ,\n(y)\ra}{|\x|^4}\n(x)\otimes \x+\frac{\la \x ,\n(x)\ra}{|\x|^4} \n(y)\otimes\x-4\frac{\la \x ,\n(y)\ra \la \x, \n(x)\ra}{|\x|^6}\x\otimes \x\\
\nm\ds &~~~~~~~~~~~~~~ ~~~~~~~~~~~~~~~~~~~~~\q~~~~~~~~~~~~~~~~~~~~~~~~~~~~~~~~~~~~~~~+\frac{\la \x ,\n(y)\ra\la \x,\n(x)\ra}{|\x|^4}\I\Big]\bphi(y)\nonumber\\
\nm\ds &\q\q+\frac{B}{2\pi}\Big[\frac{\n(x)\otimes \n(y) }{|\x|^2}-2\frac{\la \x ,\n(x)\ra}{|\x|^4}\x\otimes \n(y)+\frac{\la \n(x),\n(y)\ra}{|\x|^2}\I-2\frac{\la \n(x),\n(y)\ra}{|\x|^4}\x \otimes \x\Big]\bphi(y).\nonumber
\end{align}
Using \eqref{11}, \eqref{12},   \eqref{13}, and \eqref{dcal-sharp-presentation}, we readily get
\begin{align}\label{sharp3}
\ds \nabla_x \cdot \Big(\pd{\G(\x)}{\n(y)}\bphi(y)\Big)\n(x)=&\Big\la \nabla_x \cdot\pd{\G(\x)}{\n(y)},\bphi(y)\Big\ra \n(x)\nonumber\\
\nm\ds =&\frac{A-B}{2\pi}\Big[2 \frac{\la \x , \n(y)\ra}{|x|^4}\n(x)\otimes \x -\frac{\n(x)\otimes \n(y)}{|x|^2}\Big]\bphi(y).
\end{align}
It then follows from \eqref{conormal-derivative}, \eqref{dcal-sharp}, \eqref{sharp1}, \eqref{sharp2}, and \eqref{sharp3} that
\begin{align}\label{dcal-sharp-integral}
\ds &\frac{\p\BDcal^{\sharp}_D[\bphi]}{\p\nu} (x)\nonumber\\
\nm\ds &=\frac{1}{2\pi}\frac{A-B}{A+B}\int_{\p D}
\Big[2 \frac{\la \x, \n(y)\ra\la \x, \n(x)\ra}{|\x|^4}-\frac{\la \n(x),\n(y)\ra}{|\x|^2} \Big]\bphi(y) d\sigma(y)\nonumber\\
\nm \ds &\q+\frac{1}{2\pi}\frac{A-B}{A+B}\int_{\p D}
\Big[2 \frac{\la \x, \n(y)\ra}{|\x|^4}\Big(\x \otimes \n(x)-\n(x)\otimes \x\Big)\nonumber\\
\nm\ds &~~~~~~~~~~~~~~~~~~~~~~~~~~~~~~~~~~~~~~~~~~~~~~-\frac{\n(y) \otimes \n(x)-\n(x)\otimes \n(y)}{|\x|^2}\Big]\bphi(y) d\sigma(y)\nonumber \\
\nm \ds &\q+\frac{1}{\pi}\frac{2B}{A+B}\int_{\p D}
\Big[4 \frac{\la \x, \n(y)\ra\la \x, \n(x)\ra}{|\x|^6}\x\otimes \x-\frac{\la \n(x),\n(y)\ra}{|\x|^4}\x\otimes \x\\
\nm\ds~&~~~~~~~~~~~~~~~~~\q\q~~~~~~~~~~~~~~~~~~-\frac{\la \x,\n(x)\ra}{|\x|^4} \Big(x\otimes \n(y)+\n(y)\otimes \x\Big)\Big]\bphi(y) d\sigma(y)\nonumber.
\end{align}

\subsection*{3) Derivation of the $\nabla \nabla\cdot \BScal_D[\bphi](x)\cdot \n(x) \n(x)$}

It follows from \eqref{Divergence-SD} and \eqref{Divergence-Gamma} that
$$
\nabla \cdot \BScal_{D}[\bphi](x)=\frac{(A-B)}{2\pi}\int_{\p D}\frac{\la \x, \bphi(y)\ra}{|\x|^2} d\sigma(y).
$$
Thus
\begin{align*}
\ds \nabla \nabla \cdot \BScal_{D}[\bphi](x)=\frac{(A-B)}{2\pi}\int_{\p D}\Big(\frac{ \bphi(y)}{|\x|^2}-2\frac{\la \x, \bphi(y)\ra}{|\x|^4}\x\Big) d\sigma(y),
\end{align*}
which gives
\begin{align*}
\ds \nabla \nabla \cdot \BScal_{D}[\bphi](x)\cdot \n(x)=\frac{(A-B)}{2\pi}\int_{\p D}\Big(\frac{\la \bphi(y),\n(x)\ra}{|\x|^2}-2\frac{\la \x, \bphi(y)\ra\la \x, \n(x)\ra}{|\x|^4}\Big) d\sigma(y).
\end{align*}
Thanks to the identity in  \eqref{2}, we obtain
\begin{align}\label{Nabla-Divergence-SD-n-n}
\ds \nabla \nabla \cdot \BScal_{D}[\bphi](x)\cdot \n(x)\n(x)=\frac{(A-B)}{2\pi}\int_{\p D}\Big(\frac{ \n(x)\otimes \n(x)}{|\x|^2}-2\frac{\la \x, \n(x)\ra}{|\x|^4}\n(x)\otimes \x\Big) \bphi(y) d\sigma(y).
\end{align}

\subsection*{4)  Derivation of  the $\nabla \Big(\nabla \BScal_D [\bphi](x)+\big[\nabla \BScal_D[\bphi](x)\big]^{T}\Big)\n(x)\n(x)$}
It follows from \eqref{2} and \eqref{nabla-nablaT-SD} that
\begin{align*}
\ds &\nabla\Bigg(\nabla\Big(\frac{A}{2\pi}\log|\x|\bphi(y)-\frac{B}{2\pi}\frac{\la \x,\bphi(y)\ra}{|\x|^2}\x\Big)+\left[\nabla\Big(\frac{A}{2\pi}\log|\x|\bphi(y)-\frac{B}{2\pi}\frac{\la \x,\bphi(y)\ra}{|\x|^2}\x\Big)\right]^{T}\Bigg)\n(x)\n(x)\nonumber\\
\nm\ds &= \Bigg[\frac{(A-B)}{2\pi}\Bigg( \frac{\bphi(y)\otimes \I+\big(\I\otimes\bphi(y)\big)^{T}}{|\x|^2}-2 \frac{\bphi(y)\otimes\x\otimes \x+\x\otimes \bphi(y)\otimes \x}{|\x|^4}\Bigg)\nonumber\\
\nm\ds &\q-\frac{B}{\pi}\Bigg(\frac{\I\otimes \bphi(y)}{|\x|^2}-2\frac{\la\x,\bphi(y) \ra}{|\x|^4} \big(\I \otimes \x\big)\Bigg)\nonumber\\
\nm\ds &\q+\frac{2B}{\pi} \Bigg(\frac{\x\otimes\x\otimes \bphi(y)}{|\x|^4}-4 \frac{\la\x,\bphi(y)\ra}{|\x|^6} \Big(\x\otimes \x\otimes \x \Big)+\frac{\la \x,\bphi(y)\ra}{|\x|^4}\Big((\I\otimes \x)^{T}+\x\otimes \I\Big)\Bigg)\Bigg]\n(x)\n(x)\\
\nm\ds &= \Bigg[\frac{(A-B)}{2\pi}\Bigg( \frac{\bphi(y)\otimes \n(x)+\n(x)\otimes\bphi(y)}{|\x|^2}-2 \frac{\la \x,\n(x) \ra}{|\x|^4} \Big(\bphi(y)\otimes\x+\x\otimes \bphi(y)\Big)\Bigg)\nonumber\\
\nm\ds &\q-\frac{B}{\pi}\Bigg(\frac{\la \bphi(y),\n(x)\ra }{|\x|^2}-2\frac{\la\x,\bphi(y) \ra}{|\x^2|}\frac{\la \x,\n(x)\ra }{|\x|^2}\Bigg)\I\nonumber\\
\nm\ds &\q+\frac{2B}{\pi} \Bigg(\frac{\la \n(x),\bphi(y)\ra}{|\x|^4}\big(\x\otimes\x\big)-4 \frac{\la\x,\bphi(y)\ra\la \x,\n(x)\ra}{|\x|^6}\big(\x\otimes\x\big)\\
\nm\ds &\q\q\q\q\q\q\q\q\q\q\q\q\q\q\q\q\q\q\q\q\q +\frac{\la \x,\bphi(y)\ra}{|\x|^4}\Big(\n(x)\otimes \x+\x\otimes \n(x)\Big)\Bigg)\Bigg]\n(x)\\
\nm\ds &= \Bigg[\frac{(A-B)}{2\pi}\Bigg( \frac{ \I+\n(x)\otimes\n(x)}{|\x|^2}-2\frac{(\la \x,\n(x) \ra)^2}{|\x|^4}\I-2\frac{\la \x,\n(x) \ra}{|\x|^4} \Big(\x\otimes\n(x)\Big)\Bigg)\nonumber\\
\nm\ds &\q\q -\frac{B}{\pi}\Bigg(\frac{\n(x) \otimes \n(x) }{|\x|^2}
-2\frac{\la\x,\n(x) \ra}{|\x|^4}\Big(\n(x)\otimes \x \Big)\Bigg)\nonumber\\
\nm\ds &\q\q +\frac{2B}{\pi} \Bigg(\frac{\la\x, \n(x)\ra}{|\x|^4}\Big(\x\otimes\n(x)+\n(x)\otimes \x\Big)-4 \frac{(\la\x,\n(x)\ra)^2}{|\x|^6}\big(\x\otimes\x\big) +\frac{\x\otimes \x }{|\x|^4}\Bigg)\Bigg]\bphi(y)\\
\nm\ds &:=\LL(\x)\bphi(y).
\end{align*}
Then, we have
\begin{align*}
\ds \nabla \Big(\nabla \BScal_D[\bphi](x)+\big[\nabla \BScal_D[\bphi](x)\big]^{T}\Big)\n(x)\n(x)=\int_{\p D}\LL(x-y)\bphi(y)d\sigma(y).
\end{align*}


\begin{thebibliography}{10}

\bibitem{AEEKL} H. Ammari, E. Beretta, E. Francini, H. Kang, and M. Lim, Reconstruction of small
interface changes of an inclusion from modal measurements II: the elastic case J. Math.
Pures et Appl., 94 (2010),  3, 322--339.

\bibitem{ABGKLW} H. Ammari, E. Bretin, J. Garnier, H. Kang, H. Lee, and A. Wahab,  \textit{Mathematical methods in elasticity imaging}. Princeton Series in Applied Mathematics. Princeton University Press, Princeton, NJ, 2015.

\bibitem{AGKLS} H. Ammari, J. Garnier, H. Kang, M. Lim, and K. Solna,
Multistatic imaging of extended targets. SIAM Journal on Imaging Sciences, 5 (2012), 564--600.

\bibitem{book} H. Ammari and H. Kang, \textit{Reconstruction of small inhomogeneities from boundary
measurements, Lecture Notes in Mathematics}, Vol. 1846, Springer-Verlag, Berlin, 2004.


\bibitem{book2} -----------------,\textit{ Polarization and moment
tensors with applications to inverse problems and effective medium theory},
 Applied Mathematical Sciences, Vol. 162, Springer-Verlag, New York, 2007.

\bibitem{AKLZ1} H. Ammari, H. Kang, M. Lim, and H. Zribi, Conductivity interface problems. Part I:
Small perturbations of an interface, Trans. Amer. Math. Soc., 362 (2010), 2435--2449.




\bibitem{AKLZ2}  H. Ammari, H. Kang, M. Lim, and H. Zribi, The generalized polarization tensors for
resolved imaging. Part I: shape reconstruction of a conductivity inclusion, Math. of
Comp., 81 (2012), 367--386.


\bibitem{AKNT} H. Ammari, H. Kang, G. Nakamura, and K. Tanuma, Complete asymptotic expansions
of solutions of the system of elastostatics in the presence of an inclusion of Small
diameter and detection of an inclusion, Jour. of Elasticity, 67 (2002), 97--129.

\bibitem{BBFM}  E. Beretta, E. Bonnetier, E. Francini, and A. L. Mazzucato, An asymptotic formula for
the displacement field in the presence of small anisotropic elastic inclusions, Inverse Problems and Imaging, 6
(2012), 1--23.

\bibitem{BF}  E. Beretta and E. Francini, An asymptotic formula for the displacement field in the
presence of thin elastic inhomogeneities, SIAM J. Math. Anal., 38, (2006), 1249--1261.


\bibitem{CGHIR}  R. R. Coifman, M. Goldberg, T. Hrycak, M. Israeli, and V. Rokhlin, An improved
operator expansion algorithm for direct and inverse scattering computations, Waves
Random Media, 9 (1999), 441--457.


\bibitem{CMM82}  R. R. Coifman, A. McIntosh, and Y. Meyer, L'int\'egrale de Cauchy d\'efinit un op\'erateur
bourn\'ee sur $L^2$ pour les courbes lipschitziennes, Ann. Math., 116 (1982), 361--387.


\bibitem{DKV} B. E. Dahlberg, C.E. Kenig, and G. Verchota, Boundary value problem for the systems
of elastostatics in Lipschitz domains, Duke Math. Jour., 57 (1988), 795--818.

\bibitem{ES}  L. Escauriaza and J. K. Seo, Regularity properties of solutions to transmission problems,
Trans. Amer. Math. Soc., 338 (1993), 405--430.

\bibitem{FJR} E. B. Fabes, M. Jodeit, and N. M. Riviere, Potential techniques for boundary value problems
on $\mathcal{C}^1$ domains, Acta Math., 141 (1978), 165--186.

\bibitem{Folland76}  G. B. Folland, \textit{Introduction to partial differential equations}, Princeton University Press,
Princeton, New Jersey, 1976.


\bibitem{KKL} H. Kang, E. Kim, and J.-Y. Lee, identification of elastic inclusions and elastic moment
tensors by boundary measurements, Inverse Problems, 19 (2003), 703-724.


\bibitem{Kato} T. Kato, \textit{Perturbation theory for linear operators}, Springer-Verlag, New York, 1976.


\bibitem{KZ1}  A. Khelifi and H. Zribi, Asymptotic expansions for the voltage potentials with two and
three-dimensional thin interfaces, Math. Methods Appl. Sci., 34, (2011), 2274--2290.


\bibitem{KZ2} A. Khelifi and H. Zribi, Boundary voltage perturbations resulting from small surface
changes of a conductivity inclusion, Appl. Anal., 93, (2014), 46--64.


\bibitem{LLZ} M. Lim, K. Louati, and H. Zribi, An asymptotic formalism for reconstructing small
perturbations of scatterers from electric or acoustic far-field measurements, Math. Methods
Appl. Sci., 31 (2008), 1315--1332.


\bibitem{LY}  M. Lim and S. Yu, Reconstruction of the shape of an inclusion from elastic moment tensors,
 Contemp. Math., 548, (2011), 61--76.

\bibitem{LN} Y. Y. Li and L. Nirenberg, Estimates for elliptic systems from composite material, Comm.
Pure Appl. Math., 56 (2003),  892--925.


\bibitem{Zribi1} H. Zribi, Asymptotic expansions for currents caused by small interface changes of
an electromagnetic inclusion, Appl. Anal., 92 (2013), 172--190.


\end{thebibliography}
\end{document}